\documentclass[preprint, 11pt]{elsarticle}
 \usepackage[margin=2.5cm]{geometry}

\usepackage{lineno,hyperref}
\modulolinenumbers[5]

\journal{Computer Methods in Applied Mechanics and Engineering}

\usepackage{amsmath}
\usepackage{amsfonts}
\usepackage{amssymb}
\usepackage{bm}
\usepackage{mathrsfs}
\usepackage{physics}
\usepackage{upgreek}
\usepackage{mathtools}

\usepackage[labelformat=simple]{subcaption}

\usepackage{graphicx}
\graphicspath{ {figures/} }

\usepackage{booktabs}
\usepackage{multirow}
\usepackage{caption}
\usepackage{xcolor}
\usepackage{tablefootnote}

\usepackage{float}
\makeatletter
\newcommand\fs@boxedtop
 {\fs@boxed
  \def\@fs@mid{\vspace\abovecaptionskip\relax}%
  \let\@fs@iftopcapt\iftrue
 }
\makeatother
\floatstyle{boxedtop}
\floatname{framedbox}{Box}
\newfloat{framedbox}{hbt}{lob}[section]

\usepackage[textsize=tiny]{todonotes}

\begin{document}

\begin{frontmatter}

\title{Adaptivity for clustering-based reduced-order modeling of localized history-dependent phenomena}


\author[mymainaddress,mysecaddress]{Bernardo P. Ferreira}
\author[mymainaddress]{F.M. Andrade Pires\corref{mycor2}}
\ead{fpires@fe.up.pt}
\author[mysecaddress]{M.A. Bessa}
\ead{M.A.Bessa@tudelft.nl}

\address[mymainaddress]{DEMec - Department of Mechanical Engineering, Faculty of Engineering, University of Porto, \\ Rua Dr. Roberto Frias, 4200-465 Porto, Portugal}
\address[mysecaddress]{3mE - Faculty of Mechanical, Maritime and Materials Engineering, Delft University of Technology, \\ Mekelweg 2, 2628 CD Delft, The Netherlands}

\cortext[mycor2]{Corresponding author}

\begin{abstract}
This paper proposes a novel Adaptive Clustering-based Reduced-Order Modeling (ACROM) framework to significantly improve and extend the recent family of clustering-based reduced-order models (CROMs). 
This adaptive framework enables the clustering-based domain decomposition to evolve dynamically throughout the problem solution, ensuring optimum refinement in regions where the relevant fields present steeper gradients. It offers a new route to fast and accurate material modeling of history-dependent nonlinear problems involving highly localized plasticity and damage phenomena. The overall approach is composed of three main building blocks: target clusters selection criterion, adaptive cluster analysis, and computation of cluster interaction tensors. In addition, an adaptive clustering solution rewinding procedure and a dynamic adaptivity split factor strategy are suggested to further enhance the adaptive process. 
The coined Adaptive Self-Consistent Clustering Analysis (ASCA) is shown to perform better than its static counterpart when capturing the multi-scale elasto-plastic behavior of a particle-matrix composite and predicting the associated fracture and toughness. Given the encouraging results shown in this paper, the ACROM framework sets the stage and opens new avenues to explore adaptivity in the context of CROMs.
\end{abstract}

\end{frontmatter}



\section{Introduction \label{sec:introduction}}

The development of new materials with tailored properties and improved functionality has helped to expand the limits of human endeavor and foster achievements in a broad spectrum of scientific areas and industries (e.g., aerospace, aeronautics, electronics). To meet high-performance criteria, reducing development costs, and fulfill safety standards, materials and structures can be designed nowadays within the paradigm called Integrated Computational Materials Engineering (ICME) (\cite{horstemeyer2012, horstemeyer2018}). Fundamental to such a computational framework is the ability to perform accurate and efficient predictions of material behavior. However, fast and accurate material modeling involving highly localized plasticity or damage phenomena remains an elusive target. Addressing this critical challenge opens new possibilities for the design of materials and structures, integrating (horizontal) process-structure-property-performance relationships with (vertical) multi-scale process-structure and structure-property modeling over different scales. This article proposes a new route based on adaptive clustering that applies to any Clustering-based Reduced-Order Model (CROM) and that significantly enhances their ability to tackle such challenges -- this framework is coined Adaptive Clustering-based Reduced-Order Modeling (ACROM).

Despite the high accuracy of the direct numerical simulations (DNS) and the quantitative completeness of several multi-scale models (e.g., \cite{smit1998,feyel1998,miehe1999,terada2001,kochmann2016}), the associated computational costs may (1) call for simplified material models and suboptimal spatial discretizations, (2) become prohibitive when considering multiple scales simultaneously and (3) prevent building extensive and reliable material response datasets required by data-driven methodologies (e.g., \cite{yvonnet2009,clement2012,le2015, bessa2017}). Given the current computational resources, surpassing these limitations requires resorting to reduced-order models (ROMs) capable of achieving a key balance between accuracy and efficiency. Clustering-based ROMs (CROMs) (e.g., \cite{liu2016,wulfinghoff2017}) are a recent family of ROMs that have been particularly successful in modeling material nonlinear elasto-plastic behavior and predicted irreversible phenomena only from elastic deformations of the material in a prior learning stage. This is a distinctive characteristic when compared with other significant contributions in the literature such as the Transformation Field Analysis \cite{dvorak1992}, Nonuniform Transformation Field Analysis \cite{michel2003}, Proper Generalized Decomposition \cite{ladeveze2010}, Reduced Basis Method \cite{boyaval2008}, High-Performance Reduced-Order Model \cite{hernandez2014}, Empirical Cubature Method \cite{hernandez2017} and Wavelet-Reduced-Order Model \cite{vantuijl2019}.

The underlying idea of CROMs is to perform model reduction through a clustering-based domain decomposition relying on unsupervised machine learning (clustering algorithms). This novel approach sprouted from the work of Liu and coworkers \cite{liu2016} when they developed the Self-Consistent Clustering Analysis (SCA) method. Soon after, Wulfinghoff and coworkers \cite{wulfinghoff2017} proposed a CROM derived from the Hashin-Shtrikman variational principle that turned out to be equivalent to the SCA formulation \cite{schneider2019}. Since then, CROMs gained traction (e.g., \cite{liu2018,tang2018,yu2019,nie2019,yu2021}), mathematical foundations have been consolidated (e.g., \cite{schneider2019,li2019}) and several numerical applications have been successfully explored (e.g., \cite{yan2018,shakoor2019,cavaliere2019,he2020,gao2020,han2020}).

CROMs are two-stage algorithms. In the offline-stage, the spatial domain is decomposed such that material points with similar mechanical behavior are grouped in material clusters via a clustering algorithm. This clustering-based model reduction precedes the solution of the actual equilibrium problem occurring at the online-stage based on the so-called Lippman-Schwinger equation discretized by the material clusters. The cluster analysis is a crucial step in the process and is usually conducted on a small set of representative loadings and assuming an elastic constitutive behavior. So far, this spatial clustering discretization has remained static, i.e., it was conducted once at the offline stage and remained constant irrespective of the conditions found in the actual problem solved during the online-stage. These may include, for instance, general non-monotonic loading paths, history-dependent nonlinear constitutive behavior, and localization or damage phenomena. However, as discussed in this article, there are limits to what can be predicted without adaptively clustering the material domain of CROMs due to the onset and propagation of highly localized phenomena such as material yielding and fracture.

\begin{figure}[hbt]
\centering
		\begin{subfigure}[b]{0.49\textwidth}
                \centering
                \includegraphics[width=\textwidth]{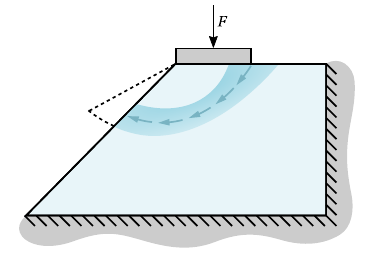}
                \caption{}
                \label{subfig:strain_localization_bench_1}
        \end{subfigure}
        \begin{subfigure}[b]{0.49\textwidth}
                \centering
                \includegraphics[width=\textwidth]{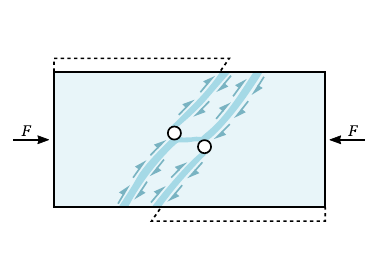}
                \caption{}
                \label{subfig:strain_localization_bench_2}
        \end{subfigure}\hfill
\caption{Schematic illustration of problems involving strain softening and localization: \subref{subfig:strain_localization_bench_1} Rigid footing placed on an elasto-plastic embankment; \subref{subfig:strain_localization_bench_2} Uniaxial compression of plane strain elasto-plastic specimen with two circular openings. Arrows within the material denote \subref{subfig:strain_localization_bench_1} plastic flow and \subref{subfig:strain_localization_bench_2} plastic slip.}
\label{fig:strain_localization_bench}
\end{figure}

The Adaptive Clustering-based Reduced-Order Modeling (ACROM) framework introduced in this article shares common features with Adaptive Finite Element Methods (AFEMs) pioneered by Babu\v{s}ka and Rheinboldt \cite{babuska1978,babuska1979}. Adaptive procedures have become essential in practical engineering analyses based on the Finite Element Method (e.g., \cite{ladeveze1998,ainsworth2000,stein2005,zienkiewicz2013}). In the context of elasto-plastic material behavior, they deal with the loss of ellipticity of the boundary value equilibrium problem as well as the consequent nonuniqueness of the solution (\cite{huerta1999,zienkiewicz2014}). Strain softening and localization (\cite{rice1976,pietruszczak1981,ortiz1991}) are among the most relevant reasons for this harmful behavior. Localization of deformation refers to the emergence of narrow regions or bands in a structure where all further deformation tends to concentrate, despite the external loading following a monotonic path \cite{deborst1993a} (see Figure~\ref{fig:strain_localization_bench}). Once localization occurs, large strains may accumulate inside the band without substantially affecting the strains in the surrounding material (\cite{hutchinson1980,ortiz1987}), the latter usually unloading elastically and behaving in a quasi-rigid manner. The localization phenomena are thus associated with displacement (e.g., sliding), strain and stress discontinuities \cite{zienkiewicz2014}. This phenomenon is observed for a wide range of materials, although the scale of localization, often associated with a `band width', may differ by some orders of magnitude. Given that it has a detrimental effect on the integrity of the structure, it often acts as a direct precursor to structural failure and fracture mechanisms \cite{deborst1993a}. To avoid the spurious dependence of the solution on the mesh refinement, the suitable handling of this phenomenon calls for a very fine mesh in the localization area that is, in general, not known a priori \cite{huerta1998}. Different strategies have been proposed to solve problems involving strain localization, such as methods involving regularization (e.g., \cite{zienkiewicz1974,deborst1993a}), non-local formulations (e.g., \cite{bazant1988a,reis2018}),  gradient-enhancement (e.g., \cite{deborst1995,geers1998,kuhl2000}),  phase-field fracture (e.g., \cite{francfort1998,bourdin2008,miehe2010,shahba2019}), concentrated discontinuities (e.g., \cite{oliver1993,oliver1999}),  the extended finite element method (e.g., \cite{moes1999,chessa2002,moes2002}) and cohesive zone models (e.g., \cite{hillerborg1976,needleman1987,tvergaard1990,paggi2017}). These and other strategies can be effectively coupled with adaptive finite element methodologies (e.g., \cite{ortiz1987,pastor1991,steinmann1992,belytschko1993,deborst1993a,peric1994,zienkiewicz1995,zienkiewicz1995a}).

AFEMs have three main ingredients \cite{huerta1998}: (1) an error estimator or indicator, employed to locate where there is a need for mesh refinement/de-refinement; (2) a procedure to adapt the spatial interpolation, increasing or decreasing the interpolation in a particular region of the computational domain; and (3) a remeshing criterion, translating the output of error analysis into a need for adaptivity and actual mesh parameters (e.g., minimum element size). A brief description of these ingredients is given in \ref{app:AFEMs} to provide a suitable background for this article. The interested reader is also referred to the extensive literature on the topic (e.g., \cite{ladeveze1998,ainsworth2000,stein2005,zienkiewicz2013}).

The authors present for the first time the ACROM framework and illustrate its application with the Adaptive Self-consistent Clustering Analyses (ASCA), leveraging relevant contributions from AFEMs and introducing new clustering adaptivity procedures to CROMs. The proposed solution avoids having the clustering-based domain decomposition dictated solely by the offline (training) stage, improving the solution accuracy by dynamically adapting the clustering-based domain decomposition in accordance with the evolution of the actual equilibrium problem. Particular focus is given to capture localization phenomena, often the precursor of material failure and fracture.

The characteristics and challenges of adaptivity in CROMs are discussed in Section~\ref{sec:methodology}, followed by the complete description of the proposed ACROM framework. Several numerical results are then shown in Section~\ref{sec:results}, aiming to demonstrate the performance of the novel methodology and discuss the associated limitations. At last, the main conclusions are summarized in Section~\ref{sec:conclusion}, where future research challenges are also addressed.

\section{Methodology \label{sec:methodology}}

\subsection{Scope of application}
The ACROM framework proposed in the present paper is applicable to most CROMs proposed in the literature (see Section~\ref{sec:introduction}). Despite requiring some intrusive operations related to the clustering-based spatial domain decomposition, the majority of the clustering adaptivity procedures are independent of the particular equilibrium problem being solved. Concerning the paradigm of multi-scale modeling, this framework can be employed at any scale whose equilibrium problem is solved through any CROM. Given its extensive application in the literature, a standard macroscale-microscale hierarchical multi-scale model is considered in this paper, focusing on the adaptive solution of the underlying microscale equilibrium problem.

The material microstructure can be defined by any number of material phases, each characterized by any class of constitutive model under infinitesimal or finite strains. In addition, the proposed framework is essentially independent of the number of dimensions of the problem, i.e., if either two-dimensional plane stress/plane strain ($n_{\mathrm{dim}}=2$) or three-dimensional ($n_{\mathrm{dim}}=3$) models are assumed.

Although the presentation of the proposed framework is kept as general as possible, the already introduced Self-Consistent Clustering Analysis (SCA) method \cite{liu2016} is adopted as the CROM of interest throughout this paper. Specifics about this particular CROM are clearly remarked amidst the framework formulation. Although not strictly necessary, the reader unfamiliar with this method is encouraged to read \ref{app:sca} for a concise description of the underlying formulation. A more complete and comprehensive description can be found in \cite{ferreira2022}.

\subsection{Nomenclature and fundamental concepts \label{ssec:nomenclature}}
Before describing the proposed clustering adaptivity framework, some fundamental concepts and nomenclature employed throughout the remainder of this paper are introduced here.

At the macroscale, the domain and its boundary are denoted by $\Omega$ and $\partial \Omega$, and the coordinates of a material point are given by $\bm{X}$ and $\bm{x}$ in the reference and deformed configurations, respectively. At the microscale, the material microstructure is assumed to be properly characterized by a representative volume element (RVE) of dimensions $(l_{\mathrm{RVE}})_{i}, \, i= 1,\dots, n_\mathrm{dim}$ and whose domain and boundary are denoted by $\Omega_{\mu}$ and $\partial \Omega_{\mu}$ (see Figure~\ref{fig:adaptivity_nomenclature}). At this scale, the coordinates of a material point are given by $\bm{Y}$ and $\bm{y}$ in the reference and deformed configurations, respectively. In addition, the reference configuration is generally denoted as $(\bullet)_{0}$ and the microscale fields are denoted as $(\bullet)_{\mu}$.

The RVE is here assumed to be discretized in a regular grid of voxels, $n_{v} = \prod_{i=1}^{n_{\mathrm{dim}}} n_{i}$, where $n_{i}$ denotes the number of voxels in the $i$th dimension (see Figure~\ref{fig:adaptivity_nomenclature}). This type of discretization is suitable for CROMs whose solution procedure is partially computed in the discrete frequency domain, but can be easily converted to a finite element mesh of quadrilateral/hexahedral finite elements (see \ref{app:femrgmesh}). In the former case, each spatial sampling point is denoted as $\bm{Y}_{s_{1},s_{2}} \equiv \bm{Y}(s_{1}, s_{2}) \in \Omega_{\mu, 0}$ (2D) or $\bm{Y}_{s_{1},s_{2},s_{3}} \equiv \bm{Y}(s_{1}, s_{2}, s_{3})\in \Omega_{\mu, 0}$ (3D), where $s_{i} = 0, 1, \dots, n_{i} - 1$ and $i = 1,\dots, n_\mathrm{dim}$. The associated coordinates can be determined in the two-dimensional case as
\begin{equation}
	\bm{Y}_{s_{1}, s_{2}} = \left( \dfrac{(l_{\mathrm{RVE}})_{1}}{n_{1}} s_{1}, \dfrac{(l_{\mathrm{RVE}})_{2}}{n_{2}} s_{2} \right), \quad s_{i} = 0, 1, \dots, n_{i} - 1, \quad i = 1, 2 \, ,
\end{equation}
and in the three-dimensional case in an analogous way. By performing a discrete Fourier transform (DFT) and jumping to the discrete frequency domain, each sampling angular frequency is then characterized by a wave vector, $\bm{\zeta}_{s_{1}, s_{2}} \equiv \bm{\zeta}(s_{1}, s_{2})$ (2D) or $\bm{\zeta}_{s_{1}, s_{2}, s_{3}} \equiv \bm{\zeta}(s_{1}, s_{2}, s_{3})$ (3D), that is defined in the two-dimensional case as
\begin{equation}
	\bm{\zeta}_{s_{1}, s_{2}} = \left( \dfrac{2 \pi}{(l_{\mathrm{RVE}})_{1}} s_{1}, \dfrac{2 \pi}{(l_{\mathrm{RVE}})_{2}} s_{2} \right), \quad s_{i} = 0, 1, \dots, n_{i} - 1, \quad i = 1, 2 \, ,
\end{equation}
and similarly in the three-dimensional case.

After the CROM clustering-based domain decomposition is performed, the RVE domain is  decomposed in $n_{\mathrm{c}}$ material clusters and a cluster-reduced representative volume element (CRVE) is obtained. The clustering-based domain decomposition that characterizes the initial CRVE is hereafter called base clustering. To take advantage of a priori knowledge about the material behavior, the cluster analysis is usually performed independently for each material phase, henceforth designated cluster-reduced material phase (CRMP). Each material cluster groups a given number of voxels and may have an arbitrary-shaped discontinuous domain, as illustrated in Figure~\ref{fig:adaptivity_nomenclature}. Usually, it is assumed that every local field, here generally denoted as $\bm{a}_{\mu}(\bm{y})$, is uniform within each material cluster,
\begin{equation}
	\bm{a}_{\mu}(\bm{y}) = \sum^{n_{c}}_{I=1} \bm{a}^{(I)}_{\mu} \chi^{(I)}(\bm{y}) \, , \qquad \chi^{(I)} (\bm{y}) = \begin{cases} 1 \quad \text{if} \; \; \; \bm{y}\in \Omega^{(I)}_{\mu} \\ 0 \quad \text{otherwise} \end{cases} \, , \label{eq:clusterdefinition}
\end{equation}
where $\bm{a}^{(I)}_{\mu}$ is the homogeneous field in the $I\text{th}$ material cluster and $\chi^{(I)} (\bm{y})$ is the characteristic function of the $I\text{th}$ material cluster. The equilibrium problem solution procedure is then formulated over the CRVE according to the chosen CROM's formulation.

\begin{figure}[hbt]
\centering
\includegraphics[width=0.9\textwidth]{images/CRVE_nomenclature}
\caption{Schematic illustration of a biphasic material (MP1 and MP2), from left to right: Representative volume element (RVE); Spatially discretized RVE in a regular grid of $8\times8$ voxels; Cluster-reduced representative volume element (CRVE) with 6 material clusters.}
\label{fig:adaptivity_nomenclature}
\end{figure}

\subsection{Characteristics and challenges of CROM adaptivity}
Despite the extensive work in the context of AFEMs, all the successful contributions do not translate directly to  CROMs. Although the main concepts and overall framework are similar to a certain extent, the several differences between both approaches raise new challenges that must be addressed. These are discussed below to set the stage for the novel framework proposed in this article.

\begin{figure}[hbt]
\centering
\includegraphics[width=0.9\textwidth]{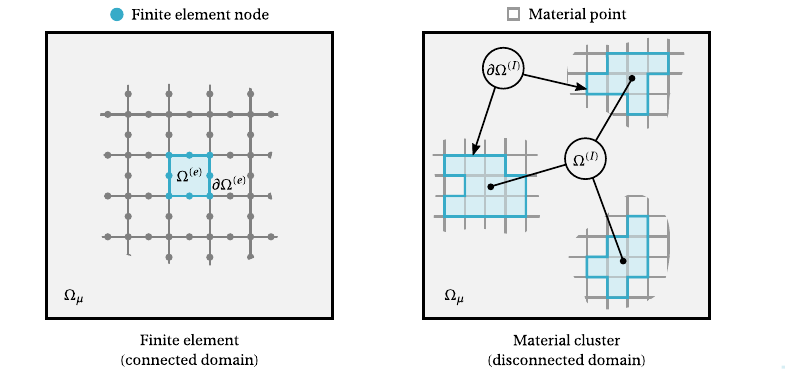}
\caption{Geometrical comparison between a finite element (left) and a material cluster (right).}
\label{fig:fe_vs_cluster}
\end{figure}

The most significant difference resides, of course, in the spatial decomposition. In FEM, the domain $\Omega$ is discretized in a finite set of subdomains called finite elements, $\Omega^{(e)}$, $e=1, ..., n_\text{elem}$. Each finite element, $e$, is a geometrically well-defined connected subspace characterized by a given number of nodes and an equal number of so-called shape functions. The latter are polynomials of a given order that perform the required field interpolations within the element domain. In contrast, in CROMs, the domain $\Omega$ is discretized in a finite set of subdomains called material clusters, $\Omega^{(I)}$, $I=1, ..., n_{\mathrm{c}}$. Each material cluster, $I$, is a geometrically (often) disconnected subspace that groups a given number of points with similar mechanical behavior (see Figure~\ref{fig:fe_vs_cluster}). In accordance, the different fields are usually assumed uniform within the cluster domain (see Equation~\eqref{eq:clusterdefinition}). Another difference that should be remarked concerns the primary unknowns of the formulation. While in FEM the primary unknowns are generally the displacements at nodes (displacement-based formulation), in CROMs they are often taken as the (uniform) strains at the material clusters (strain-based formulation).

Given that fields are assumed to be uniform within each material cluster, recovery-based error estimators that take advantage of optimal sampling points (e.g., SPR-like procedures) are not readily available. The main idea of the alternative recovery-based estimators that attempt to determine a recovered system that is smooth and continuous may be applicable. Still, the proposed contributions are closely tied to the FE formulation. Most importantly, the clusters are, in general, disconnected subdomains, being cumbersome to smoothen a given field over a patch of clusters or even a single cluster. By making use of the equilibrium residuals, the primary approach of residual-based error estimators seems to be most easily translated to the context of CROMs. However, explicit residual-based estimators usually depend on mesh-dependent parameters and often account only for the significant contribution of the so-called jump discontinuities. Given that clusters are usually disconnected subdomains, i.e., each cluster boundary is not a connected path (see Figure~\ref{fig:fe_vs_cluster}), such discontinuities require suitable treatment. In addition, cluster boundaries are not as well-defined as finite element boundaries, the latter easily characterized by a given set of boundary nodes. Implicit residual-based estimators depend significantly on the choice of suitable recovery methods and require the proper treatment of boundary fluxes, thus involving the difficulties mentioned earlier. Other types of estimators related, for instance, with the analysis of constitutive functionals, may be of interest, as the material constitutive behavior of each cluster follows standard procedures.

Due to their heuristic nature, there is a lot of flexibility in the definition of error indicators in the context of CROMs as well. While several error indicators employed in AFEMs may also be readily available in the solution procedure, some are not translatable due to the previously mentioned formulation differences. This is the case of the element aspect ratio or distortion, a geometrical indicator used effectively to deal with localization phenomena and discontinuities in finite element computations (see Figure~\ref{fig:adaptivity_element_distortion}).

Concerning the adaptive discretization procedures, $p$-adaptivity and $hp$-adaptivity are of no interest while assuming the uniformity of fields within the disconnected cluster subdomains. Focus can thus be given to $h$-adaptivity (see Figure~\ref{fig:h_adaptivity_methods}), where the size of clusters may be refined or de-refined according to the accuracy requirements. A significant difference already emerges at this point regarding the size measure. In AFEMs, the finite element size, $h$, is usually taken as the diameter of the smallest circle/sphere that contains the element domain, $\Omega^{(e)}$. Given that the domain of each material cluster is, in general, disconnected, such measure loses its significance. The closest size measure is, for instance, the number of voxels belonging to the cluster.

Besides the overall limitation in terms of engineering practicality, the simplest approach of $r$-refinement is not applicable as clusters are not well-defined by a given set of boundary nodes. At the other end of the degree of adaptivity, the complete mesh regeneration approach can be, in theory, adopted by performing a new cluster analysis and subsequent cluster domain decomposition. However, the crucial transfer of data between the old and new clusterings would be highly cumbersome, as the cluster domains are not only disconnected but can also exhibit a significantly different topology. The more conservative element subdivision seems the most suitable and natural approach to be employed in the context of CROMs. The most prominent difficulties in AFEMs are related to the placement of new nodes and mismatches between adjacent elements, inexistent issues when dealing with material clusters. In fact, the capabilities of different clustering algorithms may be explored to perform an enriched data-based subdivision effectively. Nonetheless, de-refinement procedures may be challenging in terms of data management, primarily due to the ambiguous concept of adjacency when dealing with disconnected cluster domains.

In terms of remeshing criteria, the most popular strategies can be formulated similarly. The same applies to the coupling of the adaptive procedures with the general incremental scheme to solve nonlinear problems, where some specifics of CROMs must be taken into account. A crucial aspect concerns the interaction tensors that frequently emerge in the CROMs' formulation and establish a strain-stress relationship between each pair of clusters. Given the significant costs associated with the computation of these tensors, their update must be efficiently performed every time the clustering is adapted.

Finally, a vital aspect of the reduced-order modeling paradigm is the balance between accuracy and efficiency. The primary objective of clustering adaptivity is, of course, an improvement of the solution's accuracy. As described in Section~\ref{sec:introduction}, of particular importance in modeling path-dependent nonlinear elasto-plastic materials is the phenomenon of strain softening and localization, often the precursor of material failure and fracture. Besides the shortcomings stemming from the lack of adaptivity, note that the main idea underlying the clustering-based (non-local) domain decomposition is, to a certain extent, opposite to the main features of such localized phenomena. Hence, adaptivity is deemed crucial in this context and, perhaps, more challenging. Nonetheless, it must always be present that any accuracy gains are only valuable if the computational costs of the adaptive procedures do not compromise the high efficiency of CROMs. The higher the computational cost, the lower the true value of the accuracy gains compared to standard DNS methodologies. Finding such an accuracy-efficiency equilibrium point in the development of an ACROM framework is the challenge addressed here.

\subsection{Adaptive Clustering-based Reduced-order Modeling framework \label{ssec:acromframework}}
\begin{figure}[hbt]
\centering
\includegraphics[width=0.9\textwidth]{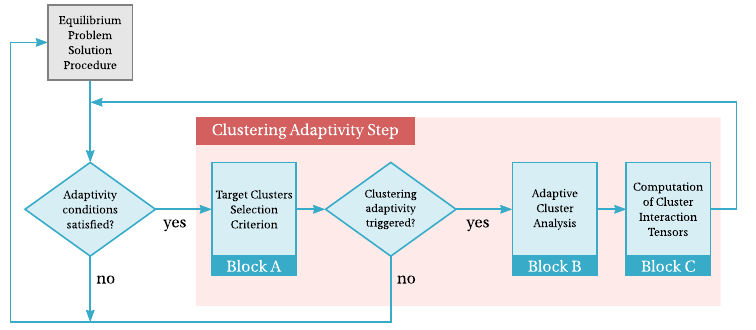}
\caption{Adaptive Clustering-based Reduced-Order Modeling (ACROM) framework. A complete clustering adaptivity step comprises three fundamental blocks: (A) target clusters selection criterion, (B) adaptive cluster analysis and (C) computation of cluster interaction tensors.}
\label{fig:clustering_adaptivity_framework}
\end{figure}

The proposed ACROM framework is illustrated in Figure~\ref{fig:clustering_adaptivity_framework}. A CRMP can be further classified as a static cluster-reduced material phase (SCRMP) if the associated clustering-based decomposition is kept constant throughout the problem's solution, or as an adaptive cluster-reduced material (ACRMP) if clustering adaptivity is allowed.
The first point to be discussed concerns the linkage between this framework and the CROM (microscale) equilibrium problem solution procedure. After the solution of each (macroscale) loading increment is obtained, a given set of adaptivity conditions is verified to check if the clustering adaptivity procedures are carried out. Some possible requirements are proposed here:
\begin{itemize}
	\item \textbf{Clustering adaptivity frequency}. The clustering adaptivity frequency, $\Delta m_{\mathrm{adapt}}$, defines if the ACROM framework is activated after each loading increment ($\Delta m_{\mathrm{adapt}} = 1$), or only every $q$ loading increments ($\Delta m_{\mathrm{adapt}} = q$). This condition could be simply set by default as $\Delta m_{\mathrm{adapt}} = 1$. However, given that the total number of clusters is often limited, a smart choice of this parameter is crucial for achieving an optimal accuracy-efficiency balance, as later discussed. In addition, this parameter can be defined independently for each ACRMP, being the ACROM framework only activated if at least one material phase needs to be evaluated. Note that a similar strategy can be adopted based on the evolution of any relevant variable other than the loading increments (e.g., a relative increase of a damage variable);
	\item \textbf{Maximum consecutive adaptivity steps}. A maximum number of consecutive clustering adaptivity steps can be set so that the equilibrium problem solution procedure continues even if further clustering adaptivity would still be triggered within the following clustering adaptivity step;
	\item \textbf{Threshold number of clusters}. A threshold may be enforced on the number of clusters of each ACRMP, $n_{\mathrm{c}}^{\mathrm{max}}$, to limit the computational cost of the equilibrium problem solution procedure. If the number of clusters surpasses this threshold, then the ACRMP's adaptivity is locked and the associated clustering is kept static during the remaining solution procedure;
	\item \textbf{Minimum localization/damage value}. A minimum significant value for a given localized or damage variable can be defined for each ACRMP. Clustering adaptivity procedures are only activated for that material phase when the variable of interest reaches the minimum significant defined value (e.g., a minimum value of accumulated plastic strain, $\bar{\varepsilon}^{\, p}_{\mathrm{min}}$, associated with the material plastic yielding).
\end{itemize}
If all adaptivity conditions are satisfied for at least one material phase, then a clustering adaptivity step is performed. Otherwise, the equilibrium problem solution procedure continues to the next macroscale loading increment in a standard way.

A clustering adaptivity step comprises three fundamental blocks. The first block (Block A) involves the evaluation of a given error estimator or indicator, i.e., a criterion that selects the target clusters to be adapted. If there is at least one target cluster, an adaptive cluster analysis (Block B) is performed for each target cluster together with a suitable transference of clusters' state-related variables. Otherwise, the equilibrium problem solution procedure continues if the clustering adaptivity is not triggered for any cluster. After performing the adaptive cluster analysis for all target clusters, the CRVE clustering-based domain decomposition has been effectively updated. However, the so-called cluster interaction tensors emerge in most CROMs' formulation and must be updated accordingly (Block C). Each of these blocks is described in detail in the following sections.

Once the clustering adaptivity step is concluded, the adaptivity conditions are re-evaluated. If these continue to be satisfied for at least one material phase, a new clustering adaptivity step is performed. Otherwise, the equilibrium problem solution procedure continues to the next macroscale loading increment. Even if a new clustering adaptivity step is performed, clustering adaptivity might not be triggered for any cluster (see Figure~\ref{fig:clustering_adaptivity_framework}). It is also remarked that the macroscale loading increment where the clustering adaptivity procedure has been performed may be repeated with the updated CRVE before proceeding to the next macroscale loading increment.

\subsection{Block A: Target clusters selection criterion \label{ssec:blockA}}
The first block in the clustering adaptivity step is to determine what material clusters need to be adapted in order to fulfill the desired accuracy requirements. Assume that the clustering adaptivity is designed to capture a given scalar microscale field, $c_{\mu}(\bm{y})$, hereafter called clustering adaptivity feature. This choice may be set independently for each CRMP and can be, for instance, a constitutive state variable associated with localization phenomena (e.g., accumulated plastic strain, $c_{\mu}(\bm{y}) = \bar{\varepsilon}^{p}(\bm{y})$) or a damage variable ($c_{\mu}(\bm{y}) = D(\bm{y})$). Irrespective of the chosen clustering adaptivity feature, having such variable directly available from the CRMP's clusters' state variables is convenient or efficiently computing it from them (e.g., the norm of cluster strain concentration tensor, $c_{\mu}(\bm{y}) = \norm{\bm{\mathsf{H}}}(\bm{y})$, the density of plastic strain energy, $c_{\mu}(\bm{y}) = U^{p}_{\mathrm{dens.}}(\bm{y})$).

\begin{figure}[hbt]
\centering
\includegraphics[width=0.9\textwidth]{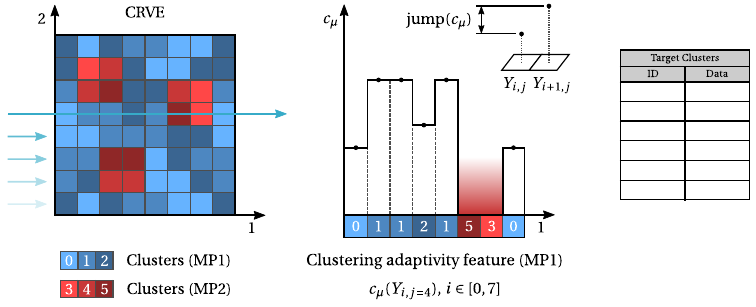}
\caption{Target clusters selection criterion based on the evaluation of the clustering adaptivity feature's, $c_{\mu}$, spatial discontinuities along clusters' boundaries, $\mathrm{jump}(c_{\mu})$. Middle plot shows the clustering adaptivity feature's profile for a given row of voxels ($j=4$) and the associated discontinuities at clusters' boundaries. Targeted clusters are stored together with any data relevant to the following adaptive cluster analysis.}
\label{fig:adaptivity_spatial_discontinuities}
\end{figure}

After setting the suitable clustering adaptivity feature, a criterion to select the target clusters to be adapted is needed. Taking inspiration from the so-called residual error estimators employed in AFEMs, a selection criterion is proposed here based on the evaluation of the spatial discontinuities (or jumps) of $c_{\mu}(\bm{y})$ along the material clusters' boundaries (see Figure~\ref{fig:adaptivity_spatial_discontinuities}). This requires a scanning procedure over the CRVE whose scan directions are collinear and equal in number to the problem dimensions, i.e., two scanning directions (2D problem) or three scanning directions (3D problem).

For simplicity, assume a 2D problem and that direction 1 is currently being scanned, as illustrated in Figure~\ref{fig:adaptivity_spatial_discontinuities}. Assume further that the $j$th row of voxels is being evaluated, i.e., the row whose voxels are defined as $\bm{Y}_{i, j}, \; i = 0, \dots, n_{1}$. Every pair of two consecutive voxels, $(\bm{Y}_{i, j}, \bm{Y}_{i + 1, j}), \; i = 0, \dots, n_{v,1}$, are evaluated as follows\footnote{The last pair of consecutive voxels to be evaluated is $(\bm{Y}_{n_{1}, j}, \bm{Y}_{0, j})$, i.e., the last voxel is paired with the first voxel ($\bm{Y}_{n_{1} + 1, j} \equiv \bm{Y}_{0, j}$).}:
\begin{enumerate}
	\item \textbf{Skip conditions}. Several conditions are defined to avoid unnecessary computations and skip straight to the next pair of voxels if:
	\begin{itemize}
		\item  Both voxels belong to the same material cluster ($\mathrm{cluster}(\bm{Y}_{i, j}) = \mathrm{cluster}(\bm{Y}_{i+1, j})$). This condition has two main outcomes: (1) ensures that only the jump in the boundaries of the clusters is evaluated; (2) computations are avoided when assuming the uniformity of fields within each cluster, i.e., when $\mathrm{jump}(c_{\mu})=0$ is known a priori;
		\item Voxels belonging to different material phases ($\mathrm{MP}(\bm{Y}_{i, j}) \neq \mathrm{MP}(\bm{Y}_{i+1, j}))$. This condition is optional, being nonetheless relevant to (1) comply with a different clustering adaptivity feature for each material phase and (2) perform the clustering adaptivity independently for each CRMP. It can be dropped if the clustering adaptivity should conform to the discontinuity of a given clustering adaptivity feature at the interface between different material phases (e.g., interface phenomena).
	\end{itemize}
	\item \textbf{Jump evaluation}. The clustering adaptivity feature jump can be evaluated and normalized as
	\begin{equation}
		\mathrm{jump}_{r}(c_{\mu}) = \dfrac{ \Big| c_{\mu}(\bm{Y}_{i, j}) - c_{\mu}(\bm{Y}_{i+1, j}) \Big|}{\Big| \mathrm{max}(c_{\mu}(\forall \bm{Y} \in \mathrm{MP})) - \mathrm{min}(c_{\mu}( \forall \bm{Y} \in \mathrm{MP})) \Big|} \, ,
	\end{equation}
being conveniently bounded as $\mathrm{jump}_{r}(c_{\mu}) \in [0,1]$. This normalized ratio quantifies the magnitude of the clustering adaptivity feature discontinuity relative to the maximum amplitude within the material phase. It can be interpreted as an error indicator in the sense that a significant discontinuity is often associated with localization or damage phenomena that demand a more refined clustering to be properly captured;
	\item \textbf{Target selection condition}. The normalized jump is now compared with a user-defined parameter called adaptivity trigger ratio, $\gamma_{\mathrm{ratio}} \in [0,1]$.
	
	If $\mathrm{jump}_{r}(c_{\mu}) >=  \gamma_{\mathrm{ratio}}$, the discontinuity is assumed to be significant and both clusters are marked as target clusters to be adapted. Otherwise, none of the clusters is marked. The adaptivity trigger ratio is thus a parameter that controls the sensitivity of the selection criterion. As $\gamma_{\mathrm{ratio}} \rightarrow 0$, the adaptivity tends to generalize to the whole ACRMP clustering, whereas $\gamma_{\mathrm{ratio}} \rightarrow 1$ tends to focus adaptivity only on clusters placed on regions exhibiting very high discontinuities.
\end{enumerate}

After evaluating all directions and every pair of consecutive voxels, clusters targeted during the CRVE scanning procedure are stored for the following adaptive cluster analysis.

\subsection{Block B: Adaptive cluster analysis \label{ssec:blockB}}
\begin{figure}[hbt]
\centering
\includegraphics[width=0.9\textwidth]{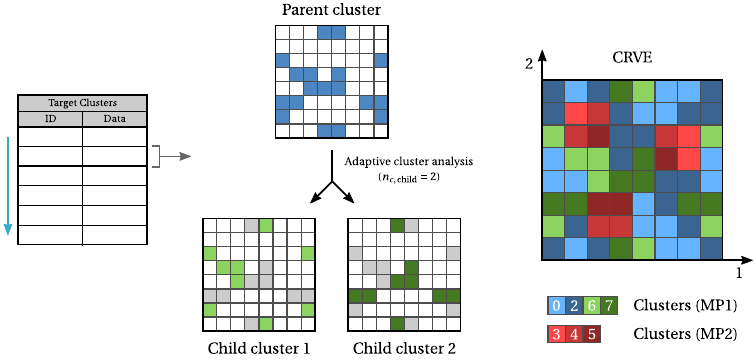}
\caption{Adaptive cluster analysis of each targeted cluster (parent cluster) into a given number of subclusters (child clusters), $n_{\mathrm{c}, \, \mathrm{child}}$, through cluster subdivision. The adaptive cluster analysis can be performed with any available clustering algorithm.}
\label{fig:adaptivity_gacrmp}
\end{figure}

If there is at least one target cluster stemming from Block A, then the clustering adaptivity step proceeds to the following block. A generalized adaptive cluster-reduced material phase (GACRMP) is proposed here. The adaptive cluster analysis of each material phase cluster is performed independently through any chosen unsupervised machine learning clustering algorithm (see Figure~\ref{fig:adaptivity_gacrmp}). The adopted strategy is similar to the $h$-refinement element subdivision approach employed in AFEMs, where target clusters are divided into smaller ones and keep their original boundaries intact. This strategy is highly convenient concerning the transfer of clusters' state-related variables to the refined clustering. After the adaptive cluster analysis and subdivision of a given cluster, each child cluster inherits its parent cluster state-related variables.

Given that any field is usually assumed uniform within each material cluster (see Equation~\eqref{eq:clusterdefinition}), the cluster-level dataset has a null cluster tendency and cannot be used to perform the adaptive cluster analysis. In this context, it is proposed to recover the offline stage data computed at the voxels belonging to the target cluster. It is remarked that the voxel-level dataset employed to perform the clustering procedure may be different from the one employed to perform the base clustering of the CRMP.

\begin{figure}[hbt]
\centering
\includegraphics[width=0.9\textwidth]{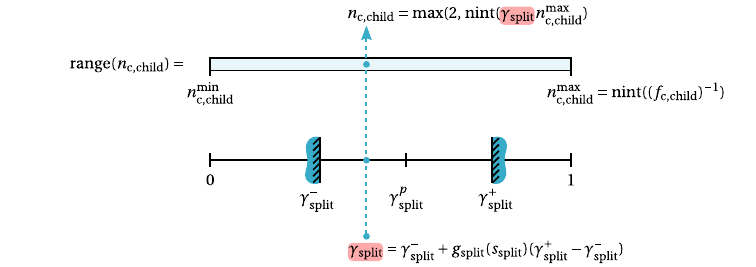}
\caption{Computation of the number of child clusters, $n_{\mathrm{c}, \, \mathrm{child}}$, in which a given target cluster is decomposed through the proposed adaptivity split factor, $\gamma_{\, \mathrm{split}}$.}
\label{fig:adaptivity_split_factor}
\end{figure}

Besides the dataset required to perform the adaptive cluster analysis, it is necessary to define the suitable degree of decomposition, i.e., the number of child clusters created from each target cluster (parent cluster) (see Figure~\ref{fig:adaptivity_split_factor}). The possible number of child clusters ranges from an absolute minimum of $n_{\mathrm{c}, \mathrm{child}}^{\mathrm{min}}=2$ up to a user-defined maximum of $n_{\mathrm{c}, \mathrm{child}}^{\mathrm{max}} = \mathrm{nint}((f_{c,\mathrm{child}})^{-1})$, where $\mathrm{nint}(x)$ denotes the nearest integer rounding function. The parameter $f_{\mathrm{c},\mathrm{child}} \in [0,1]$ denotes the volume fraction of each child cluster relative to the parent cluster by assuming a uniform subdivision, e.g., $f_{\mathrm{c},\mathrm{child}}=0.2$ would result in a maximum of $n_{\mathrm{c}, \mathrm{child}}^{\mathrm{max}} = 5$ child clusters. Given the established range of number of child clusters, the user-defined parameter called adaptivity split factor, $\gamma_{\mathrm{split}} \in [0,1]$, sets the adequate number of child clusters as\footnote{Note, however, that $n_{\mathrm{c}, \mathrm{child}}$ must be lower or equal than the number of voxels of the parent cluster, i.e., the dimension of the voxel-level dataset. Otherwise, the adaptive cluster analysis cannot be performed.}
\begin{equation}
	n_{\mathrm{c}, \mathrm{child}} = \mathrm{max}(2, \, \mathrm{nint}( \gamma_{\mathrm{split}} \, n_{\mathrm{c}, \mathrm{child}}^{\mathrm{max}} )) \, .
\end{equation}
A value of $\gamma_{\mathrm{split}} = 0$ leads to the minimum number of child clusters, $n_{\mathrm{c}, \mathrm{child}}=2$, whereas a value of $\gamma_{\mathrm{split}} = 1$ results in the maximum number of child clusters, $n_{\mathrm{c}, \mathrm{child}}= \mathrm{nint}((f_{c,\mathrm{child}})^{-1})$. Besides the bounded and comprehensible nature of both these parameters, $f_{\mathrm{c},\mathrm{child}}$ and $\gamma_{\mathrm{split}}$, the proposed strategy is convenient as it is independent of the target cluster's number of voxels and often discontinuous spatial arrangement.

The adaptive cluster analysis ends when every target cluster has been adapted and the CRVE clustering-based domain decomposition has been updated.

\paragraph{Remark} To accurately capture localized and damage phenomena, the ACROM framework counteracts the non-local nature of static CROMs. When a given set of material points is grouped into a single cluster, it is implicitly enforced that all follow exactly the same deformation history path. As soon as that cluster is adaptively subdivided, each resulting child cluster is able to follow its own deformation history path, i.e., the non-local connection with the remaining child clusters is broken from that point onwards. In this way, the evolution of material points where localized phenomena occur is not dampened by other material points assumed to have a similar behavior up to that point.

\subsection{Block C: Computation of cluster interaction tensors \label{ssec:blockC}}
\begin{figure}[hbt]
\centering
\includegraphics[width=0.9\textwidth]{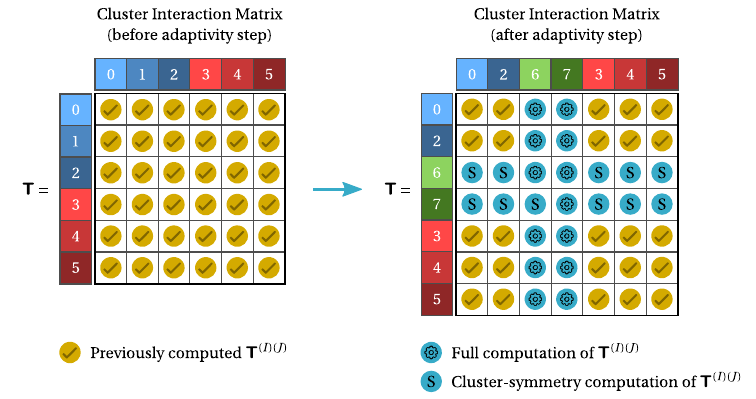}
\caption{Update of the cluster interaction matrix, $\bm{\mathsf{T}}$, consistent with the clustering adaptivity step. By taking advantage of the cluster interaction tensors cluster-symmetry, only those tensors in the columns associated with the new clusters need to be fully computed.}
\label{fig:adaptivity_cit_update}
\end{figure}

The last block of the clustering adaptivity step involves the required update of the so-called cluster interaction tensors, $\bm{\mathsf{T}}^{(I)(J)}, \, I,J = 1,2, \dots, n_{\mathrm{c}}$. These fourth-order tensors play a fundamental role in most CROM's formulation and describe a non-local strain-stress interaction between each pair of clusters. They are usually computed at CROMs' offline stage after performing the CRVE base cluster analysis. As illustrated in Figure~\ref{fig:adaptivity_cit_update}, the complete set of tensors may be conveniently stored in a cluster interaction matrix in terms of computational implementation.

In the particular case of the SCA method, a cluster interaction tensor is defined as
\begin{equation}
	\bm{\mathsf{T}}^{(I)(J)} = \dfrac{1}{f_{(I)} v_{\mu}} \int_{\Omega_{\mu, 0}} \int_{\Omega_{\mu, 0}} \chi^{(I)}(\bm{Y}) \, \chi^{(J)}(\bm{Y}) \, \bm{\Phi}^{0}(\bm{Y}-\bm{Y}') \, \mathrm{d}v' \mathrm{d}v \, , \quad I,J = 1,2,\dots,n_{\mathrm{c}},
	\label{eq:sca_cit_definition}
\end{equation}
where $f^{(I)} = v_{\mu}^{(I)}/v_{\mu}$ is the volume fraction of the $I$th cluster and $\bm{\Phi}^{0}$ is the well-known Green operator associated with the reference material. It transpires from this definition that, from a physical point of view, $\bm{\mathsf{T}}^{(I)(J)}$ represents the influence of the stress in the $J$th cluster on the strain in the $I$th cluster. It is thus understandable that the cluster interaction tensors must be updated according to the new CRVE clustering.

The most naive approach would be to compute all the cluster interaction tensors as done relative to the CRVE base clustering. However, although this option would be perfectly valid from a solution point of view, the cost associated with the computation of these fourth-order tensors is highly significant. Besides the scaling with the number of clusters, where a new cluster implies a new cluster interaction tensor with itself ($I=J$) and with all the remaining clusters ($I \neq J$), the computational cost scales with the spatial discretization of the RVE  in the spatial and/or frequency domains. If the Green operator is assumed to be known in the frequency domain, $\breve{\bm{\mathsf{\Phi}}}^{0}$, the computation of $\bm{\mathsf{T}}^{(I)(J)}$ for a given pair of clusters $I$ and $J$ (see Equation~\eqref{eq:sca_cit_definition}) involves two main steps: (1) the convolution over the cluster $J$, which can be computed in the frequency domain as
\begin{equation}
	\int_{\Omega_{\mu, \, 0}} \chi^{(J)} (\bm{Y}') \, \bm{\mathsf{\Phi}}^{0} (\bm{Y}- \bm{Y}') \, \mathrm{d} v' = \mathscr{F}^{-1} \left( \breve{\chi}^{(J)}(\bm{\zeta}) \, \breve{\bm{\mathsf{\Phi}}}^{0} (\bm{\zeta}) \right) \, ,
	\label{eq:cit_convolution_J}
\end{equation}
where $\mathscr{F}^{-1}$ denotes the inverse Fourier transform; and (2) the spatial integration over the cluster $I$ domain implied in
\begin{equation}
	\int_{\Omega_{\mu, 0}} \chi^{(I)}(\bm{Y}) \left( \int_{\Omega_{\mu, 0}}  \, \chi^{(J)}(\bm{Y}) \, \bm{\Phi}^{0}(\bm{Y}-\bm{Y}') \, \mathrm{d}v' \right) \mathrm{d}v \, .
\end{equation}
The computation of $\bm{\mathsf{T}}^{(I)(J)}$ is then concluded with the product by $(v_{\mu}^{(I)})^{-1} = (f^{(I)} v_{\mu})^{-1}$.

Note that the cluster interaction tensors exhibit a cluster-symmetry in the sense that, from Equation~\eqref{eq:sca_cit_definition},
\begin{equation}
	\bm{\mathsf{T}}^{(J)(I)} = \dfrac{f^{(I)}}{f^{(J)}} \, \bm{\mathsf{T}}^{(I)(J)} \, .
	\label{eq:cit_cluster_symmetry}
\end{equation}
Assuming a CRVE base clustering with $n_{c}$ material clusters, cluster-symmetry can be exploited by only performing a full computation of $0.5 \, n_{\mathrm{c}} ( n_{\mathrm{c}} + 1)$ cluster interaction tensors (lower triangular elements of the cluster interaction matrix). The remaining $0.5 \, n_{\mathrm{c}} (n_{\mathrm{c}} -1)$ cluster interaction tensors, representing $0.5 \, (1+ n_{\mathrm{c}}^{-1}) \times 100 \, (\%)$ of the total number of tensors, may be directly obtained from Equation~\eqref{eq:cit_cluster_symmetry}. This fraction tends to $50\%$ as the number of clusters increases.

The cluster-symmetry of the cluster interaction tensors can be further explored in the context of clustering adaptivity to minimize computational costs. In the following explanation, assume that $I$ and $J$ denote existing clusters before the clustering adaptivity step, and $I^{*}$ and $J^{*}$ denote new clusters stemming from the adaptive cluster analysis. It is proposed that the full computation is only performed    for $\bm{\mathsf{T}}^{(I)(J^{*})}$, i.e., the cluster interaction tensors located in the columns of the cluster interaction matrix associated with new clusters (see Figure~\ref{fig:adaptivity_cit_update}). Besides the cluster interaction tensors $\bm{\mathsf{T}}^{(I^{*})(J^{*})}$, it can be verified that all the remaining new cluster interaction tensors $\bm{\mathsf{T}}^{(I^{*})(J)}$ can be computed directly from cluster-symmetry. This approach not only takes advantage of the cluster-symmetry to reduce the number of full tensor computations, but also avoids repeating the computation of the convolution (see Equation~\eqref{eq:cit_convolution_J}) associated with previously existent clusters. The cluster interaction tensors $\bm{\mathsf{T}}^{(I)(J)}$ can be naturally recovered from the previous clustering. In summary, if the CRVE updated clustering is characterized by $n_{\mathrm{c}}$ previously existent clusters and $n_{\mathrm{c}}^{*}$ new clusters, the update of the cluster interaction tensors requires $0.5 \, n_{c}^{*} \, ( n_{c}^{*} + 2 n_{c} + 1)$ full tensor computations and $0.5 \, n_{\mathrm{c}}^{*} \, ( n_{\mathrm{c}}^{*} + 2 n_{\mathrm{c}} - 1)$ cluster-symmetry tensor computations. The proposed strategy is valid irrespective of the method employed to compute the cluster interaction tensors as long as the cluster-symmetry holds.

\subsection{Adaptive clustering solution rewinding \label{ssec:rewindproc}}
It is well known that the numerical solution of a general equilibrium problem depends on the spatial discretization of the domain. This is particularly important when dealing with path-dependent nonlinear material behavior and/or damage mechanisms, where the proper discretization of specific regions of the domain can be crucial to get accurate predictions. For these reasons, preliminary mesh convergence studies are a standard practice in DNS analyses of a given problem which can then be complemented by adaptive methodologies.

Given the fundamental accuracy-efficiency balance subjacent to CROMs, it is often the case that the clustering-based domain decomposition is not sufficiently refined to assume convergence with respect to spatial discretization. In the context of the ACROM framework, this means that there is a persistent material state history prior to any clustering adaptivity step that will affect the solution accuracy of the remaining loading path. In other words, even if an `optimal' clustering refinement is performed at a given instant of time, the resulting child clusters inherit a less accurate material state. For instance, it is expected that a coarser clustering leads to a rather diffuse modeling of localized phenomena, which means that child clusters inherit underestimations of the magnitude of the localized fields. More importantly, the predictions of the onset and propagation of the material failure are consequently delayed.

\begin{figure}[hbt]
\centering
		\begin{subfigure}[b]{0.9\textwidth}
                \centering
                \includegraphics[width=\textwidth]{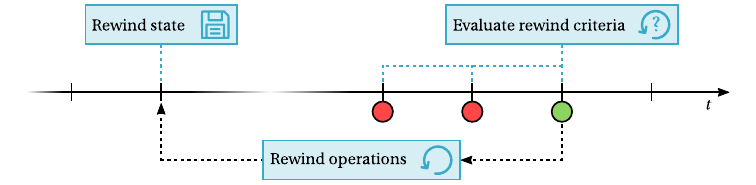}
                \caption{}
                \label{subfig:Adaptivity_rewind_increment}
        \end{subfigure} \vspace*{1.0cm}
        
        \begin{subfigure}[b]{0.95\textwidth}
                \centering
                \includegraphics[width=\textwidth]{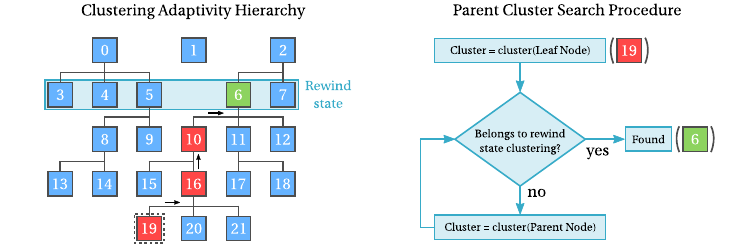}
                \caption{}
                \label{subfig:Adaptivity_rewind_state}
        \end{subfigure}\hfill
\caption{Schematic of solution rewinding procedure after clustering adaptivity: \subref{subfig:Adaptivity_rewind_increment} Solution rewinding main steps, namely (1) storage of rewind state, (2) evaluation of rewind criteria and (3) perform rewind operations; \subref{subfig:Adaptivity_rewind_state} Recovery of clusters state-related variables through a suitable search of clustering adaptivity hierarchy.}
\label{fig:Adaptivity_rewind}
\end{figure}

To address this challenge, a solution rewinding procedure compatible with the ACROM framework is proposed here. This procedure is composed of three main steps, as illustrated in Figure~\ref{subfig:Adaptivity_rewind_increment}. The first step consists in storing the rewind state whenever a given condition is met (e.g., when the material begins yielding plastically). This essentially involves taking a snapshot of the solution state, saving the loading path state, the material homogenized strain and stress, the clusters state variables, and so on. Once the rewind state is stored, the second step consists in evaluating one or more rewind criteria at each loading increment, i.e., criteria that establish when the solution must be rewound after at least one clustering adaptivity step has been performed. If the rewind criteria are met, the third step performs the actual rewind operations, returning the solution to the rewind state but with the updated clustering. This requires that the new clustering somehow recovers the material state variables associated with the rewind state.

The clusters' state variables recovery process is illustrated in Figure~\ref{subfig:Adaptivity_rewind_state}. Given the adaptive cluster analysis described in Section~\ref{ssec:blockB}, the clustering adaptivity hierarchy of the GACRMP can be built throughout the problem solution. Therefore, the recovery process requires finding the rewind state's parent cluster of each cluster from which the associated state variables are transferred. The pseudo-code of this parent cluster search procedure is presented in Figure~\ref{subfig:Adaptivity_rewind_state} and consists of an upward hierarchical search for each cluster.

\subsection{Additional adaptivity procedures \label{ssec:addadaptproc}}

Some additional procedures are proposed here to complement the ACROM framework previously described. These are entirely optional but can effectively improve the clustering adaptivity and/or provide a more practical tool to the analyst.

\begin{itemize}
	\item \textbf{Cluster adaptivity level}. The adaptivity level of a given cluster $I$, $\lambda^{(I)}$, describes its depth of refinement throughout the equilibrium problem solution procedure, being $\lambda^{(I)} = 0$ associated with the base clustering. If cluster $I$ is targeted and subjected to adaptive cluster analysis, each child cluster inherits its adaptivity level incremented by 1. In this sense, the adaptivity level measures how many clustering adaptive steps led to its creation, i.e., how far it is from the base clustering domain decomposition. This information can be useful in two different ways. In the first place, a user-defined parameter, $\lambda_{\mathrm{max}}$, may be enforced to limit the adaptivity level of all clusters belonging to a given ACRMP. This would be an additional target selection condition in the sense that any cluster $I$ could only be targeted if $\lambda^{I}<\lambda_{\mathrm{max}}$ (see Section~\ref{ssec:blockA}). In the second place, the adaptivity level may be used to enforce a certain degree of uniformity concerning the domain clustering adaptivity. A user-defined parameter, $\lambda_{\mathrm{max}}^{\Delta}$, can be defined to limit the difference between adaptivity levels of adjacent clusters. Assume the generic pair of consecutive voxels $(\bm{Y}_{i,j}, \bm{Y}_{i+1,j})$ described in Section~\ref{ssec:blockA} and that $I=\mathrm{cluster}(\bm{Y}_{i, j})$ and $J=\mathrm{cluster}(\bm{Y}_{i+1, j})$. If $\mathrm{jump}_{r}(c_{\mu}) >= \gamma_{\mathrm{ratio}}$ and $| \lambda^{(I)} - \lambda^{(J)} | <= \lambda_{\mathrm{max}}^{\Delta}$, then both clusters are marked as target clusters to be adapted. However, if $| \lambda^{(I)} - \lambda^{(J)} | > \lambda_{\mathrm{max}}^{\Delta}$, only the cluster with the lowest adaptivity level, $\mathrm{min}(\lambda^{(I)}, \lambda^{(J)})$, is marked as target cluster.
	\item \textbf{Cluster minimum number of voxels}. Irrespective of the microstructure under analysis, the smallest cluster possible is composed of a single voxel. If a proper domain discretization is performed, usually involving thousands of voxels, it is most likely that a single-voxel cluster does not contribute significantly to the accuracy of the problem solution. However, either in the computation of the base clustering or in each clustering adaptivity step, an additional cluster may contribute to a significant increase in the number of cluster interaction tensors and associated computation cost. This reasoning can be extended to a given minimum number of voxels per cluster, a user-defined parameter that naturally depends on the spatial discretization of the domain and on the minimum dimension considered significant to capture the material behavior accurately. This threshold is then enforced as an additional target selection condition, being a given cluster targetable only if the associated number of voxels is greater or equal to the prescribed minimum;
	\item \textbf{Dynamic adaptivity split factor}. In Section~\ref{ssec:blockB}, the adaptivity split factor, $\gamma_{\mathrm{split}}$, is proposed as a user-defined parameter that defines the number of child clusters of every target cluster in a given ACRMP. This concept can be further explored by recognizing that certain target clusters should be more refined on a given clustering adaptivity step than others, depending on the error estimator or indicator being employed. For instance, target clusters associated with higher discontinuities of the clustering adaptivity feature, $\mathrm{jump}_{r}(c_{\mu})$, most likely demand a greater refinement. In this sense, a magnitude associated to each target cluster, $s_\mathrm{split}$, can be conveniently defined as $s_\mathrm{split} = \mathrm{max}(\mathrm{jump}_{r}(c_{\mu})) - \gamma_{\mathrm{ratio}}$,  $s_\mathrm{split} \in [0,1-\gamma_{\mathrm{ratio}}]$, where $\mathrm{max}(\mathrm{jump}_{r}(c_{\mu}))$ is the maximum jump associated with each target cluster and stored during the CRVE scanning procedure (see Section~\ref{ssec:blockA}). This data may be used in three different adaptivity procedures:
	\begin{itemize}[$\diamond$]
		\item \textbf{Enforcement of the number of clusters.} When no metric is employed to distinguish the different target clusters, the threshold number of clusters can only be evaluated after the clustering adaptivity step (see Section~\ref{ssec:acromframework}). This means that the total number of clusters may largely exceed the prescribed threshold, given that all target clusters are split. By sorting the target clusters in descending order of magnitude, $s_\mathrm{split}$, the remaining target clusters can be discarded as soon as the total number of clusters surpasses the prescribed threshold. The deviation of the total number of clusters relative to the prescribed threshold is at most $n_{\mathrm{c}, \mathrm{child}}^{\mathrm{max}}-1$ in this case;
		\item \textbf{Dynamic split factor magnitude function.} An additional user-defined parameter called adaptivity split factor amplitude, $\gamma_{\mathrm{split}}^{\Delta}$, may be set around the prescribed adaptivity split factor, $\gamma_{\mathrm{split}}^{p}$ (see Figure~\ref{fig:adaptivity_split_factor}). The range of the adaptivity split factor is then characterized by a lower bound, $\gamma_{\mathrm{split}}^{-} = \gamma_{\mathrm{split}}^{p} - 0.5 \, \gamma_{\mathrm{split}}^{\Delta}$, and an upper bound, $\gamma_{\mathrm{split}}^{+} = \gamma_{\mathrm{split}}^{p} + 0.5 \, \gamma_{\mathrm{split}}^{\Delta}$. The suitable number of child clusters associated with each target cluster may then be set as
\begin{equation}
	\gamma_{\mathrm{split}} = \gamma_{\mathrm{split}}^{-} + g_{\mathrm{split}}(s_{\mathrm{split}}) \, \gamma_{\mathrm{split}}^{\Delta} \, ,
\end{equation}
where $g_{\mathrm{split}} = g_{\mathrm{split}}(s_{\mathrm{split}})$ can be any monotonically increasing scalar function that satisfies the boundary conditions $g_{\mathrm{split}}(0) = 0$ and  $g_{\mathrm{split}}(1-\gamma_{\mathrm{ratio}}) = 1$. This means that a target cluster that has been marked with the minimum possible magnitude ($s_{\mathrm{split}}=0$) is adapted with $\gamma_{\mathrm{split}} =  \gamma_{\mathrm{split}}^{-}$, whereas target clusters with $s_{\mathrm{split}} > 0$ are adapted with a $\gamma_{\mathrm{split}}$ proportional to the associated $s_{\mathrm{split}}$. Here a general power function is proposed as
\begin{equation}
	\gamma_{\mathrm{split}} = (1 - \gamma_{\mathrm{ratio}})^{-n} s^{n}_{\mathrm{split}} \, ,
\end{equation}
where $n \in [0,+\infty]$ (see Figure~\ref{fig:Adaptivity_split_factor_gsplit_power}). A value of $n=1.0$ can be set by default, which means that the number of child clusters increases linearly with the magnitude. A value of $n>1.0$ can be set to decrease the split factor sensitivity, i.e., to promote an increasing number of clusters in the higher magnitude range. In contrast, $n<1.0$ can be set to increase the split factor sensitivity and have an increased number of clusters from the low magnitude range;	
		\item \textbf{Lower-valued cluster split factor.} By default, it has been assumed that both target clusters associated to a given spatial discontinuity, $\mathrm{cluster}(\bm{Y}_{i, j})$ and $\mathrm{cluster}(\bm{Y}_{i+1, j})$, are evaluated in a similar way concerning the discontinuity magnitude. However, it may be of interest to differenciate the number of child clusters between the lower and higher-valued cluster, namely decreasing the number of clusters created in the lower-valued side of the jump. To achieve this behavior, the magnitude associated to the lower-valued cluster, $s_{\mathrm{split}}^{\mathrm{low}}$, can be defined as $s_{\mathrm{split}}^{\mathrm{low}} = \theta \, s_{\mathrm{split}}$, where $\theta \in [0,1]$ sets the desired differentiation.
	\end{itemize}
\end{itemize}

\begin{figure}[hbt]
\centering
\includegraphics[width=0.9\textwidth]{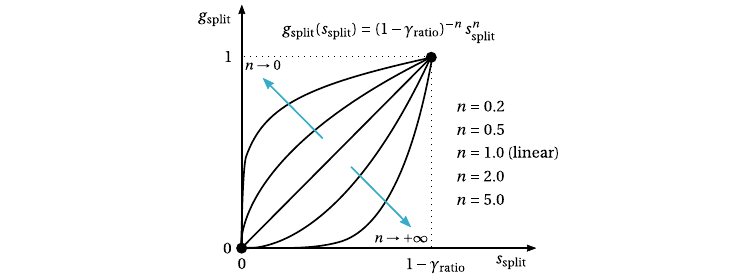}
\caption{Proposed dynamic adaptivity split factor magnitude function. Three behaviors are available: (1) $n=1.0$ means that the number of child clusters increases linearly with the magnitude associated with the parent cluster, (2) $n\rightarrow +\infty$ promotes an increasing number of child clusters in the high magnitude range (low sensitivity) and (3) $n\rightarrow 0$ raises an increasing number of child clusters from the low magnitude range (high sensitivity).}
\label{fig:Adaptivity_split_factor_gsplit_power}
\end{figure}

\subsection{Summary of hyperparameters}
It is usually the case that advanced numerical methods and/or frameworks account for several hyperparameters (or user-defined parameters) (i) whose meaning is difficult to grasp, (ii) that are cumbersome to calibrate properly, and (iii) that may significantly impact the results when slightly changed. Despite the hyperparameters introduced in the proposed ACROM framework (see summary on Table~\ref{tab:acromparameters}), it is shown in the following section that key to obtaining successful results is to adopt the proper clustering adaptivity strategy based on the physical understanding of the problem. Once the suitable clustering adaptivity strategy is established, the choice of the  hyperparameters is a consequence rather than an optimization problem. Some numerical experience shows that the role of each parameter can be easily understood and that results exhibit a lower parameter sensitivity as long as the proper strategy is adopted.

\begin{table}[hbt]
\caption{Essential hyperparameters of the ACROM framework. Insertion: target clusters selection criterion (Block A) and adaptive cluster analysis (Block B).}
\label{tab:acromparameters}
\centering
\setlength{\tabcolsep}{0.05cm}
\renewcommand{\arraystretch}{1.3}
\begin{tabular}{cccl}
	 \textbf{Parameter} & \textbf{Notation} & \textbf{Insertion} & \textbf{Meaning} \\ \midrule
	 \multirow{2}{3cm}{\centering Adaptivity frequency} & \multirow{2}{*}{$\Delta m_{\mathrm{adapt}}$} & \multirow{2}{3cm}{\centering Adaptivity conditions}  & \multirow{2}{7cm}{Controls when adaptivity procedures take place relative to loading increments} \\
	  & & & \\ \midrule
	 \multirow{2}{3cm}{\centering Adaptivity trigger ratio} & \multirow{2}{*}{$\gamma_{\mathrm{ratio}}$} & \multirow{4}{*}{Block A} & \multirow{2}{8cm}{Sets the sensitivity of the selection criterion (e.g., significant field discontinuity)} \\
	 & & & \\
	 \multirow{2}{3cm}{\centering Maximum adaptive level} & \multirow{2}{*}{$\lambda_{\mathrm{max}}$} &  & \multirow{2}{7.5cm}{Limits the cluster level of refinement relative to base clustering} \\
	 & & & \\ \midrule
	 \multirow{2}{3.5cm}{\centering Maximum number of child clusters} & \multirow{2}{*}{$n_{\mathrm{c}, \mathrm{child}}^{\mathrm{max}}$} & \multirow{6}{*}{Block B} & \multirow{2}{8cm}{Sets the  maximum allowed number of child clusters created from each target cluster} \\
	 & & & \\
	 \multirow{2}{3cm}{\centering Adaptivity split factor} & \multirow{2}{*}{$\gamma_{\mathrm{split}}$} &  & \multirow{2}{8cm}{Controls the number of child clusters created from each target cluster${}^{(a)}$}  \\
	 & & & \\
	 \multirow{2}{3cm}{\centering Adaptivity split factor amplitude} & \multirow{2}{*}{$\gamma_{\mathrm{split}}^{\Delta}$} &  & \multirow{2}{8cm}{Controls the variation of the number of child clusters created from each target cluster${}^{(b)}$} \\
	 & & & \\ 
	  \bottomrule
	    \multicolumn{4}{l}{\footnotesize $^{(a)}$ Set as $\gamma_{\mathrm{split}} = 1.0$ if using a static adaptivity split factor, i.e., $n_{\mathrm{c}, \mathrm{child}} = \mathrm{max}(2, \, \mathrm{nint}( \, n_{\mathrm{c}, \mathrm{child}}^{\mathrm{max}} ))$.} \\[-2pt]
	    \multicolumn{4}{l}{\footnotesize $^{(b)}$ Only required if using a dynamic adaptivity split factor.}
\end{tabular}
\end{table}

\section{Numerical results and discussion \label{sec:results}}
This section presents some numerical examples aiming to illustrate the performance of the proposed ACROM framework. As described in Section~\ref{sec:methodology}, the Self-Consistent Clustering Analysis (SCA) method \cite{liu2016} is extended to its adaptive counterpart, coined Adaptive Self-Consistent Clustering Analysis (ASCA), and explored throughout the following numerical assessment.

\subsection{Computational setting}
The numerical simulations are performed with CRATE (Clustering-based Nonlinear Analyses of Materials), an object-oriented Python program that performs multi-scale nonlinear analyses through clustering-based reduced-order models. This program has been fully designed and coded by Bernardo P. Ferreira, and its initial version will soon be released. A thorough description of this program's object-oriented design and application will also be published so that it can be easily exploited and extended by the interested research community.

\begin{figure}[hbt]
\centering
\includegraphics[width=0.6\textwidth]{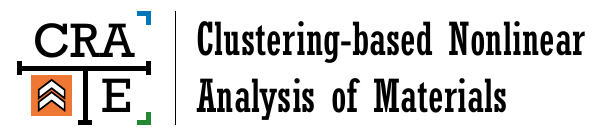}
\caption{Logo of CRATE (Clustering-based Nonlinear Analyses of Materials).}
\label{fig:crate_logo}
\end{figure} 

Unless otherwise stated, all the numerical simulations presented in this section are run with 1 CPU core in a personal computer with the following specifications: CPU Intel Core i7-6800K (3.40GHz $\times$ 6 cores/12 threads, Broadwell), RAM 32GB ($4 \times 8$GB) DDR4 @ 2400MHz. Nonetheless, some third-party (low-level) packages may perform multi-processor computations at certain operations.

\subsection{Particle-matrix composite}
\begin{figure}[hbt]
\centering
\includegraphics[width=0.9\textwidth]{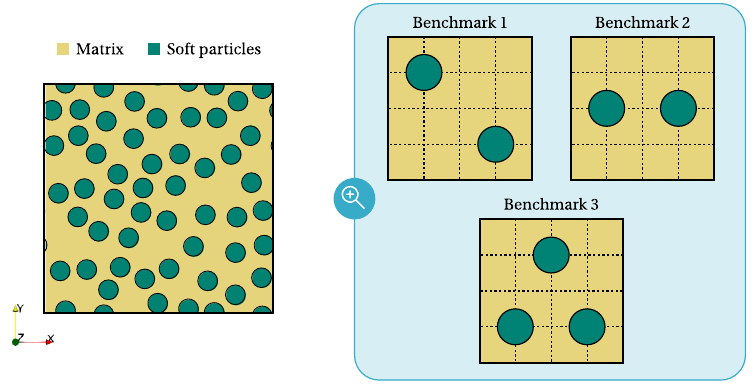}
\caption{2D RVE of biphasic material characterized by randomly distributed circular particles ($f=30\%$) embedded in a matrix ($f=70\%$), and three SVEs representing different particle spatial arrangements found in the biphasic material RVE.}
\label{fig:softparticlescomposite}
\end{figure}

\begin{figure}[hbt]
\centering
\includegraphics[width=0.9\textwidth]{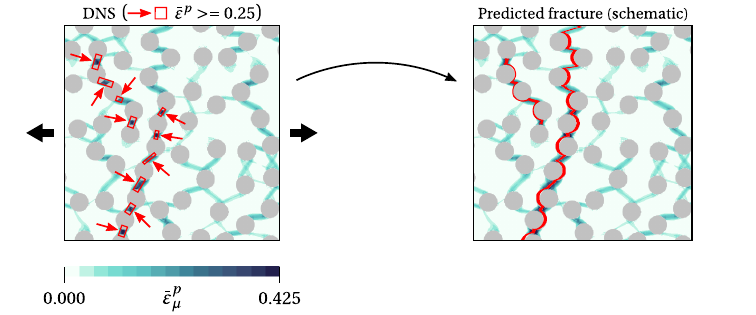}
\caption{Simplified fracture criterion assumed for the particle-matrix composite: fracture occurs when 0.5\% of the matrix material phase surpasses an accumulated plastic strain of $\bar{\varepsilon}^{p}=0.25$. On the left, the accumulated plastic strain field at fracture was predicted with the DNS solution. On the right, a schematic illustration of the fracture propagation and failure of the material's load-bearing capacity.}
\label{fig:2d_frc_acc_p_strain_fracture}
\end{figure}

In order to follow the numerical examples of the original SCA paper \cite{liu2016} and avoid unnecessary complexity in the demonstration of the ACROM framework, a biphasic material characterized by randomly distributed circular particles ($f=30\%$) embedded in a matrix ($f=70\%$) is here considered (see Figure~\ref{fig:softparticlescomposite}). The matrix material phase constitutive behavior is assumed isotropic and elasto-plastic, governed by a standard von Mises associative flow rule with isotropic strain hardening, while the particle material phase is assumed isotropic and elastic. However, in contrast with the SCA paper \cite{liu2016}, particles are assumed soft to promote localized shear yielding and failure in the matrix. The material phase´s properties are summarized in Table~\ref{tab:mat_properties}. In addition, a simplified ductile fracture criterion is formulated according to which fracture occurs when $0.5\%$ of the matrix material phase surpasses an accumulated plastic strain of $0.25$ (see Figure~\ref{fig:2d_frc_acc_p_strain_fracture}). It is remarked that such a simple fracture criterion does not affect the constitutive behavior of the material and is only employed here to compare the failure prediction between different solution methods.

\begin{table}[hbt]
\caption{Constitutive model and properties of particle-matrix composite's material phases. Notation: Young's Modulus ($E$), Poisson's ratio ($\nu$), Yield stress ($\sigma_{y}$) and accumulated plastic strain ($\bar{\varepsilon}^{p}$).}
\label{tab:mat_properties}
\centering
\setlength{\tabcolsep}{0.35cm}
\renewcommand{\arraystretch}{1.3}
\begin{tabular}{ccccc}
	 Phase & Model & $E$(MPa) & $\nu$ & $\sigma_{y}$(MPa)  \\ \midrule
	 Matrix & von Mises (isotropic) & 100 & 0.3 & $\sigma_{y} = 0.5 + 0.2(\bar{\varepsilon}^{p})^{0.4}$ \\
	 Particles & Elastic (isotropic) & 1 & 0.19 & - \\
	  \bottomrule 
\end{tabular}
\end{table}

Given the illustrative purposes of this numerical assessment and without any loss concerning the ACROM framework application, the 2D case under infinitesimal strains is explored throughout this section. Problems involving more complex loading paths and material constitutive behavior are left for future work, as discussed in Section~\ref{sec:conclusion}. The RVE of the particle-matrix composite is shown in Figure~\ref{fig:softparticlescomposite}, together with three statistical volume elements (SVEs) that resemble different particles' spatial arrangements found in the actual RVE. The material is assumed to undergo a macroscale strain-driven uniaxial traction loading defined as $\varepsilon_{xx}=5.0 \times 10^{-2}$ ($\varepsilon_{ij}=0.0$ for $ij \neq xx$) and prescribed in a total of 200 increments. Under these loading conditions, the composite's toughness is computed as the integral of the composite stress-strain curve before fracture prediction.

\paragraph{Remark} The authors' are aware that performing a first-order homogenization and assuming microscale periodic boundary conditions in the presence of localized phenomena is not accurate from a physical point of view. Although such a topic has received a lot of attention in recent years, such considerations can be disregarded in this paper, where the DNS is assumed to deliver the `high-fidelity' (reference) solution to evaluate the CROMs numerical accuracy.

\subsection{Solution methods and error assessment}
Three solution methods are considered in the numerical assessment presented in the following sections: a direct numerical solution (DNS), here adopted as the reference (or `high-fidelity') solution, a Self-Consistent Clustering Analysis (SCA) solution, following the original paper \cite{liu2016}, and the Adaptive Self-Consistent Clustering Analysis (ASCA), the adaptive counterpart of the SCA obtained with the ACROM framework proposed in this paper.

The direct numerical solution (DNS) is obtained with a Finite Element Method (FEM) first-order hierarchical multi-scale model based on computational homogenization. This solution is carried out with LINKS (Large Strain Implicit Non-linear Analysis of Solids Linking Scales), an implicit multi-scale finite element Fortran code developed by the CM2S research group at the Faculty of Engineering of the University of Porto\footnote{Bernardo P. Ferreira is an active member of the CM2S research group, led by F.M. Andrade Pires, and one of the developers of LINKS.}. LINKS is also employed to perform the linear elastic DNS simulations required in the offline-stage of both SCA and ASCA (see \ref{app:femrgmesh} for finite element mesh compatibility and element averaging procedures). Both SCA and ASCA solutions are obtained with CRATE as previously described. Periodic boundary conditions are assumed and a refined spatial discretization of $n_{v}=400 \times 400$ voxels is considered for all numerical examples. Without loss of generality, only the matrix material phase is subject to clustering adaptivity in the ASCA solution. The base clustering of both material phases is performed with the well-known $K$-Means clustering (standard Lloyds' algorithm) with 10 centroid initializations computed with $K$-Means++. In turn, the adaptive cluster analysis of the matrix material phase relies on the Mini-Batch $K$-Means with the same 10 centroid initializations computed with $K$-Means++. Both of these algorithms implementations are readily available from Scikit-learn, a Python module integrating a wide range of state-of-the-art machine-learning algorithms for medium-scale supervised and unsupervised problems \cite{pedregosa2012}.

Concerning the error assessment, two different metrics are adopted. On the one hand, the relative error, $\epsilon$, is employed to compare the macroscale homogenized response and toughness, being always defined as
\begin{equation}
	\epsilon \, (\%) = \left| \dfrac{a_{\mathrm{CROM}} - a_{\mathrm{DNS}}}{a_{\mathrm{DNS}}} \right| \times 100 \, ,
\end{equation}
where $a$ denotes the scalar quantity of interest. On the other hand, the root-mean-square error (RMSE) is employed to compare the microscale field solutions (see \ref{app:femrgmesh}), being defined as
\begin{equation}
	\mathrm{RMSE} \, (a_{\mu}) = \sqrt{\dfrac{1}{n_{v}} \sum_{i=1}^{n_{v}} \left( a_{\mu, \mathrm{CROM}} - a_{\mu, \mathrm{DNS}} \right)_{i}^{2}} \, ,
\end{equation}
where $a_{\mu}$ denotes any microscale scalar quantity. Note that the RMSE has the units of the actual quantity being evaluated. Two microscale local scalar quantities are evaluated in the following analyses: (1) the accumulated plastic strain, $\bar{\varepsilon}^{p}$, a state variable of the von Mises constitutive model governing the behavior of the matrix material phase; and (2) the accumulated plastic strain energy density defined as $U^{p}_{\mathrm{dens}} = \int \Delta \bm{\sigma} \, \mathrm{d} \bm{\varepsilon}^{p}$. Both quantities are intrinsically associated with the highly localized plastic yielding of the composite as well as with the adopted fracture criterion.

\subsection{Benchmark analysis}
The benchmarks presented in Figure~\ref{fig:softparticlescomposite} promote simple plastic strain localization patterns that can be easily visualized. Despite this simplicity, accurately capturing such localized phenomena already poses important challenges to standard CROMs. For these reasons, they are selected as the first illustration of the ACROM framework application and performance without introducing additional complexities.

\begin{figure}[hbt]
\centering
	\includegraphics[width=0.95\textwidth]{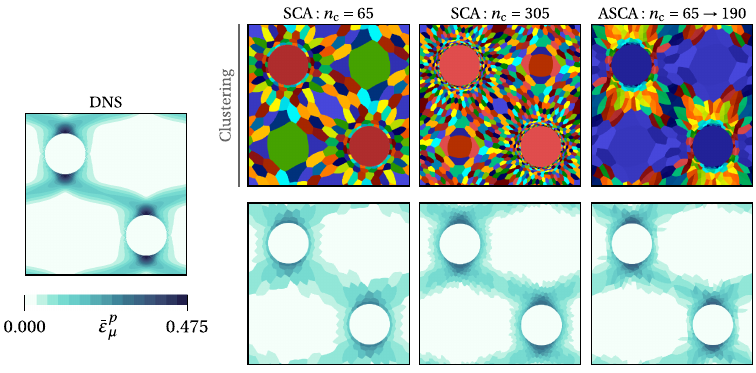}
\caption{Comparison of the clustering-based domain decomposition and local accumulated plastic strain field at the end of the deformation path of Benchmark 1.}
\label{fig:bm_1_local_acc_p_strain}
\end{figure}

As shown in Figure~\ref{fig:bm_1_local_acc_p_strain}, benchmark 1 promotes a single plastic strain localized band between both particles. The first important task consists in identifying, as well as possible, what are the main evolution steps of the field of interest. The evolution of the accumulated plastic strain field in this benchmark has essentially two main steps: (1) the initial plastic yielding of the matrix at the particle-matrix interface (transverse plane relative to loading direction), and (2) the propagation of plastic yielding that leads to the band connecting both particles. In accordance, it is expected that the `optimal' clustering should exhibit a more refined domain decomposition on these regions to accurately capture the material response.

The second task consists in choosing a suitable number of initial clusters. Such choice is by no means unique and should ensure a minimal accuracy upon which the adaptivity can be properly developed. As a rule of thumb, the base clustering should at the very least capture the main regions where the field of interest shows some significant evolution. In this benchmark, the minimal number of clusters has been set to $n_{\mathrm{c}}=65$ (60 clusters in the matrix ACRMP and 5 clusters in the particles SCRMP). From Figure~\ref{fig:bm_1_local_acc_p_strain}, it can be verified that the SCA solution with $n_{\mathrm{c}}=65$ can capture the main regions of plastic strain but in a rather diffuse manner.

Last but not least, it remains to define how should the ACROM framework be effectively employed and set the associated parameters accordingly. In the first place, the clustering adaptivity feature for the matrix ACRMP is set as the accumulated plastic strain, $c_{\mu}(\bm{y})=\bar{\varepsilon}^{p}$. At this point, the clustering adaptivity strategy must be defined based on the available knowledge about the problem physics, namely the main evolution steps already described. Here, clustering adaptivity should follow the plastic strain field front propagation, a similar strategy to that commonly used in phase-field methodologies. It is assumed that a 10\% normalized accumulated plastic strain discontinuity is significant by setting $\gamma_{\mathrm{ratio}}=0.1$ and the number of child clusters is simply set constant and equal to $n_{\mathrm{c}, \mathrm{child}}=2$ ($f_{\mathrm{c}, \mathrm{child}}=0.5$, $\gamma_{\mathrm{split}}=1.0$). Clustering adaptivity procedures should be activated as soon as the matrix begins to yield plastically and the clustering adaptivity frequency is set to $\Delta m_{\mathrm{adapt}} = 15$ so that adaptivity steps may efficiently accompany the plastic strain field front propagation.

Besides the DNS solution, three different SCA solutions are selected to provide insightful comparisons: a `coarse' solution with $n_{\mathrm{c}}=65$, a `medium' solution with $n_{\mathrm{c}}=185$ and a `fine' solution\footnote{The number of clusters of the SCA `fine' solution has been selected so that the relative error of the macroscale homogenized response is below 5\% with respect to the FEM DNS solution.} with $n_{\mathrm{c}}=305$. To focus the attention on the matrix material phase and perform a fair comparison between the different solution methods, the number of clusters of the particles SCRMP is always set to 5 (approximately proportional to the particles volume fraction in the SCA `coarse' solution). To evaluate the performance of the ASCA solution, the initial number of clusters is set to $n_{\mathrm{c}}=65$, similar to the SCA `coarse' solution, and the threshold number of clusters is set to $n_{\mathrm{c}}=185$, similar to the SCA `medium' solution. The maximum allowed adaptivity level is set as $\lambda_{\mathrm{max}}=2$, i.e., each cluster can only be subject to a maximum of two adaptive cluster analyses. This limit is important to avoid exhausting the available number of clusters before all the critical domain regions are properly adapted.

\begin{figure}[hbt]
\centering
		\begin{subfigure}[b]{0.495\textwidth}
                \centering
                \includegraphics[width=\textwidth]{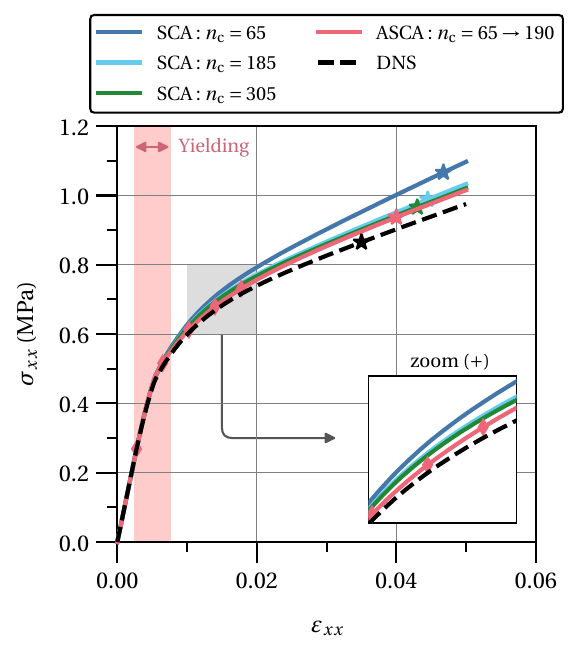}
                \caption{}
                \label{subfig:bm_1_stress_xx_strain_xx}
        \end{subfigure}
        \begin{subfigure}[b]{0.495\textwidth}
                \centering
                \includegraphics[width=\textwidth]{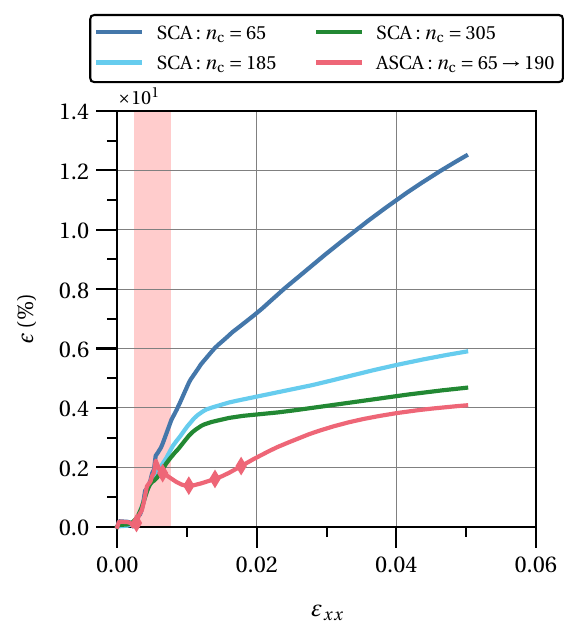}
                \caption{}
                \label{subfig:bm_1_stress_xx_strain_xx_error}
        \end{subfigure}\hfill
\caption{Comparison of the macroscale homogenized response of Benchmark 1 predicted by different solution methods: \subref{subfig:bm_1_stress_xx_strain_xx} Homogenized stress; \subref{subfig:bm_1_stress_xx_strain_xx_error} Relative error with respect to the FEM DNS solution. Symbology: clustering adaptivity step ($\blacklozenge$), fracture criterion prediction ($\bigstar$).}
\label{fig:bm_1_stress_xx}
\end{figure}

\begin{figure}[hbt]
\centering
		\begin{subfigure}[b]{0.490\textwidth}
                \centering
                \includegraphics[width=\textwidth]{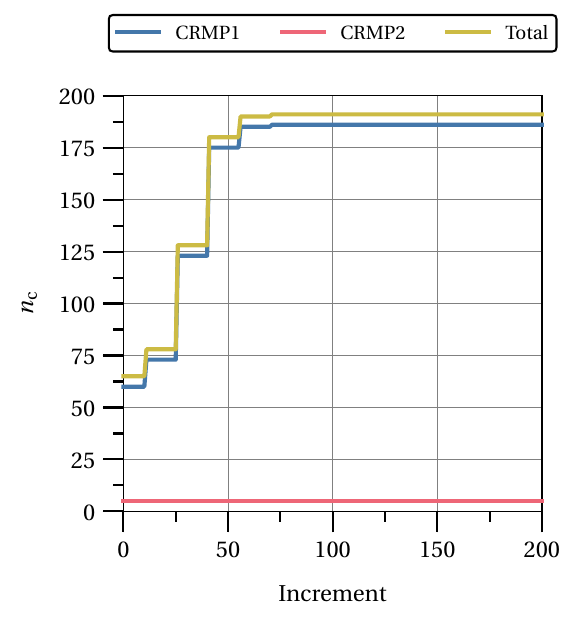}
                \caption{}
                \label{subfig:bm_1_adapt_n_clusters}
        \end{subfigure}
        \begin{subfigure}[b]{0.495\textwidth}
                \centering
                \includegraphics[width=\textwidth]{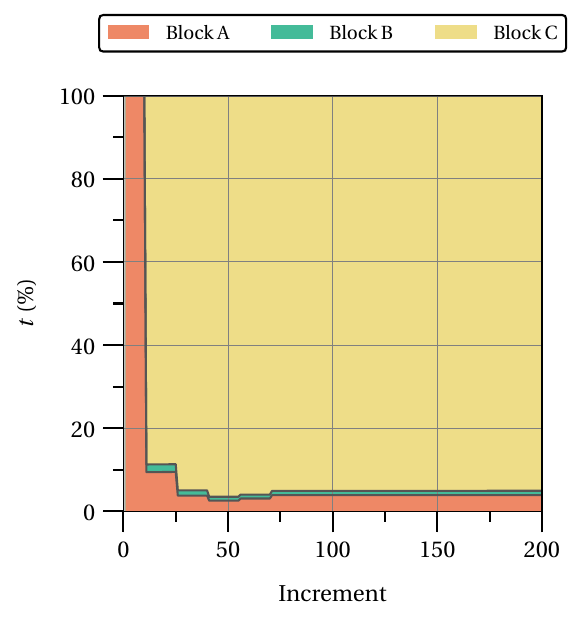}
                \caption{}
                \label{subfig:bm_1_adapt_phases_time_perc}
        \end{subfigure}\hfill
\caption{Evolution of clustering adaptivity related metrics in Benchmark 1: \subref{subfig:bm_1_adapt_n_clusters} Number of clusters; \subref{subfig:bm_1_adapt_phases_time_perc} Relative time of each block of the clustering adaptivity step with respect to the total time spent in clustering adaptivity procedures.}
\label{fig:bm_1_adaptivity_metrics}
\end{figure}

The comparison between the macroscale homogenized response predicted by the different solution methods is shown in Figure~\ref{fig:bm_1_stress_xx}. The clustering has been dynamically adapted throughout the ASCA solution procedure, where the number of clusters has increased from $n_{\mathrm{c}}=65$ to $n_{\mathrm{c}}=190$ (see Figure~\ref{subfig:bm_1_adapt_n_clusters}). A total of 5 clustering adaptivity steps have been performed (signaled by the diamond-shaped markers), two of them in the yielding region (particle-matrix interface) and the remaining on the post-yielding region (plastic band). It is possible to see that ASCA solution's accuracy outperforms the SCA `medium' solution with a similar number of clusters as well as the SCA `fine' solution with $n_{\mathrm{c}}=305$, exhibiting an overall relative error below $4\%$. The ASCA solution departs from the SCA `coarse' solution with $n_{\mathrm{c}}=65$ occurs in the yielding region soon after the first clustering adaptivity step, hence highlighting the importance of such procedure to accurately capture the initial yielding of the matrix material phase. In addition, ASCA is able to predict the material fracture (signaled by the star-shaped marker) and toughness with significantly greater accuracy (see Table~\ref{tab:bm_1_failure_toughness}).

\begin{table}[hbt]
\caption{Ductile failure and associated toughness of Benchmark 1, predicted by the different solution methods. Failure assumed when $0.5\%$ of the matrix phase surpasses an accumulated plastic strain of $0.25$. Number of clusters: ($\mathrm{c}$) $n_{c}=65$ (coarse), ($\mathrm{m}$) $n_{c}=185$ (medium), ($\mathrm{f}$) $n_{c}=305$ (fine) and (*) $n_{c}=65 \rightarrow 190$.}
\label{tab:bm_1_failure_toughness}
\centering
\setlength{\tabcolsep}{0.3cm}
\renewcommand{\arraystretch}{1.3}
\begin{tabular}{ccccc}
	 & \multicolumn{2}{c}{\textbf{Failure}} & \multicolumn{2}{c}{\textbf{Toughness}} \\ \cmidrule(l){2-3} \cmidrule(l){4-5}
	 \textbf{Method} & $\varepsilon_{xx}$ & $\sigma_{xx}$(MPa) & $J$(mJ/mm$^{3}$) & $\epsilon (\%)$ \\ \midrule
	 DNS & $3.50\times 10^{-2}$ & $8.65 \times 10^{-1}$ & $2.23 \times 10^{-2}$ & - \\
	 SCA${}^{(\mathrm{c})}$ & $4.68\times 10^{-2}$ & $10.67 \times 10^{-1}$ & $3.62 \times 10^{-2}$ & 62.33 \\
	 SCA${}^{(\mathrm{m})}$ & $4.45\times 10^{-2}$ & $9.89 \times 10^{-1}$ & $3.27 \times 10^{-2}$ & 46.64 \\
	 SCA${}^{(\mathrm{f})}$ & $4.30\times 10^{-2}$ & $8.65 \times 10^{-1}$ & $3.10 \times 10^{-2}$ & 39.01 \\
	 ASCA${}^{(*)}$ & $4.00\times 10^{-2}$ & $9.37 \times 10^{-1}$ & $2.78 \times 10^{-2}$ & 24.66
 \\ \bottomrule 
\end{tabular}
\end{table}

\begin{figure}[hbt]
\centering
		\begin{subfigure}[b]{0.495\textwidth}
                \centering
                \includegraphics[width=\textwidth]{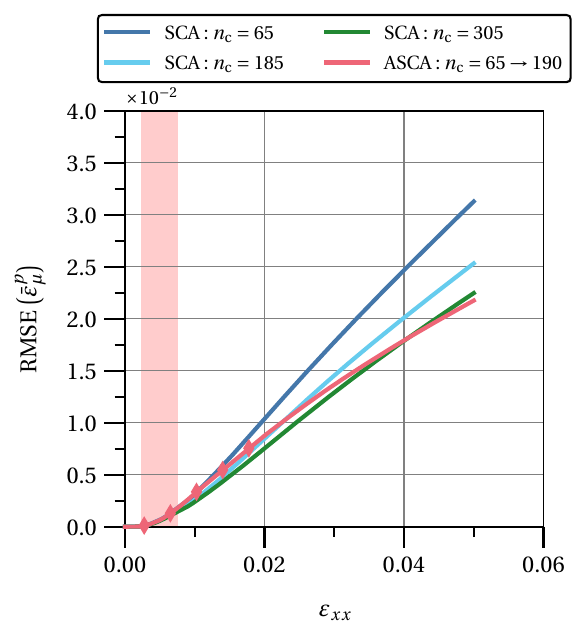}
                \caption{}
                \label{subfig:bm_1_rmse_acc_p_strain}
        \end{subfigure}
        \begin{subfigure}[b]{0.495\textwidth}
                \centering
                \includegraphics[width=\textwidth]{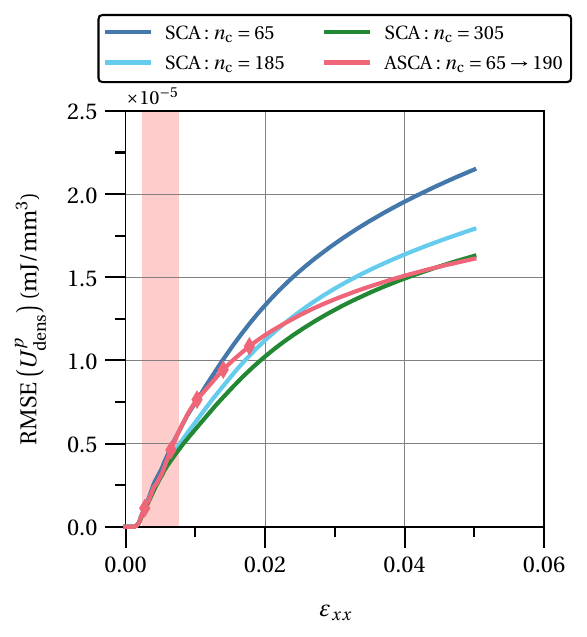}
                \caption{}
                \label{subfig:bm_1_rmse_acc_p_energy_dens}
        \end{subfigure}\hfill
\caption{Root-mean-square error (RMSE) of Benchmark 1 microscale fields solutions relative to the FEM DNS solution: \subref{subfig:bm_1_rmse_acc_p_strain} Local accumulated plastic strain field; \subref{subfig:bm_1_rmse_acc_p_energy_dens} Local accumulated plastic strain energy density field. Symbology: clustering adaptivity step ($\blacklozenge$).}
\label{fig:bm_1_rmse}
\end{figure}

The numerical assessment proceeds with the comparison between the microscale solution fields, as shown by the RMSE of the local accumulated plastic strain and accumulated plastic strain energy density fields in Figure~\ref{fig:bm_1_rmse}. Note that such error metric accounts for the errors stemming from the whole domain, hence the improvement from the SCA solution with $n_{\mathrm{c}}=65$ is delayed in comparison with macroscale homogenized response. Nonetheless, the ASCA solution's accuracy surpasses the SCA `medium' solution with a similar number of clusters and even reaches the SCA `fine' solution with $n_{\mathrm{c}}=305$. These and the macroscale homogenized results can be further comprehended by analyzing the clustering and local accumulated plastic strain field at the end of the deformation path (see Figure~\ref{fig:bm_1_local_acc_p_strain}). While the SCA `fine' solution with $n_{\mathrm{c}}=305$ scatters the clusters over the whole domain, the ASCA solution places the majority of clusters where they are most need. As the magnitude of the plastic strain localization increases, the improved accuracy of the ASCA on these regions dominates over the less significant errors in the remaining domain, in accordance with Figure~\ref{fig:bm_1_rmse}.

Despite the encouraging results from an accuracy point of view, the efficiency of the ACROM framework is fundamental in terms of engineering utility. The partial and total computational times of the different solution methods are presented in Table~\ref{tab:bm_1_comp_times} (see Figure~\ref{fig:scaschemecompact} for a description of each step). In the first place, the DNS solution computational time is approximately one order of magnitude greater than the most expensive SCA solution with $n_{c}=305$. Second, the total time of the ASCA solution is faster than the SCA `medium' solution with a similar number of clusters. It is instructive to analyze the partial computational times to understand the comparison between the SCA and ASCA. The computational time associated with the DNS linear elastic analyses (step 1 of the offline-stage) is independent of the number of clusters, hence similar among the different solutions. The computational cost of the base cluster analysis (step 2 of the offline-stage) scales with the number of clusters. Given that ASCA starts with an initial number of clusters $n_{\mathrm{c}}=65$, this step is faster than both SCA solutions with a higher number of clusters. The same applies to the computation of the cluster interaction tensors (step 3 of the offline-stage), the most expensive step of the offline-stage. It is thus clear that ASCA results in significant savings associated with the offline-stage overhead costs. Concerning the actual solution procedure of the online-stage, ASCA is much faster than the SCA with a similar number of clusters. This results from the fact that approximately one-fourth of the deformation path is solved with $n_{c}<185$ (see Figure~\ref{fig:bm_1_adaptivity_metrics}). Finally, there is a computational cost associated with all clustering adaptivity procedures exclusive to the ASCA solution. It is noticeable that this cost is relatively high when compared with the actual online solution procedure. Figure~\ref{subfig:bm_1_adapt_phases_time_perc} shows that the update of the cluster interaction tensors (Block C) dominates the computational cost of the clustering adaptivity after the first clustering adaptivity step is performed.

\begin{table}[hbt]
\caption{Partial and total computational times (s) of the different solution methods associated to Benchmark 1. Number of clusters: ($\mathrm{c}$) $n_{c}=65$ (coarse), ($\mathrm{m}$) $n_{c}=185$ (medium), ($\mathrm{f}$) $n_{c}=305$ (fine) and (*) $n_{c}=65 \rightarrow 190$. Offline-stage: linear elastic DNS simulations (step 1), base cluster analysis (step 2) and computation of cluster-interaction tensors (step 3).}
\label{tab:bm_1_comp_times}
\centering
\setlength{\tabcolsep}{0.25cm}
\renewcommand{\arraystretch}{1.3}
\begin{tabular}{ccccccc}
	 & \multicolumn{6}{c}{\textbf{Computational time (s)}} \\ \cmidrule(l){2-7}
	 \multirow{2}{*}{\textbf{Method}} & \textbf{Offline} & \textbf{Offline} & \textbf{Offline} & \textbf{Online} & \textbf{Online} & \multirow{2}{*}{\textbf{Total}} \\[-4pt]
	 & \textbf{(step 1)} & \textbf{(step 2)} & \textbf{(step 3)} & \textbf{(solution)} & \textbf{(adapt.)} & \\ \midrule
	 DNS & - & - & - & $2.38 \times 10^{4}$ & - & $2.38 \times 10^{4}$ \\
	 SCA${}^{(\mathrm{c})}$ & $2.19 \times 10^{2}$ & $2.20 \times 10^{1}$ & $5.20 \times 10^{1}$ & $1.32 \times 10^{2}$ & - & $4.25 \times 10^{2}$ \\
	 SCA${}^{(\mathrm{m})}$ & $2.19 \times 10^{2}$ & $1.08 \times 10^{2}$ & $3.99 \times 10^{2}$ & $1.35 \times 10^{3}$ & - & $2.08 \times 10^{3}$ \\
	 SCA${}^{(\mathrm{f})}$ & $2.19 \times 10^{2}$ & $1.54 \times 10^{2}$ & $1.10 \times 10^{3}$ & $2.46 \times 10^{3}$ & - & $3.93 \times 10^{3}$ \\
	 ASCA${}^{(*)}$ & $2.19 \times 10^{2}$ & $2.20 \times 10^{1}$ & $5.20 \times 10^{1}$ & $3.24 \times 10^{2}$ & $6.37 \times 10^{2}$ & $1.01 \times 10^{3}$ \\
	  \bottomrule 
\end{tabular}
\end{table}

\paragraph{Remark} Despite promoting different plastic strain localization patterns, the same clustering adaptivity strategy has been successfully employed to solve benchmarks 2 and 3. Given that similar conclusions have been obtained, the corresponding results and discussion are omitted to avoid overextending this paper without any apparent benefit.

\clearpage

\subsection{Particle-matrix RVE analysis}
The analysis of the particle-matrix RVE is far more complex than the previously described benchmarks. Figure~\ref{fig:2d_frc_acc_p_strain_fracture} shows that numerous plastic strain localized bands are developed between particles throughout the deformation path. It is possible to conclude that the composite fracture results from two proeminent plastic bands whose propagation occurs transversally to the loading direction. It is thus comprehensible that the clustering adaptivity strategy employed in the previous section may not be the most suitable choice to tackle this problem. Moreover, some of the additional adaptivity procedures proposed in Section~\ref{ssec:addadaptproc} can be explored here as well as the adaptive clustering solution rewinding procedure described in Section~\ref{ssec:rewindproc}.

The identification of the accumulated plastic strain field's main evolution steps is not straightforward as in the benchmark cases. At most, it can be a priori expected that the composite macroscale yielding and posterior fracture is dictated by one or more continuous bands where plastic strain localizes with a higher magnitude. The `optimal' clustering should thus prioritize the refinement of the primary plastic bands but cannot disregard the development of secondary localized phenomena with lower magnitude.

\begin{figure}[hbt]
\centering
	\includegraphics[width=0.9\textwidth]{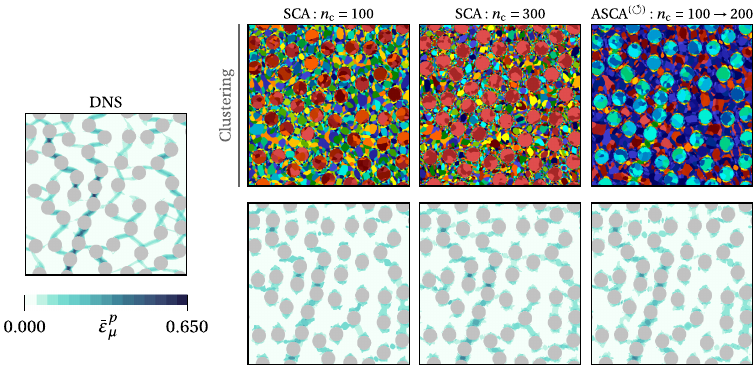}
\caption{Comparison of the clustering-based domain decomposition and local accumulated plastic strain field at the end of the deformation path of the particle-matrix composite 2D RVE. Symbology: solution rewinding ($\circlearrowleft$).}
\label{fig:2d_frc_local_acc_p_strain}
\end{figure}

The initial number of clusters is set to $n_{c}=100$ (70 clusters in the matrix ACRMP and 30 clusters in the particles SCRMP). Figure~\ref{fig:2d_frc_local_acc_p_strain} shows that the corresponding SCA solution satisfies a minimal required accuracy despite failing to capture several secondary plastic bands.

Concerning the ACROM framework, the clustering adaptivity feature for the matrix ACRMP is again set as the accumulated plastic strain, $c_{\mu}(\bm{y})=\bar{\varepsilon}^{p}$. Given that several plastic strain localized bands are not developed at the same time and/or rate, attempting to follow all plastic strain fronts efficiently is cumbersome. Instead, the strategy adopted here is to perform a single clustering adaptivity step right after the yielding region ($\bar{\varepsilon}^{p}=0.05$), i.e., once the primary plastic bands have developed as well as some of the secondary ones. In order to do so, it is assumed that a 5\% normalized accumulated plastic strain discontinuity is significant by setting $\gamma_{\mathrm{ratio}}=0.05$. In addition, the number of child clusters is set dynamically as $n_{\mathrm{c}, \mathrm{child}} \in [3,8]$ ($f_{\mathrm{c}, \mathrm{child}}=0.125$, $\gamma_{\mathrm{split}}=0.7$, $\gamma_{\mathrm{split}}^{\Delta}=0.6$) and assuming the default linear magnitude function ($n=1.0$). In this way, the regions associated with the primary plastic bands (higher discontinuities) are refined with a greater number of clusters without wasting too many clusters in the secondary plastic bands (lower discontinuities). For the same reason, $s_{\mathrm{split}}^{\mathrm{low}}=0.0$ ($\theta=0.0$) so that the matrix lower-valued clusters surrounding all plastic bands do not exhaust the available number of clusters.

\paragraph{Remark} Although the refined clustering of the primary plastic bands is crucial to accurately predict the material's structural degradation and fracture, the suitable refinement of the surrounding matrix clustering plays an important role in the macroscale homogenized response.

Besides the DNS solution, three different SCA solutions are again selected for comparison following the same reasoning: a `coarse' solution with $n_{\mathrm{c}}=100$, a `medium' solution with $n_{\mathrm{c}}=200$ and a `fine' solution with $n_{\mathrm{c}}=300$. The number of clusters of the particles SCRMP\footnote{It is remarked that if the number of clusters of the particles' SCRMP is increased proportionally to the corresponding volume fraction, then the convergence of the SCA and ASCA solutions towards the DNS is significantly improved.} is always set to 30. To evaluate the performance of the ASCA solution, the initial number of clusters is set to $n_{\mathrm{c}}=100$, the threshold number of clusters is set to $n_{\mathrm{c}}=200$ and the maximum allowed adaptivity level is kept as $\gamma_{\mathrm{max}}=2$. Moreover, following the discussion in Section~\ref{ssec:rewindproc}, two ASCA solutions are evaluated: one where the ASCA solution proceeds after the clustering adaptivity step is performed, and the other where the ASCA solution is rewound back to the beginning of the deformation path after the clustering adaptivity step. The following results demonstrate that solving the yielding region with a proper clustering refinement and overwriting the less accurate state variable history is fundamental to the solution's accuracy. This approach can then be characterized by two steps: (1) a first fast step where a coarse clustering is solved up to the point where a proper clustering adaptivity has been performed and (2) a second step where the deformation path and material state history is partially or totally solved with adequate domain decomposition. Unless otherwise stated, the ASCA solution with rewinding is evaluated in the following discussion.

\begin{figure}[hbt]
\centering
		\begin{subfigure}[b]{0.495\textwidth}
                \centering
                \includegraphics[width=\textwidth]{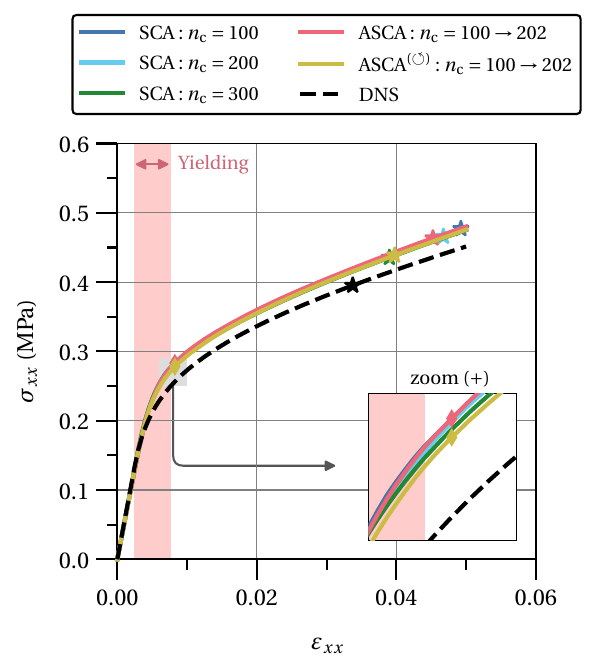}
                \caption{}
                \label{subfig:2d_frc_stress_xx_strain_xx}
        \end{subfigure}
        \begin{subfigure}[b]{0.495\textwidth}
                \centering
                \includegraphics[width=\textwidth]{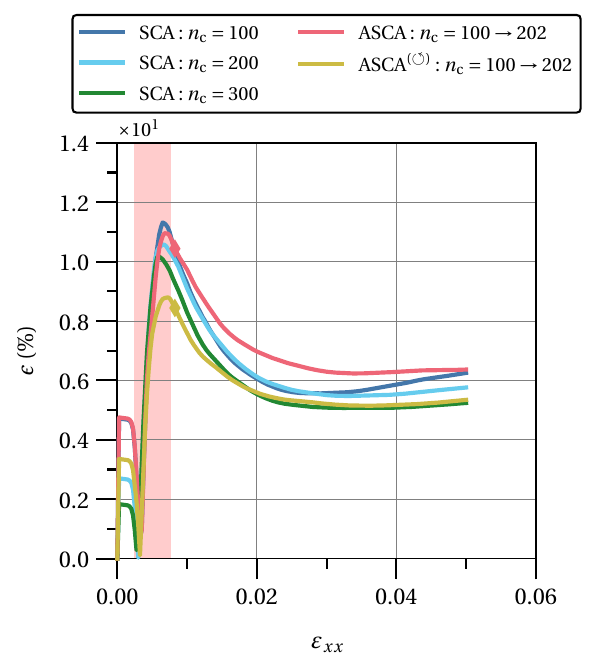}
                \caption{}
                \label{subfig:2d_frc_stress_xx_strain_xx_error}
        \end{subfigure}\hfill
\caption{Comparison of the macroscale homogenized response of the particle-matrix composite 2D RVE predicted by different solution methods: \subref{subfig:2d_frc_stress_xx_strain_xx} Homogenized stress; \subref{subfig:2d_frc_stress_xx_strain_xx_error} Relative error with respect to the FEM DNS solution. Symbology: clustering adaptivity step ($\blacklozenge$), fracture criterion prediction ($\bigstar$), solution rewinding ($\circlearrowleft$).}
\label{fig:2d_frc_stress_xx}
\end{figure}

\begin{figure}[hbt]
\centering
		\begin{subfigure}[b]{0.494\textwidth}
                \centering
                \includegraphics[width=\textwidth]{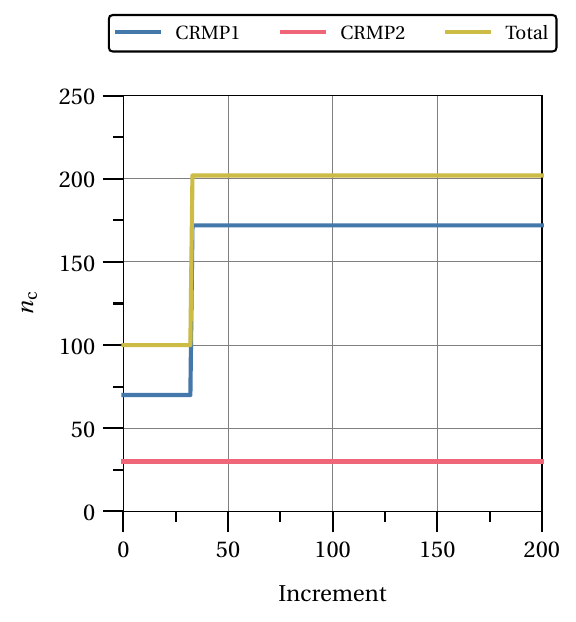}
                \caption{}
                \label{subfig:2d_frc_adapt_n_clusters}
        \end{subfigure}
        \begin{subfigure}[b]{0.498\textwidth}
                \centering
                \includegraphics[width=\textwidth]{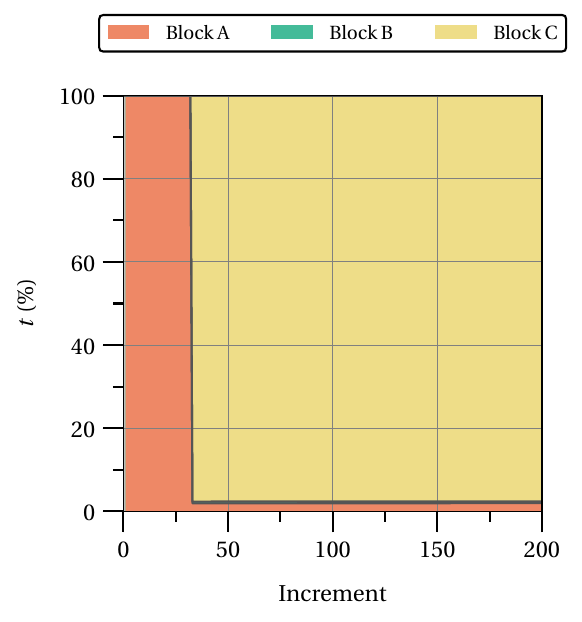}
                \caption{}
                \label{subfig:2d_frc_adapt_phases_time_perc}
        \end{subfigure}\hfill
\caption{Clustering adaptivity of particle-matrix composite 2D RVE (without rewinding): \subref{subfig:2d_frc_adapt_n_clusters} Evolution of the number of clusters; \subref{subfig:2d_frc_adapt_phases_time_perc} Evolution of the cluster adaptivity framework's blocks computation time relative to the total adaptivity computation time up to a given increment.}
\label{fig:2d_frc_adaptivity_evolution}
\end{figure}

The comparison between the macroscale homogenized response predicted by the different solution methods is shown in Figure~\ref{fig:2d_frc_stress_xx}. The clustering has been dynamically adapted from $n_{\mathrm{c}}=100$ to $n_{\mathrm{c}}=202$ in a single clustering adaptivity step (signaled by the diamond-shaped marker) right after the yielding region (see Figure~\ref{fig:2d_frc_adaptivity_evolution}). It is possible to see that the ASCA solution's accuracy outperforms the SCA solution with a similar number of clusters, namely in the yielding region where the divergence between both solutions occurs. The outcome of a clustering adaptivity step performed right after the yielding region is even more highlighted in comparison with the SCA `fine' solution with $n_{\mathrm{c}}=300$. The ASCA solution has greater accuracy in the yielding and early post-yielding region, after which the SCA `fine' solution catches up due to a more refined clustering dispersed over the whole matrix material phase. A similar comparison results from the prediction of the material fracture (signaled by the star-shaped marker) and toughness (see Table~\ref{tab:2d_frc_failure_toughness}). ASCA prediction's accuracy is almost equal to the SCA's `fine' solution, while significantly better than the SCA `medium' solution with a similar number of clusters. Despite the reasonable fracture and toughness predictions of the ASCA solution without rewinding, the accuracy of the macroscale homogenized response is severely affected and significantly worse than SCA `medium' solution with a similar number of clusters.

\begin{table}[hbt]
\caption{Ductile failure and associated toughness of the particle-matrix composite 2D RVE predicted by the different solution methods. Failure assumed when $0.5\%$ of the matrix phase surpasses an accumulated plastic strain of $0.25$. Number of clusters: ($\mathrm{c}$) $n_{c}=100$ (coarse), ($\mathrm{m}$) $n_{c}=200$ (medium), ($\mathrm{f}$) $n_{c}=300$ (fine), (*) $n_{c}=100 \rightarrow 202$ and ($\circlearrowleft$) $n_{c}=100 \rightarrow 202$ (with rewinding).}
\label{tab:2d_frc_failure_toughness}
\centering
\setlength{\tabcolsep}{0.3cm}
\renewcommand{\arraystretch}{1.3}
\begin{tabular}{ccccc}
	 & \multicolumn{2}{c}{\textbf{Failure}} & \multicolumn{2}{c}{\textbf{Toughness}} \\ \cmidrule(l){2-3} \cmidrule(l){4-5}
	 \textbf{Method} & $\varepsilon_{xx}$ & $\sigma_{xx}$(MPa) & $J$(mJ/mm$^{3}$) & $\epsilon (\%)$ \\ \midrule
	 DNS & $3.38\times 10^{-2}$ & $3.95 \times 10^{-1}$ & $1.00 \times 10^{-2}$ & - \\
	 SCA${}^{(\mathrm{c})}$ & $4.93\times 10^{-2}$ & $4.78 \times 10^{-1}$ & $1.76 \times 10^{-2}$ & 76.00 \\
	 SCA${}^{(\mathrm{m})}$ & $4.68\times 10^{-2}$ & $4.66 \times 10^{-1}$ & $1.64 \times 10^{-2}$ & 64.00 \\
	 SCA${}^{(\mathrm{f})}$ & $3.90\times 10^{-2}$ & $4.35 \times 10^{-1}$ & $1.28 \times 10^{-2}$ & 28.00 \\
	 ASCA${}^{(*)}$ & $4.53\times 10^{-2}$ & $4.64 \times 10^{-1}$ & $1.58 \times 10^{-2}$ & 58.00 \\
	 ASCA${}^{(\circlearrowleft)}$ & $3.98\times 10^{-2}$ & $4.39 \times 10^{-1}$ & $1.31 \times 10^{-2}$ & 31.00
 \\ \bottomrule 
\end{tabular}
\end{table}

\begin{figure}[hbt]
\centering
		\begin{subfigure}[b]{0.495\textwidth}
                \centering
                \includegraphics[width=\textwidth]{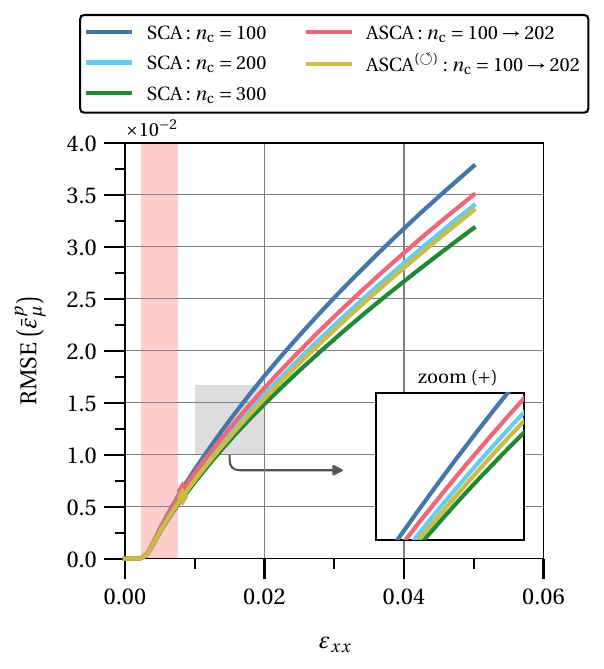}
                \caption{}
                \label{subfig:2d_frc_rmse_acc_p_strain}
        \end{subfigure}
        \begin{subfigure}[b]{0.480\textwidth}
                \centering
                \includegraphics[width=\textwidth]{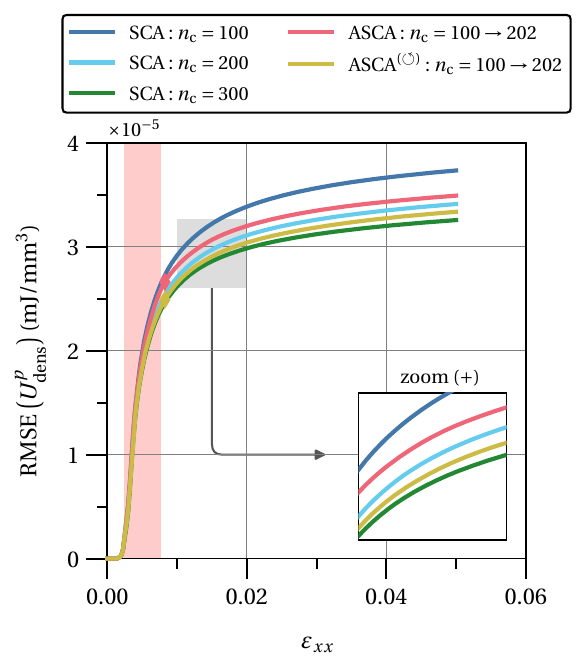}
                \caption{}
                \label{subfig:2d_frc_rmse_acc_p_energy_dens}
        \end{subfigure}\hfill
\caption{Root-mean-square error (RMSE) of the particle-matrix composite 2D RVE microscale fields solutions relative to the FEM DNS solution: \subref{subfig:2d_frc_rmse_acc_p_strain} Local accumulated plastic strain field; \subref{subfig:2d_frc_rmse_acc_p_energy_dens} Local accumulated plastic strain energy density field. Symbology: clustering adaptivity step ($\blacklozenge$), solution rewinding ($\circlearrowleft$).}
\label{fig:2d_frc_rmse}
\end{figure}

The RMSE of the local accumulated plastic strain and accumulated plastic strain energy density fields is shown in Figure~\ref{fig:2d_frc_rmse}. The ASCA solution's accuracy surpasses the SCA `medium' solution with a similar number of clusters and closely follows the SCA `fine' solution in the yielding and early post-yielding regions. These results follow the same reasoning of the macroscale homogenized response and are complemented by the clustering and local accumulated plastic strain field at the end of the deformation path (see Figure~\ref{fig:2d_frc_local_acc_p_strain}). It is important to highlight the contrast between the ASCA solution with and without performing the rewinding procedure. While the former diverges from the SCA `coarse' solution after the clustering adaptivity step, the latter diverges from the beginning and improves the solution's accuracy throughout the whole deformation path.

The partial and total computational times of the different methods are presented in Table~\ref{tab:2d_frc_comp_times}. In the first place, the ASCA solution is much faster than the DNS solution concerning the total computational time, faster than the SCA `fine' solution with $n_{c}=300$ and slightly slower than the SCA `medium' solution with the same number of clusters. The later comparison results from the balance between different costs: (i) ASCA base cluster analysis is faster due to the lower number of initial clusters; (ii) the computation of ASCA base clustering cluster interaction tensors is much faster for the same reason; (iii) ASCA must re-solve the macroscale loading increments prior to the solution rewinding; (iv) ASCA must spend additional computational time associated to the clustering adaptivity procedures. Hence, by rewinding the solution back to the beginning of the deformation path, the ASCA solution not only loses the benefit of partially solving the deformation path with a lower number of clusters but also needs to re-solve the overwritten loading increments. This is a clear trade-off between accuracy and efficiency, which is perfectly reasonable in this particular problem. At last, Figure~\ref{subfig:2d_frc_adapt_phases_time_perc} shows once again that the update of the cluster interaction tensors (Block C) dominates the computational cost of the clustering adaptivity.

\begin{table}[hbt]
\caption{Partial and total computational times (s) of the different solution methods associated to the particle-matrix 2D RVE. Number of clusters: ($\mathrm{c}$) $n_{c}=100$ (coarse), ($\mathrm{m}$) $n_{c}=200$ (medium), ($\mathrm{f}$) $n_{c}=300$ (fine), (*) $n_{c}=100 \rightarrow 202$ and ($\circlearrowleft$) $n_{c}=100 \rightarrow 202$ (with rewinding). Offline-stage: linear elastic DNS simulations (step 1), base cluster analysis (step 2) and computation of cluster-interaction tensors (step 3).}
\label{tab:2d_frc_comp_times}
\centering
\setlength{\tabcolsep}{0.25cm}
\renewcommand{\arraystretch}{1.3}
\begin{tabular}{ccccccc}
	 & \multicolumn{6}{c}{\textbf{Computational time (s)}} \\ \cmidrule(l){2-7}
	 \multirow{2}{*}{\textbf{Method}} & \textbf{Offline} & \textbf{Offline} & \textbf{Offline} & \textbf{Online} & \textbf{Online} & \multirow{2}{*}{\textbf{Total}} \\[-4pt]
	 & \textbf{(step 1)} & \textbf{(step 2)} & \textbf{(step 3)} & \textbf{(solution)} & \textbf{(adapt.)} & \\ \midrule
	 DNS & - & - & - & $3.08 \times 10^{4}$ & - & $3.08 \times 10^{4}$ \\
	 SCA${}^{(\mathrm{c})}$ & $2.19 \times 10^{2}$ & $3.44 \times 10^{1}$ & $1.15 \times 10^{2}$ & $1.84 \times 10^{2}$ & - & $5.52 \times 10^{2}$ \\
	 SCA${}^{(\mathrm{m})}$ & $2.19 \times 10^{2}$ & $1.45 \times 10^{2}$ & $4.38 \times 10^{2}$ & $3.97 \times 10^{2}$ & - & $1.20 \times 10^{3}$ \\
	 SCA${}^{(\mathrm{f})}$ & $2.19 \times 10^{2}$ & $2.16 \times 10^{2}$ & $9.69 \times 10^{2}$ & $6.38 \times 10^{2}$ & - & $2.04 \times 10^{3}$ \\
	 ASCA${}^{(*)}$ & $2.19 \times 10^{2}$ & $3.44 \times 10^{1}$ & $1.15 \times 10^{2}$ & $3.67 \times 10^{2}$ & $4.32 \times 10^{2}$ & $1.17 \times 10^{3}$ \\
	 ASCA${}^{(\circlearrowleft)}$ & $2.19 \times 10^{2}$ & $3.44 \times 10^{1}$ & $1.15 \times 10^{2}$ & $4.75 \times 10^{2}$ & $4.32 \times 10^{2}$ & $1.28 \times 10^{3}$ \\
	  \bottomrule 
\end{tabular}
\end{table}

\clearpage

\subsection{Update of the cluster interaction matrix}
The previous results show that the computation of the cluster interaction tensors (Block C) plays a major role in the total computational cost of the clustering adaptivity procedures. In this regard, two important aspects have been discussed in Section~\ref{ssec:blockC}: (i) taking advantage of the cluster-symmetry property of cluster interaction tensors (see Equation~\eqref{eq:cit_cluster_symmetry}) and (ii) update the cluster interaction matrix such that only new clusters require a full computation of the cluster interaction tensors. It is important to remark that the computational cost of the cluster interaction tensors scales with the number of problem dimensions, the total number of voxels and the number of independent strain components. Following the numerical examples mostly explored in the previous sections, it is assumed a 2D problem, a spatial discretization in a regular grid of $n_{v} = 400 \times 400$ voxels and an infinitesimal strain formulation with 3 independent strain components. 

\begin{figure}[hbt]
\centering
	\includegraphics[width=0.5\textwidth]{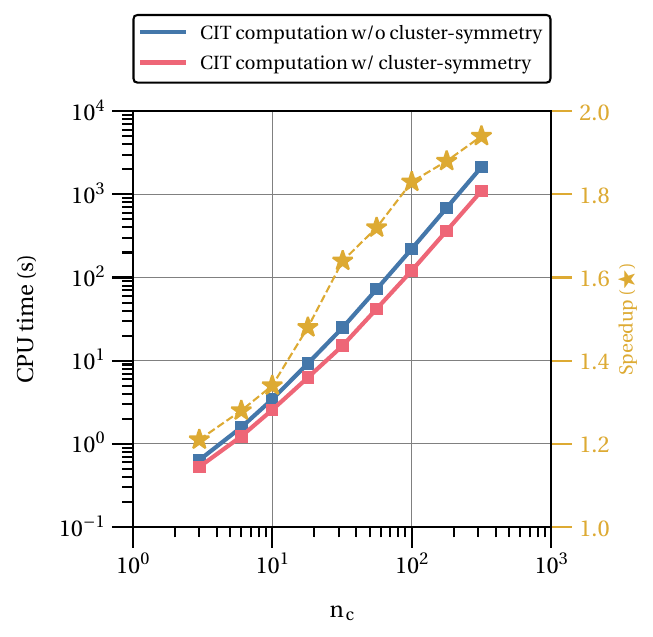}
\caption{Comparison between the computational cost of the cluster interaction tensors with and without taking advantage of the cluster-symmetry property. Conditions: 2D problem, $n_{v}=400\times400$, infinitesimal strain formulation with 3 independent strain components.}
\label{fig:CIT_time_full_vs_sym_2d}
\end{figure}

To illustrate the advantage of the cluster-symmetry property, the computation of all cluster interaction tensors is performed for a different number of clusters comprised between 3 and 316. Figure~\ref{fig:CIT_time_full_vs_sym_2d} compares the computational time with and without considering the cluster-symmetry property as well as the associated speed-up. As expected, the cluster-symmetry property accelerates the computation and the speed-up increases significantly with the number of clusters. This scaling is a direct consequence of the increase in the number of cluster interaction tensors that are obtained directly through cluster-symmetry.

The analysis of the computational cost associated with the update of the cluster interaction matrix is more intricate. For a given clustering adaptivity step, this cost depends on the initial number of clusters, $n_{\mathrm{c}}^{\mathrm{init}}$, the number of unchanged clusters, $n_{\mathrm{c}}^{\mathrm{old}}$, and the number of new clusters, $n_{\mathrm{c}}^{\mathrm{new}}$. Assume that $n_{\mathrm{c}}^{\mathrm{old}}$ and $n_{\mathrm{c}}^{\mathrm{new}}$ are both defined as functions of $n_{\mathrm{c}}^{\mathrm{init}}$ as
\begin{equation}
    n_{\mathrm{c}}^{\mathrm{old}} = \mathrm{nint}( \alpha \,  n_{\mathrm{c}}^{\mathrm{init}}) \, ,
\end{equation}
and
\begin{equation}
    n_{\mathrm{c}}^{\mathrm{new}} = \mathrm{nint}( (1 + \beta) \, n_{\mathrm{c}}^{\mathrm{init}}) - n_{\mathrm{c}}^{\mathrm{old}} \, ,
\end{equation}
where $\alpha \in [0, 1]$ and $\beta \in [0, \infty[$. From these definitions, the parameter $\alpha$ sets the relative amount of initial clusters that are kept unchanged and parameter $\beta$ sets the relative increase of the number of clusters resulting from the clustering adaptivity step. For the sake of clarity, assume that at the beginning of a clustering adaptivity step, one has $n_{\mathrm{c}}^{\mathrm{init}} = 10$ and $(\alpha, \beta) = (0.5, 1.0)$. This means that $n_{\mathrm{c}}^{\mathrm{old}} = 5$ clusters are not adapted and the remaining 5 clusters are targeted, originating $n_{\mathrm{c}}^{\mathrm{new}} = 10$ new clusters and leading to a total of 20 clusters.

\begin{figure}[hbt]
\centering
		\begin{subfigure}[b]{0.495\textwidth}
                \centering
                \includegraphics[width=\textwidth]{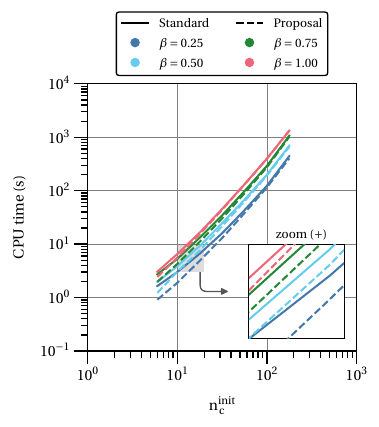}
                \caption{}
                \label{subfig:CIT_time_naive_vs_proposal_2d_inkscape}
        \end{subfigure}
        \begin{subfigure}[b]{0.48\textwidth}
                \centering
                \includegraphics[width=\textwidth]{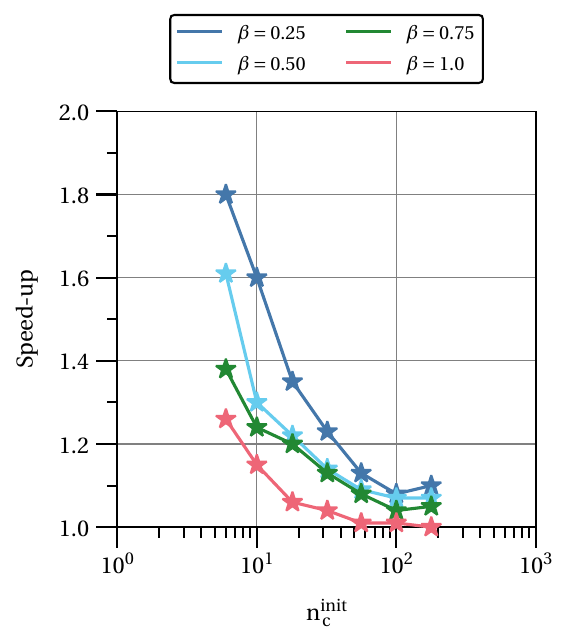}
                \caption{}
                \label{subfig:CIT_time_naive_vs_proposal_2d_speed_up}
        \end{subfigure}\hfill
\caption{Computational cost associated with the update of the cluster interaction matrix on a given clustering adaptivity step ($\alpha=0.75$): \subref{subfig:CIT_time_naive_vs_proposal_2d_inkscape} Comparison between standard and proposed approaches; \subref{subfig:CIT_time_naive_vs_proposal_2d_speed_up} Speed-up resulting from the proposed approach. Conditions: 2D problem, $n_{v}=400\times400$, infinitesimal strain formulation with 3 independent strain components.}
\label{fig:CIT_time_naive_vs_proposal_2d}
\end{figure}

To illustrate the advantage of the procedure proposed in Section~\ref{ssec:blockC}, Figure~\ref{fig:CIT_time_naive_vs_proposal_2d} shows the comparison of the computational cost associated with the update of the cluster interaction matrix on a given clustering adaptivity step. The standard procedure involves cycling through the cluster interaction matrix rows and columns in a sorted manner, computing the clustering interaction tensors and taking advantage of the cluster-symmetry property whenever possible. It is assumed that $\alpha=0.75$, i.e., $n_{\mathrm{c}}^{\mathrm{old}} = \mathrm{nint}(0.75  \, n_{\mathrm{c}}^{\mathrm{init}})$ are kept unchanged, and the relative increase of the number of clusters is varied according to $\beta \in [0.25, 1.0]$.

In the first place, it is possible to observe that the proposed strategy is more efficient than the standard approach in the whole spectrum of the number of initial clusters, $n_{\mathrm{c}}^{\mathrm{init}}$. In the second place, it is seen that the speed-up decreases as the number of initial clusters increases. This is expected given that the computational savings associated with the convolution of previously existing clusters interaction tensors becomes less relevant as the total number of computed tensors increases. Finally, the speed-up also decreases as the relative number of new clusters increases and the full computation of the associated cluster interaction tensors becomes dominant.

\begin{figure}[hbt]
\centering
	\includegraphics[width=0.5\textwidth]{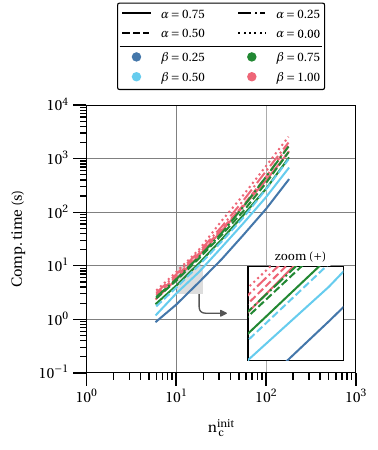}
\caption{Computational cost associated with the update of the cluster interaction matrix on a given clustering adaptivity step. Conditions: 2D problem, $n_{v}=400\times400$, infinitesimal strain formulation with 3 independent strain components.}
\label{fig:CIT_time_update_2d_inkscape}
\end{figure}

At last, Figure~\ref{fig:CIT_time_update_2d_inkscape} compiles the computational cost associated with the update of the cluster interaction matrix for $n_{\mathrm{c}}^{\mathrm{init}} \in [6, 178]$, $\alpha \in [0, 0.75]$ and $\beta \in [0.25, 1.0]$\footnote{It is possible to verify that the constraint $\beta >= 1 - \alpha$ applies irrespective of the number of initial clusters.}. 

As expected, the computational cost increases significantly with the increase of the number of initial clusters, i.e., the update of the cluster interaction matrix tends to become more expensive as the clustering adaptivity advances. For a given number of initial clusters, $n_{\mathrm{c}}^{\mathrm{init}}$, and a given number of unchanged clusters (fixed $\alpha$), the computational cost increases with the number of new clusters as more cluster interaction tensors must be computed. For a given number of initial clusters, $n_{\mathrm{c}}^{\mathrm{init}}$, and a given number of new clusters (fixed $\beta$), the computational cost increases as the number of unchanged clusters decreases. Maximum computational cost occurs in the particular case where all clusters are adapted ($\alpha = 0$, $\beta > 0$), i.e., a whole new cluster interaction matrix must be computed.

The main observations and tendencies taken from the results in this section  should hold irrespective of the problem dimension, the spatial discretization and strain formulation. Nonetheless, the absolute value of the computational cost is expected to increase as the spatial discretization is refined (an increase of the total number of voxels) and the number of independent strain components increases (3D problem and/or finite strain formulation). Similar studies should be carried out in future contributions to investigate the actual computational costs stemming from these conditions.

\clearpage

\section{Conclusion \label{sec:conclusion}}
This paper proposes a novel Adaptive Clustering-based Reduced-Order Modeling (ACROM) framework to significantly improve the recent family of clustering-based reduced-order models (CROMs). In contrast with the static nature of most CROMs, where the clustering-based domain decomposition is solely defined by a preliminary offline-stage, this framework allows the clustering to evolve dynamically throughout the actual problem solution. This adaptive enhancement offers a new route to fast and accurate material modeling involving highly localized plasticity and damage phenomena. The framework is suitable for infinitesimal or finite strain-based formulations and is composed of three main blocks: (i) target clusters selection criterion (Block A), to determine what clusters need to be adapted; (ii) adaptive cluster analysis (Block B), to perform the actual clustering adaptivity; and (iii) computation of cluster interaction tensors (Block C), required to complete the clustering characterization. Several complementary procedures are further proposed to improve the overall performance of the adaptivity process, namely an adaptive clustering solution rewinding procedure and a dynamic adaptivity split factor strategy.

To illustrate the ACROM framework, the pioneering CROM called Self-Consistent Clustering Analysis (SCA) \cite{liu2016} is extended to its adaptive counterpart coined Adaptive Self-Consistent Clustering Analysis (ASCA). An elasto-plastic particle-matrix composite exhibiting highly localized plasticity is numerically modeled and a simple ductile fracture criterion is defined. Two different clustering adaptivity strategies are successfully demonstrated and a thorough analysis shows that ASCA's accuracy is able to outperform SCA with a similar number of clusters concerning (i) the macroscale homogenized response, (ii) the microscale accumulated plastic strain and accumulated plastic strain energy density fields and (iii) the composite's fracture and toughness prediction. From an efficiency point of view, ASCA can deliver solutions with greater accuracy than SCA with a lower or similar computational cost.

Given the promising results shown in this paper, the proposed framework sets the stage and opens up new avenues to explore adaptivity in the context of CROMs. From a development standpoint, improved selection criteria shall be designed together with appropriate error estimators (Block A), the capabilities of different cluster analysis algorithms may be investigated (Block B), and additional strategies should be employed to accelerate the computation of the cluster interaction tensors (Block C) (e.g., \cite{zhang2019a}). In addition, the effectiveness of the ACROM framework shall be analyzed when dealing with different scales, finite strains, complex materials in terms of topology and/or constitutive behavior (e.g., anisotropy, damage models), non-monotonic loading paths, and failure criteria.

\section{Acknowledgements}

Bernardo P. Ferreira acknowledges the support provided by Funda\c{c}\~{a}o para a Ci\^{e}ncia e a Tecnologia through the scholarship with reference SFRH/BD/130593/2017. This research has also been supported by Instituto de Ci\^{e}ncia e Inova\c{c}\~{a}o em Engenharia Mec\^{a}nica e Engenharia Industrial (INEGI). Miguel A. Bessa acknowledges the support from the project `Artificial intelligence towards a sustainable future: ecodesign of recycled polymers and composites' (with project number 17260 of the research programme Applied and Engineering Sciences) which is financed by the Dutch Research Council (NWO).

\appendix

\section{Key ingredients in Adaptive Finite Element Methods (AFEMs) \label{app:AFEMs}}

As mentioned in the main text, AFEMs have three main ingredients \cite{huerta1998}: (1) an error estimator or indicator; (2) a procedure to adapt the spatial interpolation; and (3) a remeshing criterion.

\subsection{Error estimators and indicators}

Both error estimators and error indicators consist of measures that may be used to assess the error\footnote{In this paper, the focus is given to the spatial discretization error associated with the computational domain decomposition. It is assumed that the remaining sources of error (e.g., time integration, iterative solution methods) are sufficiently small and deemed secondary.}. The distinction between the nature of these quantities is described in \cite{huerta1999} together with a comparison of the associated advantages and limitations. On the one hand, error estimators approximate a measure of the actual error in a given norm or behave as equivalent norms. On the other hand, error indicators are quantities that are chosen, in an ad hoc manner, as an indicator of the error, often based on heuristic considerations. The so-called a posteriori error estimators, which extract information about the error of the approximate solution without invoking the often unknown exact solution, may be classified into two prominent families: recovery-based estimators and residual-based estimators. Recovery-based estimators take advantage of an improved solution by a posterior treatment of the finite element data to estimate the error (e.g., \cite{zienkiewicz1987,zienkiewicz1992,boroomand1997}). In contrast, residual-based estimators make use of the equilibrium residuals of the finite element approximation, either explicitly or implicitly. The simplest explicit residual-based estimators are based on the equilibrium residuals in the element interior and on the flux discontinuities at the element boundary, the later often named jump discontinuities (e.g., \cite{babuska1978,kelly1983,verfurth1994,ainsworth2000}). Among all the existing implicit residual-based estimators, the equilibrated residual estimator (e.g., \cite{bank1985,ainsworth1993}) is the most robust. Other types of estimators involve the analysis of constitutive functionals (e.g., \cite{ladeveze1985,ladeveze1989,gallimard1996}) or general error functionals derived by duality arguments (e.g., \cite{eriksson1988,johnson1992,rannacher1998}).

\begin{figure}[hbt]
	\centering
	\includegraphics[width=0.9\textwidth]{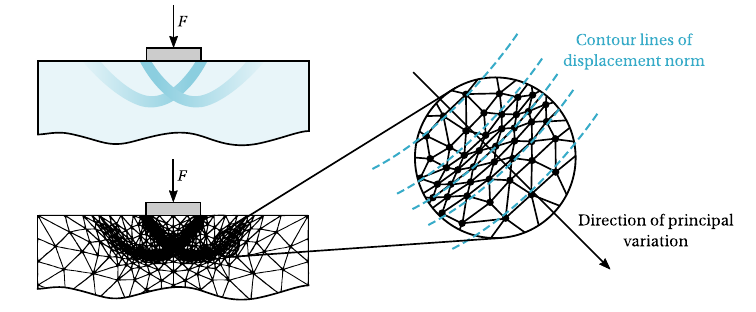}
	\caption{Rigid footing placed on elasto-plastic foundation. Mesh $h$-refinement and elements aspect ratio based on gradient and curvature of displacement norm.}
	\label{fig:adaptivity_element_distortion}
\end{figure}

Attending to their heuristic nature, a great variety of error indicators can be found in the literature and a structured classification similar to the error estimators is not available. Given that the choice of these measures is closely tied with the precise nature of the problem, error indicators are usually quantities that are already available in the finite element computation. For instance, in nonlinear solid mechanics, some common choices are the equivalent plastic strain or its gradient \cite{huerta1999}. In \cite{ortiz1991}, an error indicator is proposed tailored to problems of strain localization and based on the variation of the solution within each element. The adaptive strategy consists of equidistributing the variation of the velocity field over the elements of the mesh and a heuristic justification for the use of variations as indicators is provided. A similar approach is proposed in \cite{jin1990}, being the isolines of effective strain used to perform mesh adaptivity. In the context of thermodynamics with internal variables, error indicators based on a generalized energy norm, the plastic dissipation functional, and the rate of plastic work are explored by Peri\'{c} and coworkers (e.g., \cite{peric1994,peric1998a}). These are found to be an appropriate choice for the adaptive solution of finite strains elasto-plastic problems of practical interest. In the same context, \cite{lee1994} propose pointwise error indicators for stresses and plastic strain increments. In order to deal effectively with material failure and the associated high gradients of the state variables, \cite{peric1999} proposes error indicators based on the rate of fracture indicators. From a geometrical point of view, the element aspect ratio or the distortion can be used, an approach that has proven to be effective when dealing with localization phenomena and discontinuities (e.g., \cite{zienkiewicz1994,zienkiewicz1995a,zienkiewicz1995}). Some methodologies are based on recognizing that the unknown function that we are attempting to model exhibits higher gradients or curvatures in specific directions. A high degree of refinement can be achieved economically in high gradient areas with elongated elements with a suitable orientation \cite{zienkiewicz2014} (see Figure~\ref{fig:adaptivity_element_distortion}). Due to the complexities arising in the error estimation of problems involving path-dependent nonlinear material behavior and, in particular, strain softening and localization, several authors proposed various methodologies based on error indicators (e.g., \cite{pastor1991,yu1991,peric1994,pijaudier-cabot1995,stein1998}).

\subsection{Procedure to adapt spatial interpolation}

\begin{figure}[hbt]
	\centering
	\includegraphics[width=0.9\textwidth]{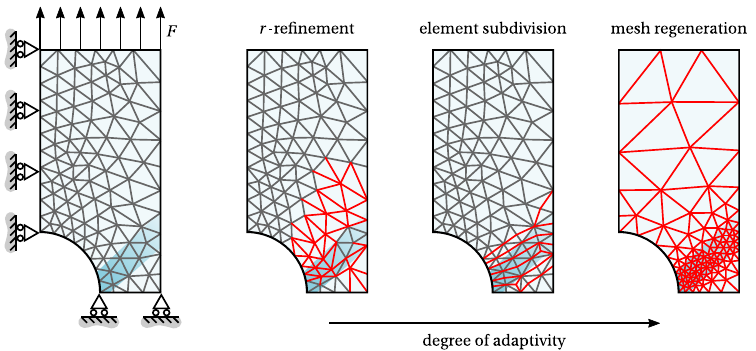}
	\caption{Schematic of the different $h$-refinement adaptive procedures in the analysis of a elasto-plastic perforated specimen (symmetry conditions) under uniaxial traction.}
	\label{fig:h_adaptivity_methods}
\end{figure}

There are essentially two leading families of procedures to perform the
adaptivity of spatial discretization in finite element computations: $h$-adaptivity and $p$-adaptivity. Adaptivity in $h$-adaptive procedures is focused on the size of elements in the mesh, keeping the type of elements unchanged. In the regions where a more accurate solution is needed, the element size is decreased in order to enrich the spatial interpolation (mesh refinement). Likewise, in areas where the accuracy is exceeding the prescribed requirements, the element size may be increased to provide maximum economy in reaching the desired accuracy (mesh de-refinement). As illustrated in Figure~\ref{fig:h_adaptivity_methods}, this family can be further divided into different strategies: (i) $r$-refinement, based on the relocation of existing nodes and keeping the mesh connectivities constant; (ii) element subdivision, based on the division of existing elements into smaller ones and keeping the boundaries of the original elements intact; and (iii) mesh regeneration, consisting in building an entirely new mesh where the elements' optimal size are specified throughout the spatial domain. In contrast, adaptivity in $p$-adaptive procedures is focused on the order of the polynomial that defines each type of element, keeping the elements' size and boundaries unchanged. Therefore, in the regions where great accuracy is required, the polynomial order is increased to enrich the spatial interpolation (mesh refinement). Much work has been put forth in an efficient coupling between $h$- and $p$-adaptive approaches, the so-called $hp$-adaptivity, where both the size of elements and the order of the associated polynomial are simultaneously changed (e.g., \cite{babuska1988,guo1986,guo1986a,demkowicz1989,rachowicz1989,zienkiewicz1989}). A fairly complete comparison between the advantages and disadvantages of these different approaches can be found in \cite{zienkiewicz2013}.

\subsection{Remeshing criterion}

Finally, the remeshing criterion consists of a given strategy that, based on the error assessment, defines when, where and how much should the spatial discretization be adapted through a given procedure. It effectively establishes the bridge between the error estimators/indicators and the adaptive discretization procedure to complete the overall finite element adaptive strategy. For instance, in the context of $h$-refinement, the remeshing criterion should provide information about the required element size throughout the mesh to fulfill the prescribed accuracy requirements. The criteria can be conveniently formulated through a given normalized error estimator or indicator that is compared in some way with a user-defined threshold. Some of the most popular strategies are described in \cite{nochetto2009a}, namely: (i) maximum strategy, where the error is normalized by the maximum value found in the whole mesh and elements whose normalized error is greater than the prescribed threshold are marked for refinement; (ii) equidistribution strategy, similar to the maximum strategy but the error is normalized by the average value over the whole mesh; and (iii) D\"orflers's strategy, where a group of elements holding the largest errors is marked if the group's total error normalized by the total mesh error is greater than the prescribed threshold. Similar criteria to predict the best polynomial order in the case of $p$-refinement are cumbersome to derive.

\section{Self-Consistent Clustering Analysis (SCA) \label{app:sca}}
The Self-Consistent Clustering Analysis (SCA) method proposed by Liu and coworkers \cite{liu2016} can be considered the pioneering clustering-based reduced-order model (CROM) in computational solid mechanics. Given its application in the present paper, a concise description of this method is provided here,\footnote{The nomenclature presented in Section~\ref{ssec:nomenclature} is adopted throughout this Appendix.} following the formulation of the original publication \cite{liu2016}. A more comprehensive description can be found therein and also in \cite{ferreira2022}.

\begin{figure}[hbt]
\centering
\includegraphics[width=0.9\textwidth]{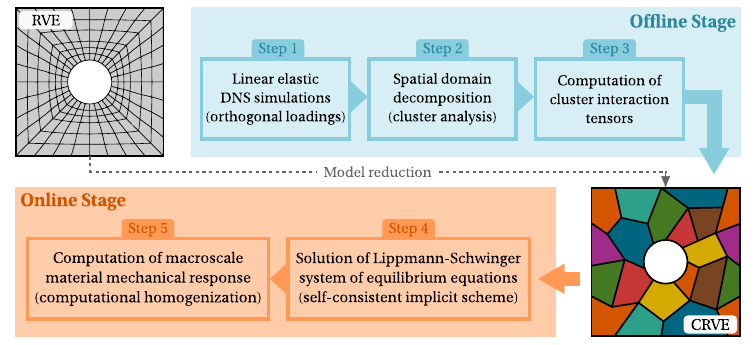}
\caption{Schematic of the Self-Consistent Clustering Analysis (SCA) clustering-based reduced-order model \cite{liu2016}.}
\label{fig:scaschemecompact}
\end{figure}

Similar to most CROMs, the SCA is classified as a two-stage algorithm, meaning that it comprises an offline training/learning stage and an online prediction stage (see Figure~\ref{fig:scaschemecompact}). The offline stage aims to compress the high-fidelity representative volume element (RVE) into a cluster-reduced RVE (CRVE). This model reduction is carried out by performing a clustering-based decomposition of the spatial domain into a given number of material clusters (step 2), where each cluster group points with similar mechanical behavior. In order to do so, the cluster analysis is based on the fourth-order local elastic strain concentration tensor, $\bm{\mathsf{H}}$,
\begin{equation}
	\bm{\varepsilon}_{\mu}^{e} (\bm{Y}) = \bm{\mathsf{H}}^{e} : \bm{\varepsilon}^{e} \, , \quad \forall \, \bm{Y} \in \Omega_{\mu, 0} \, ,
\end{equation}
which establishes the relation between the macroscale strain tensor, $\bm{\varepsilon}^{e}$, and the microscale strain tensor, $\bm{\varepsilon}_{\mu}^{e}$, at each point of the domain. Inspection of the previous relation in matricial form (Voigt's notation),
\begin{equation}
	\mathbf{\upvarepsilon}^{e}_{\mu} = \mathbf{H}^{e} \mathbf{\upvarepsilon}^{e} \quad \rightarrow \quad
	\begin{bmatrix}
		\varepsilon^{e}_{\mu,11} \\[5pt]
		\varepsilon^{e}_{\mu,22} \\[5pt]
		\varepsilon^{e}_{\mu,33} \\[5pt]
		2\varepsilon^{e}_{\mu,12} \\[5pt]
		2\varepsilon^{e}_{\mu,23} \\[5pt]
		2\varepsilon^{e}_{\mu,13}
	\end{bmatrix}
	=
	\begin{bmatrix}
		\mathsf{H}^{e}_{1111} & \mathsf{H}^{e}_{1122} & \mathsf{H}^{e}_{1133} & \mathsf{H}^{e}_{1112} & \mathsf{H}^{e}_{1123} & \mathsf{H}^{e}_{1113} \\[5pt]
		\mathsf{H}^{e}_{2211} & \mathsf{H}^{e}_{2222} & \mathsf{H}^{e}_{2233} & \mathsf{H}^{e}_{2212} & \mathsf{H}^{e}_{2223} & \mathsf{H}^{e}_{2213} \\[5pt]
		\mathsf{H}^{e}_{3311} & \mathsf{H}^{e}_{3322} & \mathsf{H}^{e}_{3333} & \mathsf{H}^{e}_{3312} & \mathsf{H}^{e}_{3323} & \mathsf{H}^{e}_{3313} \\[5pt]
		\mathsf{H}^{e}_{1211} & \mathsf{H}^{e}_{1222} & \mathsf{H}^{e}_{1233} & \mathsf{H}^{e}_{1212} & \mathsf{H}^{e}_{1223} & \mathsf{H}^{e}_{1213} \\[5pt]
		\mathsf{H}^{e}_{2311} & \mathsf{H}^{e}_{2322} & \mathsf{H}^{e}_{2333} & \mathsf{H}^{e}_{2312} & \mathsf{H}^{e}_{2323} & \mathsf{H}^{e}_{2313} \\[5pt]
		\mathsf{H}^{e}_{1311} & \mathsf{H}^{e}_{1322} & \mathsf{H}^{e}_{1333} & \mathsf{H}^{e}_{1312} & \mathsf{H}^{e}_{1323} & \mathsf{H}^{e}_{1313}
	\end{bmatrix}
	\begin{bmatrix}
		\varepsilon^{e}_{11} \\[5pt]
		\varepsilon^{e}_{22} \\[5pt]
		\varepsilon^{e}_{33} \\[5pt]
		2\varepsilon^{e}_{12} \\[5pt]
		2\varepsilon^{e}_{23} \\[5pt]
		2\varepsilon^{e}_{13}
	\end{bmatrix} \, ,
\end{equation}
reveals that this tensor can be determined from the direct numerical solution (DNS) of 6 linear elastic microscale equilibrium problems under orthogonal loading conditions (step 1). For instance, the first column of the matrix $\mathbf{H}^{e}$ is computed by imposing a macroscale strain tensor defined as $\bm{\upvarepsilon}^{e} = [1,0,0,0,0,0]^{T}$, the second column by imposing $\bm{\upvarepsilon}^{e} = [0,1,0,0,0,0]^{T}$, and so on. To conclude the complete characterization of the CRVE, the so-called cluster interaction tensors must be computed between every pair of clusters (step 3). The cluster interaction tensor between clusters $I$ and $J$ is defined as
\begin{equation}
	\bm{\mathsf{T}}^{(I)(J)} = \dfrac{1}{f^{(I)} v_{\mu}} \int_{\Omega_{\mu, \, 0}} \int_{\Omega_{\mu, \, 0}} \chi^{(I)}(\bm{Y}) \, \chi^{(J)} (\bm{Y}') \, \bm{\mathsf{\Phi}}^{0} (\bm{Y}-\bm{Y}') \, \mathrm{d} v' \mathrm{d} v \, , \quad I,J = 1,2,\dots, n_{c} \, , \label{eq:clusterinttensors} 
\end{equation}
where $f^{(I)}$ is the volume fraction of the $I$th cluster,
\begin{equation}
	f^{(I)} = \dfrac{v_{\mu}^{(I)}}{v_{\mu}} \, , \label{eq:clustervolfraction}
\end{equation}
and $\bm{\mathsf{\Phi}}^{0}$ is the Green operator associated with a fictitious elastic reference material. Each of these fourth-order tensors has a non-local nature in the sense that it physically represents the influence of the stress in the $J$th cluster on the strain in the $I$th cluster.

In the following online-stage, the CRVE is loaded with the macroscale strain and/or stress constraints and the microscale equilibrium problem is formulated based on the well-known Lippmann-Schwinger equation,
\begin{equation}
	\bm{\varepsilon}_{\mu} (\bm{Y}) = - \int_{\Omega_{\mu,{0}}} \bm{\mathsf{\Phi}}^{0} (\bm{Y} - \bm{Y}') \,  \left( \bm{\sigma}_{\mu}(\bm{Y}') - \bm{\mathsf{D}}^{e,\,0}  : \bm{\varepsilon}_{\mu}(\bm{Y}') \right)   \, \mathrm{d} v' + \bm{\varepsilon}_{\mu}^{0} \, , \quad \forall \bm{Y} \in \Omega_{\mu, \, 0} \, ,
	\label{eq:LippmannSchwingerFarfield2}
\end{equation}
where $\bm{\mathsf{D}}^{e,\,0}$ is the elasticity tensor of the elastic reference material and $\bm{\varepsilon}_{\mu}^{0}$ is a homogeneous far-field strain. After a suitable time discretization and considering the general (pseudo-)time increment $[t_{m}, \, t_{m+1}]$, the incremental Lippmann-Schwinger integral equilibrium equation can be averaged over each material cluster $I$ as
\begin{equation}
	\Delta \bm{\varepsilon}_{\mu,\,m+1}^{(I)} = - \sum^{n_{c}}_{J=1} \bm{\mathsf{T}}^{(I)(J)} : \left( \Delta \hat{\bm{\sigma}}_{\mu,\,m+1}^{(J)} - \bm{\mathsf{D}}^{e,\, 0}: \Delta\bm{\varepsilon}_{\mu,\,m+1}^{(J)} \right)  + \Delta \bm{\varepsilon}_{\mu,\,m+1}^{0} \, , \quad I = 1,2,\dots, n_{c} \, , \label{eq:clusterLP5}
\end{equation}
being the macroscale strain and stress constraints expressed as
\begin{equation}
	\sum^{n_{c}}_{I=1} f^{(I)}  \Delta \bm{\varepsilon}_{\mu,\,m+1}^{(I)} = \Delta \bm{\varepsilon}_{m+1} (\bm{X}) \, , \qquad \sum^{n_{c}}_{I=1} f^{(I)}  \Delta \hat{\bm{\sigma}}_{\mu,\,m+1}^{(I)} = \Delta \bm{\sigma}_{m+1} (\bm{X}) \, . \label{eq:macroscaleconstraints}
\end{equation}
The Lippmann-Schwinger system of equilibrium equations that must be solved in the online-stage (step 4) is then composed by (i) $n_{c}$ Lippmann-Schwinger integral equilibrium equations (see Equation~\eqref{eq:clusterLP5}) and (ii) macroscale strain and/or stress constraints (see Equation~\eqref{eq:macroscaleconstraints}). Due to the general nonlinearity stemming from the material phase's constitutive models, the equilibrium problem is generally nonlinear and can be efficiently solved through the well-known Newton-Raphson Method. The definition of suitable residual functions and solution procedure of a given macroscale load increment is summarized in Box~\ref{box:scanr}. Once the solution is obtained, the macroscale material response can be computed by computational homogenization (step 5).

Given that the solution of the Lippmann-Schwinger system of equilibrium equations depends on the choice of the elastic reference material properties, the well-known self-consistent micromechanical approach is adopted to determine the `optimal' properties at each macroscale loading increment. By assuming an isotropic elastic reference material, the regression-based self-consistent scheme is then mathematically formulated as an optimization problem,
\begin{equation}
		\left\{ \lambda^{0}_{m+1}, \, \mu^{0}_{m+1} \right\} = \underset{ \left\{ \lambda ', \, \, \mu ' \right\} }{\mathrm{argmin}} \,  \norm{\Delta \bm{\sigma}_{m+1} - \bm{\mathsf{D}}^{e, \, 0}_{m+1}(\lambda ',\mu ') : \Delta \bm{\varepsilon}_{m+1} }^{2} \, ,
\end{equation}
from which the reference material elastic properties, $\lambda^{0}_{m+1}$ and $\mu^{0}_{m+1}$, can be determined at each macroscale loading increment. The solution procedure of the Lippmann-Schwinger system of equilibrium equations summarized in Box~\ref{box:scanr} is thus enriched and embedded within a self-consistent iterative scheme.

\begin{framedbox}[hbt]
\caption{The solution procedure of a given load increment $m+1$ of the Lippmann-Schwinger system of equilibrium equations through the Newton-Raphson Method.}
\begin{center}
  \begin {minipage}{0.90\textwidth}
  	\footnotesize 
  	\smallskip
    \begin{enumerate}
    	\renewcommand{\labelenumi}{(\roman{enumi})}
    	\item Initialize iterative counter, $k \coloneqq 1$
    	\item Set initial guess for incremental strains, $\Delta \bm{\varepsilon}_{\mu, \, m+1 \vphantom{-1}}^{(0)} = \Delta \bm{\varepsilon}_{\mu, \, m}$
    	\item Perform the material state update for each cluster
		\begin{equation*}
			\Delta \bm{\sigma}_{\mu, \, m+1}^{(I)(k-1)} = \Delta \hat{\bm{\sigma}}_{\mu, \, m+1} \left( \Delta \bm{\varepsilon}_{\mu, \, m+1}^{(I)(k-1)}, \, \bm{\alpha}^{(I)}_{m} \right) \, , \qquad \forall I = 1,2, \, \dots, \, n_{c}
		\end{equation*}
    	\item Compute global residual functions
    		\begin{gather*}
	\bm{R}^{(I)}_{m+1} \left( \Delta \bm{\varepsilon}_{\mu, \, m+1 }^{(k-1)} \right)  = \Delta \bm{\varepsilon}_{\mu, \, m+1}^{(I)(k-1)} + \sum^{n_{c}}_{J=1} \bm{\mathsf{T}}^{(I)(J)} : \left[ \Delta \hat{\bm{\sigma}}_{\mu, \, m+1}^{(J)(k-1)}  - \bm{\mathsf{D}}^{e,\, 0}: \Delta\bm{\varepsilon}_{\mu, \, m+1}^{(J)(k-1)} \right]  - \Delta \bm{\varepsilon}_{\mu, \, m+1}^{0, \, (k-1)}\, , \\
	\bm{R}^{(n_{c} + 1)}_{m+1} \left( \Delta \bm{\varepsilon}_{\mu, \, m+1 }^{(k-1)} \right) = \sum^{n_{c}}_{I=1} f^{(I)}  \Delta \bm{\varepsilon}_{\mu, \, m+1}^{(I)} - \Delta \bm{\varepsilon}_{m+1}(\bm{X}) \, ; \\[5pt] \bm{R}^{(n_{c} + 1)}_{m+1} \left( \Delta \bm{\sigma}_{\mu, \, m+1 }^{(k-1)} \right) = \sum^{n_{c}}_{I=1} f^{(I)}  \Delta \bm{\sigma}_{\mu, \, m+1}^{(I)} - \Delta \bm{\sigma}_{m+1}(\bm{X}) \, , \\[5pt] \forall I = 1, 2,\, \dots, \, n_{c}
			\end{gather*}
		\item Check for convergence
		
			if (converged) then
				\begin{itemize}
					\item Update incremental solution, $( \, \bullet \, )_{m+1} = ( \, \bullet \, )_{m+1}^{(k)}$, and exit
				\end{itemize}
					
    	\item Compute the Jacobian matrix 
    		\begin{equation*}
	\bm{J} \left( \Delta \bm{\varepsilon}_{\mu, \, m+1 }^{(k-1)} \right) = \begin{bmatrix}
		\dfrac{\partial \bm{R}^{(I)}_{m+1}}{\partial \Delta \bm{\varepsilon}_{\mu, \, m+1}^{(K)}} &  \dfrac{\partial \bm{R}^{(I)}_{m+1}}{\partial \Delta \bm{\varepsilon}_{\mu, \, m+1}^{0}} \\[20pt]
		\dfrac{\partial \bm{R}^{(n_{c} + 1)}_{m+1}}{\partial \Delta \bm{\varepsilon}_{\mu, \, m+1}^{(K)}} &  \dfrac{\partial \bm{R}^{(n_{c} + 1)}_{m+1}}{\partial \Delta \bm{\varepsilon}_{\mu, \, m+1}^{0}}
		\end{bmatrix}_{ \left( \Delta \bm{\varepsilon}_{\mu, \, m+1 }^{(k-1)} \right)} \, , \qquad \forall I,K=1,2, \, \dots, \, n_{c}
			\end{equation*}
		\item Solve the system of linear equations
			\begin{equation*}
	 			\delta \Delta \bm{\upvarepsilon}_{\mu}^{(k)} = - \mathbf{J}^{-1} \left( \Delta \bm{\varepsilon}_{\mu, \, m+1}^{(k-1)} \right) \, \mathbf{R} \left(\Delta \bm{\varepsilon}_{\mu, \, m+1}^{(k-1)} \right) 
			\end{equation*}
		\item Update incremental strains
			\begin{equation*}
				\Delta \bm{\varepsilon}_{\mu,\, m+1}^{(k)} = \Delta \bm{\varepsilon}_{\mu,\, m+1}^{(k-1)} + \delta \Delta \bm{\varepsilon}_{\mu}^{(k)}
			\end{equation*}
		\item Increment iteration counter, $k:=k+1$
		\item Go to step (iii)
    \end{enumerate}
  \end{minipage}
  \end{center}
\label{box:scanr}
\end{framedbox}

\clearpage

\section{Compatible FEM regular mesh and element averaging \label{app:femrgmesh}}

As described in Section~\ref{ssec:nomenclature}, the RVE is usually discretized in regular voxels in CROMs whose solution procedure is partially computed in the discrete frequency domain. Given that the Finite Element Method (FEM) is selected as the direct numerical solution (DNS) method throughout this paper, two essential procedures are required to perform a suitable transfer and/or comparison of data between both types of spatial discretizations. For instance, these arise when (1) performing the offline-stage DNS simulations with FEM and (2) quantitatively comparing CROMs and FEM microscale field solutions.

\begin{figure}[hbt]
\centering
\includegraphics[width=0.9\textwidth]{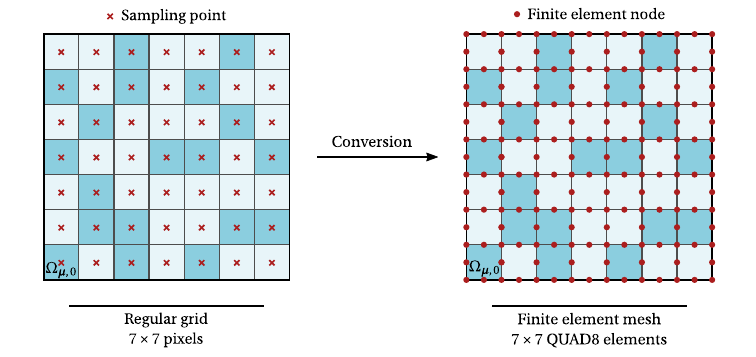}
\caption{Conversion of a regular grid of pixels into a finite element mesh with quadrilateral 8-noded quadratic elements.}
\label{fig:FFT_FEM_Conversion}
\end{figure}

The first procedure involves the generation of a FEM mesh compatible with the regular grid of voxels. A simple way to achieve this is to convert the pixels (2D) and voxels (3D) into quadrilateral (2D) and brick (3D) finite elements, respectively (see Figure~\ref{fig:FFT_FEM_Conversion}). The suitable order and type of finite elements should be defined according to the problem under analysis.

The second procedure allows the transfer and/or comparison between both types of spatial discretizations. The integrations over the finite element domain are generally performed by means of the Gaussian Quadrature Method and the state variables are computed at a given set of Gauss sampling points. Because each pixel (2D) or voxel (3D) only contains one sampling point, the data from the different Gauss sampling points of the associated finite element must be `compressed' in some way. An element volumetric averaging procedure is proposed as follows: for a generic field $a_{\mu}(\bm{Y})$ and the local normalized domain $\Upsilon$ of the element, the volumetric average of $a_{\mu}(\bm{Y})$ among the Gauss sampling points of each element $e$ is taken as
	\begin{equation}
		\bar{a}_{e, \, \mu} = \dfrac{1}{v_{e,\,\mu}} \, \sum_{i=1}^{n_{\mathrm{gp}}} w_{i} a_{\mu}(\bm{\xi}_{i}) j(\bm{\xi}_{i}) \, ,
	\end{equation}	
where $v_{e,\,\mu}$ is the volume of the finite element, $\bm{\xi}_{i}$ and $w_{i}$, $i=1,2,\dots, n_{\mathrm{gp}}$, are the positions and weights of the Gauss sampling points in the domain $\Upsilon$ and $j(\bm{\xi})$ is the determinant of the mapping ($\bm{Y} : \Upsilon \rightarrow \Omega^{(e)}$) Jacobian (see Figure~\ref{fig:FFT_FEM_Conversion}).

\begin{figure}[hbt]
\centering
\includegraphics[width=0.9\textwidth]{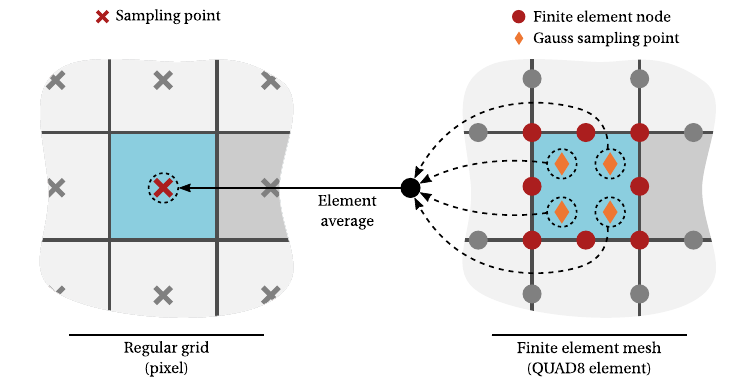}
\caption{Element volumetric averaging as a means to transfer and/or compare data between a finite element mesh and a regular grid of pixels or voxels.}
\label{fig:element_vol_average}
\end{figure}

\clearpage

\section{Accelerating the CRVE scanning procedure \label{app:scanning}}
The target clusters selection criterion described in Section~\ref{ssec:blockA} involves a scanning procedure over the CRVE dimensions. The computational cost of this operation scales with the spatial discretization of the RVE as the number of voxels along with each scanning direction increases. Although the update of the cluster interaction tensors (Block C) dominates the overall cost of the clustering adaptivity procedures, the cost of the selection criterion (Block A) is by no means neglectable. 

\begin{figure}[hbt]
\centering
\includegraphics[width=0.9\textwidth]{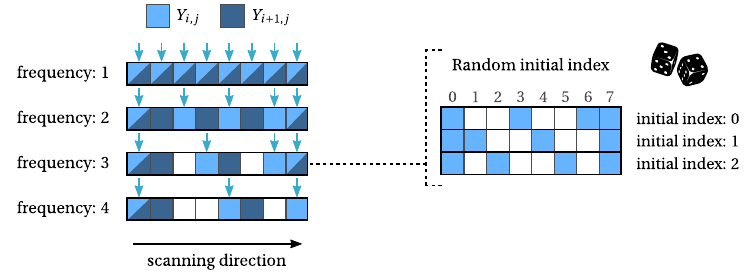}
\caption{Randomized scanning frequency aiming to accelerate the CRVE scanning procedure underlying the target clusters selection criterion based on spatial discontinuities.}
\label{fig:swipefrequency}
\end{figure}

In this context, a simple but effective strategy is proposed here to accelerate the CRVE scanning procedure. In the standard approach, every pair of consecutive voxels is scanned along each direction. However, despite the completeness of this procedure, it can be argued that scanning every voxel is not absolutely necessary to perform a proper assessment and cluster selection. This is especially true when considering a spatial discretization where the voxels' dimensions are small in comparison with the clusters' sizes and/or the characteristic length of the relevant phenomena. The main idea is thus to set a scanning frequency (user-defined parameter) associated with each scanning direction as schematically illustrated in Figure~\ref{fig:swipefrequency}. A frequency of $X$ means that the scanning is performed every $X$ voxels, with $X=1$ recovering the standard (complete) approach, being enforced that the first and last voxels of each scanning direction are always evaluated. In addition, in order to avoid missing any particular domain regions, the initial voxel index can be randomly picked at each clustering adaptive step. Note that a frequency of $X$ leads to $X$ distinct initial voxel indexes defined as $0, 1, \dots, X-1$.

\clearpage


\bibliography{mybibfile}

\begin{thebibliography}{104}
\expandafter\ifx\csname natexlab\endcsname\relax\def\natexlab#1{#1}\fi
\providecommand{\url}[1]{\texttt{#1}}
\providecommand{\href}[2]{#2}
\providecommand{\path}[1]{#1}
\providecommand{\DOIprefix}{doi:}
\providecommand{\ArXivprefix}{arXiv:}
\providecommand{\URLprefix}{URL: }
\providecommand{\Pubmedprefix}{pmid:}
\providecommand{\doi}[1]{\href{http://dx.doi.org/#1}{\path{#1}}}
\providecommand{\Pubmed}[1]{\href{pmid:#1}{\path{#1}}}
\providecommand{\bibinfo}[2]{#2}
\ifx\xfnm\relax \def\xfnm[#1]{\unskip,\space#1}\fi
\bibitem[{Horstemeyer(2012)}]{horstemeyer2012}
\bibinfo{author}{M.~F. Horstemeyer}, \bibinfo{title}{Integrated {{Computational
  Materials Engineering}} ({{ICME}}) for {{Metals}}: {{Using Multiscale
  Modeling}} to {{Invigorate Engineering Design}} with {{Science}}},
  \bibinfo{publisher}{{John Wiley \& Sons}}, \bibinfo{year}{2012}.
\bibitem[{Horstemeyer(2018)}]{horstemeyer2018}
\bibinfo{author}{M.~F. Horstemeyer}, \bibinfo{title}{Integrated {{Computational
  Materials Engineering}} ({{ICME}}) for {{Metals}}: {{Concepts}} and {{Case
  Studies}}}, \bibinfo{publisher}{{John Wiley \& Sons}}, \bibinfo{year}{2018}.
\bibitem[{Smit et~al.(1998)Smit, Brekelmans, and Meijer}]{smit1998}
\bibinfo{author}{R.~J.~M. Smit}, \bibinfo{author}{W.~A.~M. Brekelmans},
  \bibinfo{author}{H.~E.~H. Meijer},
\newblock \bibinfo{title}{Prediction of the mechanical behavior of nonlinear
  heterogeneous systems by multi-level finite element modeling},
\newblock \bibinfo{journal}{Computer Methods in Applied Mechanics and
  Engineering} \bibinfo{volume}{155} (\bibinfo{year}{1998})
  \bibinfo{pages}{181--192}.
\bibitem[{Feyel(1998)}]{feyel1998}
\bibinfo{author}{F.~Feyel}, \bibinfo{title}{Application Du Calcul Parall\`ele
  Aux Mod\`eles \`a Grand Nombre de Variables Internes}, \bibinfo{type}{{{PhD
  Thesis}}}, \'Ecole Nationale Sup\'erieure des Mines de Paris,
  \bibinfo{year}{1998}.
\bibitem[{Miehe et~al.(1999)Miehe, Schr{\"o}der, and Schotte}]{miehe1999}
\bibinfo{author}{C.~Miehe}, \bibinfo{author}{J.~Schr{\"o}der},
  \bibinfo{author}{J.~Schotte},
\newblock \bibinfo{title}{Computational homogenization analysis in finite
  plasticity {{Simulation}} of texture development in polycrystalline
  materials},
\newblock \bibinfo{journal}{Computer Methods in Applied Mechanics and
  Engineering} \bibinfo{volume}{171} (\bibinfo{year}{1999})
  \bibinfo{pages}{387--418}.
\bibitem[{Terada and Kikuchi(2001)}]{terada2001}
\bibinfo{author}{K.~Terada}, \bibinfo{author}{N.~Kikuchi},
\newblock \bibinfo{title}{A class of general algorithms for multi-scale
  analyses of heterogeneous media},
\newblock \bibinfo{journal}{Computer Methods in Applied Mechanics and
  Engineering} \bibinfo{volume}{190} (\bibinfo{year}{2001})
  \bibinfo{pages}{5427--5464}.
\bibitem[{Kochmann et~al.(2016)Kochmann, Wulfinghoff, Reese, Mianroodi, and
  Svendsen}]{kochmann2016}
\bibinfo{author}{J.~Kochmann}, \bibinfo{author}{S.~Wulfinghoff},
  \bibinfo{author}{S.~Reese}, \bibinfo{author}{J.~R. Mianroodi},
  \bibinfo{author}{B.~Svendsen},
\newblock \bibinfo{title}{Two-scale {{FE}}\textendash{{FFT}}- and
  phase-field-based computational modeling of bulk microstructural evolution
  and macroscopic material behavior},
\newblock \bibinfo{journal}{Computer Methods in Applied Mechanics and
  Engineering} \bibinfo{volume}{305} (\bibinfo{year}{2016})
  \bibinfo{pages}{89--110}.
\bibitem[{Yvonnet et~al.(2009)Yvonnet, Gonzalez, and He}]{yvonnet2009}
\bibinfo{author}{J.~Yvonnet}, \bibinfo{author}{D.~Gonzalez},
  \bibinfo{author}{Q.~C. He},
\newblock \bibinfo{title}{Numerically explicit potentials for the
  homogenization of nonlinear elastic heterogeneous materials},
\newblock \bibinfo{journal}{Computer Methods in Applied Mechanics and
  Engineering} \bibinfo{volume}{198} (\bibinfo{year}{2009})
  \bibinfo{pages}{2723--2737}.
\bibitem[{Cl{\'e}ment et~al.(2012)Cl{\'e}ment, Soize, and
  Yvonnet}]{clement2012}
\bibinfo{author}{A.~Cl{\'e}ment}, \bibinfo{author}{C.~Soize},
  \bibinfo{author}{J.~Yvonnet},
\newblock \bibinfo{title}{Computational nonlinear stochastic homogenization
  using a nonconcurrent multiscale approach for hyperelastic heterogeneous
  microstructures analysis},
\newblock \bibinfo{journal}{International Journal for Numerical Methods in
  Engineering} \bibinfo{volume}{91} (\bibinfo{year}{2012})
  \bibinfo{pages}{799--824}.
\bibitem[{Le et~al.(2015)Le, Yvonnet, and He}]{le2015}
\bibinfo{author}{B.~A. Le}, \bibinfo{author}{J.~Yvonnet},
  \bibinfo{author}{Q.-C. He},
\newblock \bibinfo{title}{Computational homogenization of nonlinear elastic
  materials using neural networks},
\newblock \bibinfo{journal}{International Journal for Numerical Methods in
  Engineering} \bibinfo{volume}{104} (\bibinfo{year}{2015})
  \bibinfo{pages}{1061--1084}.
\bibitem[{Bessa et~al.(2017)Bessa, Bostanabad, Liu, Hu, Apley, Brinson, Chen,
  and Liu}]{bessa2017}
\bibinfo{author}{M.~Bessa}, \bibinfo{author}{R.~Bostanabad},
  \bibinfo{author}{Z.~Liu}, \bibinfo{author}{A.~Hu}, \bibinfo{author}{D.~W.
  Apley}, \bibinfo{author}{C.~Brinson}, \bibinfo{author}{W.~Chen},
  \bibinfo{author}{W.~K. Liu},
\newblock \bibinfo{title}{A framework for data-driven analysis of materials
  under uncertainty: {{Countering}} the curse of dimensionality},
\newblock \bibinfo{journal}{Computer Methods in Applied Mechanics and
  Engineering} \bibinfo{volume}{320} (\bibinfo{year}{2017})
  \bibinfo{pages}{633--667}.
\bibitem[{Liu et~al.(2016)Liu, Bessa, and Liu}]{liu2016}
\bibinfo{author}{Z.~Liu}, \bibinfo{author}{M.~Bessa}, \bibinfo{author}{W.~K.
  Liu},
\newblock \bibinfo{title}{Self-consistent clustering analysis: {{An}} efficient
  multi-scale scheme for inelastic heterogeneous materials},
\newblock \bibinfo{journal}{Computer Methods in Applied Mechanics and
  Engineering} \bibinfo{volume}{306} (\bibinfo{year}{2016})
  \bibinfo{pages}{319--341}.
\bibitem[{Wulfinghoff et~al.(2017)Wulfinghoff, Cavaliere, and
  Reese}]{wulfinghoff2017}
\bibinfo{author}{S.~Wulfinghoff}, \bibinfo{author}{F.~Cavaliere},
  \bibinfo{author}{S.~Reese},
\newblock \bibinfo{title}{Model order reduction of nonlinear homogenization
  problems using a {{Hashin}}-{{Shtrikman}} type finite element method},
\newblock \bibinfo{journal}{Computer Methods in Applied Mechanics and
  Engineering} \bibinfo{volume}{330} (\bibinfo{year}{2017}).
\bibitem[{Dvorak(1992)}]{dvorak1992}
\bibinfo{author}{G.~J. Dvorak},
\newblock \bibinfo{title}{Transformation field analysis of inelastic composite
  materials},
\newblock \bibinfo{journal}{Proceedings of the Royal Society of London. Series
  A: Mathematical and Physical Sciences} \bibinfo{volume}{437}
  (\bibinfo{year}{1992}) \bibinfo{pages}{311--327}.
\bibitem[{Michel and Suquet(2003)}]{michel2003}
\bibinfo{author}{J.~C. Michel}, \bibinfo{author}{P.~Suquet},
\newblock \bibinfo{title}{Nonuniform transformation field analysis},
\newblock \bibinfo{journal}{International Journal of Solids and Structures}
  \bibinfo{volume}{40} (\bibinfo{year}{2003}) \bibinfo{pages}{6937--6955}.
\bibitem[{Ladev{\`e}ze et~al.(2010)Ladev{\`e}ze, Passieux, and
  N{\'e}ron}]{ladeveze2010}
\bibinfo{author}{P.~Ladev{\`e}ze}, \bibinfo{author}{J.~C. Passieux},
  \bibinfo{author}{D.~N{\'e}ron},
\newblock \bibinfo{title}{The {{LATIN}} multiscale computational method and the
  {{Proper Generalized Decomposition}}},
\newblock \bibinfo{journal}{Computer Methods in Applied Mechanics and
  Engineering} \bibinfo{volume}{199} (\bibinfo{year}{2010})
  \bibinfo{pages}{1287--1296}.
\bibitem[{Boyaval(2008)}]{boyaval2008}
\bibinfo{author}{S.~Boyaval},
\newblock \bibinfo{title}{Reduced-{{Basis Approach}} for {{Homogenization}}
  beyond the {{Periodic Setting}}},
\newblock \bibinfo{journal}{Multiscale Modeling \& Simulation}
  \bibinfo{volume}{7} (\bibinfo{year}{2008}) \bibinfo{pages}{466--494}.
\bibitem[{Hern{\'a}ndez et~al.(2014)Hern{\'a}ndez, Oliver, Huespe, Caicedo, and
  Cante}]{hernandez2014}
\bibinfo{author}{J.~Hern{\'a}ndez}, \bibinfo{author}{J.~Oliver},
  \bibinfo{author}{A.~Huespe}, \bibinfo{author}{M.~Caicedo},
  \bibinfo{author}{J.~Cante},
\newblock \bibinfo{title}{High-performance model reduction techniques in
  computational multiscale homogenization},
\newblock \bibinfo{journal}{Computer Methods in Applied Mechanics and
  Engineering} \bibinfo{volume}{276} (\bibinfo{year}{2014})
  \bibinfo{pages}{149--189}.
\bibitem[{Hern{\'a}ndez et~al.(2017)Hern{\'a}ndez, Caicedo, and
  Ferrer}]{hernandez2017}
\bibinfo{author}{J.~A. Hern{\'a}ndez}, \bibinfo{author}{M.~A. Caicedo},
  \bibinfo{author}{A.~Ferrer},
\newblock \bibinfo{title}{Dimensional hyper-reduction of nonlinear finite
  element models via empirical cubature},
\newblock \bibinfo{journal}{Computer Methods in Applied Mechanics and
  Engineering} \bibinfo{volume}{313} (\bibinfo{year}{2017})
  \bibinfo{pages}{687--722}.
\bibitem[{{van Tuijl} et~al.(2019){van Tuijl}, Harnish, Matou{\v s}, Remmers,
  and Geers}]{vantuijl2019}
\bibinfo{author}{R.~A. {van Tuijl}}, \bibinfo{author}{C.~Harnish},
  \bibinfo{author}{K.~Matou{\v s}}, \bibinfo{author}{J.~J.~C. Remmers},
  \bibinfo{author}{M.~G.~D. Geers},
\newblock \bibinfo{title}{Wavelet based reduced order models for
  microstructural analyses},
\newblock \bibinfo{journal}{Computational Mechanics} \bibinfo{volume}{63}
  (\bibinfo{year}{2019}) \bibinfo{pages}{535--554}.
\bibitem[{Schneider(2019)}]{schneider2019}
\bibinfo{author}{M.~Schneider},
\newblock \bibinfo{title}{On the mathematical foundations of the
  self-consistent clustering analysis for non-linear materials at small
  strains},
\newblock \bibinfo{journal}{Computer Methods in Applied Mechanics and
  Engineering} \bibinfo{volume}{354} (\bibinfo{year}{2019})
  \bibinfo{pages}{783--801}.
\bibitem[{Liu et~al.(2018)Liu, Fleming, and Liu}]{liu2018}
\bibinfo{author}{Z.~Liu}, \bibinfo{author}{M.~Fleming}, \bibinfo{author}{W.~K.
  Liu},
\newblock \bibinfo{title}{Microstructural material database for self-consistent
  clustering analysis of elastoplastic strain softening materials},
\newblock \bibinfo{journal}{Computer Methods in Applied Mechanics and
  Engineering} \bibinfo{volume}{330} (\bibinfo{year}{2018})
  \bibinfo{pages}{547--577}.
\bibitem[{Tang et~al.(2018)Tang, Zhang, and Liu}]{tang2018}
\bibinfo{author}{S.~Tang}, \bibinfo{author}{L.~Zhang}, \bibinfo{author}{W.~K.
  Liu},
\newblock \bibinfo{title}{From virtual clustering analysis to self-consistent
  clustering analysis: A mathematical study},
\newblock \bibinfo{journal}{Computational Mechanics} \bibinfo{volume}{62}
  (\bibinfo{year}{2018}) \bibinfo{pages}{1443--1460}.
\bibitem[{Yu et~al.(2019)Yu, Kafka, and Liu}]{yu2019}
\bibinfo{author}{C.~Yu}, \bibinfo{author}{O.~L. Kafka}, \bibinfo{author}{W.~K.
  Liu},
\newblock \bibinfo{title}{Self-consistent clustering analysis for multiscale
  modeling at finite strains},
\newblock \bibinfo{journal}{Computer Methods in Applied Mechanics and
  Engineering} \bibinfo{volume}{349} (\bibinfo{year}{2019})
  \bibinfo{pages}{339--359}.
\bibitem[{Nie et~al.(2019)Nie, Cheng, Li, Xu, and Li}]{nie2019}
\bibinfo{author}{Y.~Nie}, \bibinfo{author}{G.~Cheng}, \bibinfo{author}{X.~Li},
  \bibinfo{author}{L.~Xu}, \bibinfo{author}{K.~Li},
\newblock \bibinfo{title}{Principle of cluster minimum complementary energy of
  {{FEM}}-cluster-based reduced order method: Fast updating the interaction
  matrix and predicting effective nonlinear properties of heterogeneous
  material},
\newblock \bibinfo{journal}{Computational Mechanics} \bibinfo{volume}{64}
  (\bibinfo{year}{2019}) \bibinfo{pages}{323--349}.
\bibitem[{Yu et~al.(2021)Yu, Kafka, and Liu}]{yu2021}
\bibinfo{author}{C.~Yu}, \bibinfo{author}{O.~L. Kafka}, \bibinfo{author}{W.~K.
  Liu},
\newblock \bibinfo{title}{Multiresolution clustering analysis for efficient
  modeling of hierarchical material systems},
\newblock \bibinfo{journal}{Computational Mechanics} \bibinfo{volume}{67}
  (\bibinfo{year}{2021}) \bibinfo{pages}{1293--1306}.
\bibitem[{Li et~al.(2019)Li, Kafka, Gao, Yu, Nie, Zhang, Tajdari, Tang, Guo,
  Li, Tang, Cheng, and Liu}]{li2019}
\bibinfo{author}{H.~Li}, \bibinfo{author}{O.~L. Kafka},
  \bibinfo{author}{J.~Gao}, \bibinfo{author}{C.~Yu}, \bibinfo{author}{Y.~Nie},
  \bibinfo{author}{L.~Zhang}, \bibinfo{author}{M.~Tajdari},
  \bibinfo{author}{S.~Tang}, \bibinfo{author}{X.~Guo}, \bibinfo{author}{G.~Li},
  \bibinfo{author}{S.~Tang}, \bibinfo{author}{G.~Cheng}, \bibinfo{author}{W.~K.
  Liu},
\newblock \bibinfo{title}{Clustering discretization methods for generation of
  material performance databases in machine learning and design optimization},
\newblock \bibinfo{journal}{Computational Mechanics} \bibinfo{volume}{64}
  (\bibinfo{year}{2019}) \bibinfo{pages}{281--305}.
\bibitem[{Yan et~al.(2018)Yan, Lin, Kafka, Yu, Liu, Lian, Wolff, Cao, Wagner,
  and Liu}]{yan2018}
\bibinfo{author}{W.~Yan}, \bibinfo{author}{S.~Lin}, \bibinfo{author}{O.~L.
  Kafka}, \bibinfo{author}{C.~Yu}, \bibinfo{author}{Z.~Liu},
  \bibinfo{author}{Y.~Lian}, \bibinfo{author}{S.~Wolff},
  \bibinfo{author}{J.~Cao}, \bibinfo{author}{G.~J. Wagner},
  \bibinfo{author}{W.~K. Liu},
\newblock \bibinfo{title}{Modeling process-structure-property relationships for
  additive manufacturing},
\newblock \bibinfo{journal}{Frontiers of Mechanical Engineering}
  \bibinfo{volume}{13} (\bibinfo{year}{2018}) \bibinfo{pages}{482--492}.
\bibitem[{Shakoor et~al.(2019)Shakoor, Kafka, Yu, and Liu}]{shakoor2019}
\bibinfo{author}{M.~Shakoor}, \bibinfo{author}{O.~L. Kafka},
  \bibinfo{author}{C.~Yu}, \bibinfo{author}{W.~K. Liu},
\newblock \bibinfo{title}{Data science for finite strain mechanical science of
  ductile materials},
\newblock \bibinfo{journal}{Computational Mechanics} \bibinfo{volume}{64}
  (\bibinfo{year}{2019}) \bibinfo{pages}{33--45}.
\bibitem[{Cavaliere et~al.(2019)Cavaliere, Reese, and
  Wulfinghoff}]{cavaliere2019}
\bibinfo{author}{F.~Cavaliere}, \bibinfo{author}{S.~Reese},
  \bibinfo{author}{S.~Wulfinghoff},
\newblock \bibinfo{title}{Efficient two\textendash scale simulations of
  engineering structures using the {{Hashin}}\textendash{{Shtrikman}} type
  finite element method},
\newblock \bibinfo{journal}{Computational Mechanics}  (\bibinfo{year}{2019}).
\bibitem[{He et~al.(2020)He, Gao, Li, Ge, Chen, Liu, and Fang}]{he2020}
\bibinfo{author}{C.~He}, \bibinfo{author}{J.~Gao}, \bibinfo{author}{H.~Li},
  \bibinfo{author}{J.~Ge}, \bibinfo{author}{Y.~Chen}, \bibinfo{author}{J.~Liu},
  \bibinfo{author}{D.~Fang},
\newblock \bibinfo{title}{A data-driven self-consistent clustering analysis for
  the progressive damage behavior of {{3D}} braided composites},
\newblock \bibinfo{journal}{Composite Structures} \bibinfo{volume}{249}
  (\bibinfo{year}{2020}) \bibinfo{pages}{112471}.
\bibitem[{Gao et~al.(2020)Gao, Shakoor, Domel, Merzkirch, Zhou, Zeng, Su, and
  Liu}]{gao2020}
\bibinfo{author}{J.~Gao}, \bibinfo{author}{M.~Shakoor},
  \bibinfo{author}{G.~Domel}, \bibinfo{author}{M.~Merzkirch},
  \bibinfo{author}{G.~Zhou}, \bibinfo{author}{D.~Zeng},
  \bibinfo{author}{X.~Su}, \bibinfo{author}{W.~K. Liu},
\newblock \bibinfo{title}{Predictive multiscale modeling for {{Unidirectional
  Carbon Fiber Reinforced Polymers}}},
\newblock \bibinfo{journal}{Composites Science and Technology}
  \bibinfo{volume}{186} (\bibinfo{year}{2020}) \bibinfo{pages}{107922}.
\bibitem[{Han et~al.(2020)Han, Gao, Fleming, Xu, Xie, Meng, and Liu}]{han2020}
\bibinfo{author}{X.~Han}, \bibinfo{author}{J.~Gao},
  \bibinfo{author}{M.~Fleming}, \bibinfo{author}{C.~Xu},
  \bibinfo{author}{W.~Xie}, \bibinfo{author}{S.~Meng}, \bibinfo{author}{W.~K.
  Liu},
\newblock \bibinfo{title}{Efficient multiscale modeling for woven composites
  based on self-consistent clustering analysis},
\newblock \bibinfo{journal}{Computer Methods in Applied Mechanics and
  Engineering} \bibinfo{volume}{364} (\bibinfo{year}{2020})
  \bibinfo{pages}{112929}.
\bibitem[{Babu{\v s}ka and Rheinboldt(1978)}]{babuska1978}
\bibinfo{author}{I.~Babu{\v s}ka}, \bibinfo{author}{W.~C. Rheinboldt},
\newblock \bibinfo{title}{A-posteriori error estimates for the finite element
  method},
\newblock \bibinfo{journal}{International Journal for Numerical Methods in
  Engineering} \bibinfo{volume}{12} (\bibinfo{year}{1978})
  \bibinfo{pages}{1597--1615}.
\bibitem[{Babu{\v s}ka and Rheinboldt(1979)}]{babuska1979}
\bibinfo{author}{I.~Babu{\v s}ka}, \bibinfo{author}{W.~C. Rheinboldt},
\newblock \bibinfo{title}{Adaptive approaches and reliability estimations in
  finite element analysis},
\newblock \bibinfo{journal}{Computer Methods in Applied Mechanics and
  Engineering} \bibinfo{volume}{17-18} (\bibinfo{year}{1979})
  \bibinfo{pages}{519--540}.
\bibitem[{Ladev{\`e}ze and Oden(1998)}]{ladeveze1998}
\bibinfo{editor}{P.~Ladev{\`e}ze}, \bibinfo{editor}{J.~T. Oden} (Eds.),
  \bibinfo{title}{Advances in Adaptive Computational Methods in Mechanics},
  number~\bibinfo{number}{47} in \bibinfo{series}{Studies in Applied
  Mechanics}, \bibinfo{publisher}{{Elsevier}}, \bibinfo{address}{{Amsterdam ;
  New York}}, \bibinfo{year}{1998}.
\bibitem[{Ainsworth and Oden(2000)}]{ainsworth2000}
\bibinfo{author}{M.~Ainsworth}, \bibinfo{author}{J.~T. Oden}, \bibinfo{title}{A
  {{Posteriori Error Estimation}} in {{Finite Element}} Analysis},
  \bibinfo{publisher}{{John Wiley \& Sons, Ltd}}, \bibinfo{year}{2000}.
\bibitem[{Stein(2005)}]{stein2005}
\bibinfo{editor}{E.~Stein} (Ed.), \bibinfo{title}{Adaptive {{Finite Elements}}
  in {{Linear}} and {{Nonlinear Solid}} and {{Structural Mechanics}}}, {{CISM
  International Centre}} for {{Mechanical Sciences}},
  \bibinfo{publisher}{{Springer-Verlag}}, \bibinfo{address}{{Wien}},
  \bibinfo{year}{2005}. \DOIprefix\doi{10.1007/3-211-38060-4}.
\bibitem[{Zienkiewicz et~al.(2013)Zienkiewicz, Taylor, and
  Zhu}]{zienkiewicz2013}
\bibinfo{author}{O.~C. Zienkiewicz}, \bibinfo{author}{R.~L. Taylor},
  \bibinfo{author}{J.~Z. Zhu}, \bibinfo{title}{The Finite Element Method: Its
  Basis and Fundamentals}, \bibinfo{edition}{seventh edition} ed.,
  \bibinfo{publisher}{{Elsevier, Butterworth-Heinemann}},
  \bibinfo{address}{{Amsterdam}}, \bibinfo{year}{2013}.
\bibitem[{Huerta et~al.(1999)Huerta, Rodr{\'i}guez-Ferran, D{\'i}ez, and
  Sarrate}]{huerta1999}
\bibinfo{author}{A.~Huerta}, \bibinfo{author}{A.~Rodr{\'i}guez-Ferran},
  \bibinfo{author}{P.~D{\'i}ez}, \bibinfo{author}{J.~Sarrate},
\newblock \bibinfo{title}{Adaptive finite element strategies based on error
  assessment},
\newblock \bibinfo{journal}{International Journal for Numerical Methods in
  Engineering} \bibinfo{volume}{46} (\bibinfo{year}{1999})
  \bibinfo{pages}{1803--1818}.
\bibitem[{Zienkiewicz et~al.(2014)Zienkiewicz, Taylor, and
  Fox}]{zienkiewicz2014}
\bibinfo{author}{O.~C. Zienkiewicz}, \bibinfo{author}{R.~L. Taylor},
  \bibinfo{author}{D.~Fox}, \bibinfo{title}{The Finite Element Method for Solid
  and Structural Mechanics}, \bibinfo{edition}{7th ed} ed.,
  \bibinfo{publisher}{{Elsevier/Butterworth-Heinemann}},
  \bibinfo{address}{{Amsterdam ; Boston}}, \bibinfo{year}{2014}.
\bibitem[{Rice(1976)}]{rice1976}
\bibinfo{author}{J.~R. Rice},
\newblock \bibinfo{title}{The localization of plastic deformation},
\newblock in: \bibinfo{booktitle}{In: {{W}}.{{T}}. {{Koiter}} ({{Ed}}.),
  {{Theoretical}} and {{Applied Mechanics}}},
  \bibinfo{publisher}{{North-Holland Publishing Company}},
  \bibinfo{year}{1976}, pp. \bibinfo{pages}{207--220}.
\bibitem[{Pietruszczak and Mr{\'o}z(1981)}]{pietruszczak1981}
\bibinfo{author}{S.~Pietruszczak}, \bibinfo{author}{Z.~Mr{\'o}z},
\newblock \bibinfo{title}{Finite element analysis of deformation of
  strain-softening materials},
\newblock \bibinfo{journal}{International Journal for Numerical Methods in
  Engineering} \bibinfo{volume}{17} (\bibinfo{year}{1981})
  \bibinfo{pages}{327--334}.
\bibitem[{Ortiz and Quigley(1991)}]{ortiz1991}
\bibinfo{author}{M.~Ortiz}, \bibinfo{author}{J.~J. Quigley},
\newblock \bibinfo{title}{Adaptive mesh refinement in strain localization
  problems},
\newblock \bibinfo{journal}{Computer Methods in Applied Mechanics and
  Engineering} \bibinfo{volume}{90} (\bibinfo{year}{1991})
  \bibinfo{pages}{781--804}.
\bibitem[{De~Borst et~al.(1993)De~Borst, Sluys, Muhlhaus, and
  Pamin}]{deborst1993a}
\bibinfo{author}{R.~De~Borst}, \bibinfo{author}{L.~J. Sluys},
  \bibinfo{author}{H.~B. Muhlhaus}, \bibinfo{author}{J.~Pamin},
\newblock \bibinfo{title}{Fundamental issues in finite element analyses of
  localization of deformation},
\newblock \bibinfo{journal}{Engineering Computations 10(2), 99-121. (1993)}
  (\bibinfo{year}{1993}).
\bibitem[{Hutchinson and Tvergaard(1980)}]{hutchinson1980}
\bibinfo{author}{J.~W. Hutchinson}, \bibinfo{author}{V.~Tvergaard},
\newblock \bibinfo{title}{Surface instabilities on statically strained plastic
  solids},
\newblock \bibinfo{journal}{International Journal of Mechanical Sciences}
  \bibinfo{volume}{22} (\bibinfo{year}{1980}) \bibinfo{pages}{339--354}.
\bibitem[{Ortiz et~al.(1987)Ortiz, Leroy, and Needleman}]{ortiz1987}
\bibinfo{author}{M.~Ortiz}, \bibinfo{author}{Y.~Leroy},
  \bibinfo{author}{A.~Needleman},
\newblock \bibinfo{title}{A finite element method for localized failure
  analysis},
\newblock \bibinfo{journal}{Computer Methods in Applied Mechanics and
  Engineering} \bibinfo{volume}{61} (\bibinfo{year}{1987})
  \bibinfo{pages}{189--214}.
\bibitem[{Huerta et~al.(1998)Huerta, D{\'i}ez, {Rodr{\'i}gues-Ferran}, and
  {Pijaudier-Cabot}}]{huerta1998}
\bibinfo{author}{A.~Huerta}, \bibinfo{author}{P.~D{\'i}ez},
  \bibinfo{author}{A.~{Rodr{\'i}gues-Ferran}},
  \bibinfo{author}{G.~{Pijaudier-Cabot}},
\newblock \bibinfo{title}{Error estimation and adaptive finite element analysis
  of softening solids},
\newblock in: \bibinfo{booktitle}{Advances in Adaptive Computational Methods in
  Mechanics}, number~\bibinfo{number}{47} in \bibinfo{series}{Studies in
  Applied Mechanics}, \bibinfo{publisher}{{Elsevier}},
  \bibinfo{address}{{Amsterdam ; New York}}, \bibinfo{year}{1998}, pp.
  \bibinfo{pages}{333--347}.
\bibitem[{Zienkiewicz and Cormeau(1974)}]{zienkiewicz1974}
\bibinfo{author}{O.~C. Zienkiewicz}, \bibinfo{author}{I.~C. Cormeau},
\newblock \bibinfo{title}{Visco-plasticity\textemdash plasticity and creep in
  elastic solids\textemdash a unified numerical solution approach},
\newblock \bibinfo{journal}{International Journal for Numerical Methods in
  Engineering} \bibinfo{volume}{8} (\bibinfo{year}{1974})
  \bibinfo{pages}{821--845}.
\bibitem[{Bazant and {Pijaudier-Cabot}(1988)}]{bazant1988a}
\bibinfo{author}{Z.~P. Bazant}, \bibinfo{author}{G.~{Pijaudier-Cabot}},
\newblock \bibinfo{title}{Nonlocal {{Continuum Damage}}, {{Localization
  Instability}} and {{Convergence}}},
\newblock \bibinfo{journal}{Journal of Applied Mechanics} \bibinfo{volume}{55}
  (\bibinfo{year}{1988}) \bibinfo{pages}{287--293}.
\bibitem[{Reis et~al.(2018)Reis, Rodrigues~Lopes, Andrade~Pires, and
  Andrade}]{reis2018}
\bibinfo{author}{F.~J.~P. Reis}, \bibinfo{author}{I.~A. Rodrigues~Lopes},
  \bibinfo{author}{F.~M. Andrade~Pires}, \bibinfo{author}{F.~X.~C. Andrade},
\newblock \bibinfo{title}{Microscale analysis of heterogeneous ductile
  materials with nonlocal damage models of integral type},
\newblock \bibinfo{journal}{Computers \& Structures} \bibinfo{volume}{201}
  (\bibinfo{year}{2018}) \bibinfo{pages}{37--57}.
\bibitem[{{de Borst} et~al.(1995){de Borst}, Pamin, Peerlings, and
  Sluys}]{deborst1995}
\bibinfo{author}{R.~{de Borst}}, \bibinfo{author}{J.~Pamin},
  \bibinfo{author}{R.~H.~J. Peerlings}, \bibinfo{author}{L.~J. Sluys},
\newblock \bibinfo{title}{On gradient-enhanced damage and plasticity models for
  failure in quasi-brittle and frictional materials},
\newblock \bibinfo{journal}{Computational Mechanics} \bibinfo{volume}{17}
  (\bibinfo{year}{1995}) \bibinfo{pages}{130--141}.
\bibitem[{Geers et~al.(1998)Geers, R.~Borst, Brekelmans, and
  Peerlings}]{geers1998}
\bibinfo{author}{M.~G.~D. Geers}, \bibinfo{author}{D.~R.~Borst},
  \bibinfo{author}{W.~a.~M. Brekelmans}, \bibinfo{author}{R.~H.~J. Peerlings},
\newblock \bibinfo{title}{Strain-based transient-gradient damage model for
  failure analyses},
\newblock \bibinfo{journal}{Computer Methods in Applied Mechanics and
  Engineering} \bibinfo{volume}{160} (\bibinfo{year}{1998})
  \bibinfo{pages}{133--153}.
\bibitem[{Kuhl et~al.(2000)Kuhl, Ramm, and Borst}]{kuhl2000}
\bibinfo{author}{E.~Kuhl}, \bibinfo{author}{E.~Ramm},
  \bibinfo{author}{R.~Borst},
\newblock \bibinfo{title}{An anisotropic gradient damage model for
  quasi-brittle materials},
\newblock \bibinfo{journal}{Computer Methods in Applied Mechanics and
  Engineering} \bibinfo{volume}{183} (\bibinfo{year}{2000})
  \bibinfo{pages}{87--103}.
\bibitem[{Francfort and Marigo(1998)}]{francfort1998}
\bibinfo{author}{G.~A. Francfort}, \bibinfo{author}{J.-J. Marigo},
\newblock \bibinfo{title}{Revisiting brittle fracture as an energy minimization
  problem},
\newblock \bibinfo{journal}{Journal of Mechanics Physics of Solids}
  \bibinfo{volume}{46} (\bibinfo{year}{1998}) \bibinfo{pages}{1319--1342}.
\bibitem[{Bourdin et~al.(2008)Bourdin, Francfort, and Marigo}]{bourdin2008}
\bibinfo{author}{B.~Bourdin}, \bibinfo{author}{G.~A. Francfort},
  \bibinfo{author}{J.-J. Marigo},
\newblock \bibinfo{title}{The {{Variational Approach}} to {{Fracture}}},
\newblock \bibinfo{journal}{Journal of Elasticity} \bibinfo{volume}{91}
  (\bibinfo{year}{2008}) \bibinfo{pages}{5--148}.
\bibitem[{Miehe et~al.(2010)Miehe, Hofacker, and Welschinger}]{miehe2010}
\bibinfo{author}{C.~Miehe}, \bibinfo{author}{M.~Hofacker},
  \bibinfo{author}{F.~Welschinger},
\newblock \bibinfo{title}{A phase field model for rate-independent crack
  propagation: {{Robust}} algorithmic implementation based on operator splits},
\newblock \bibinfo{journal}{Computer Methods in Applied Mechanics and
  Engineering} \bibinfo{volume}{199} (\bibinfo{year}{2010})
  \bibinfo{pages}{2765--2778}.
\bibitem[{Shahba and Ghosh(2019)}]{shahba2019}
\bibinfo{author}{A.~Shahba}, \bibinfo{author}{S.~Ghosh},
\newblock \bibinfo{title}{Coupled phase field finite element model for crack
  propagation in elastic polycrystalline microstructures},
\newblock \bibinfo{journal}{International Journal of Fracture}
  (\bibinfo{year}{2019}).
\bibitem[{Oliver et~al.(1993)Oliver, Sim{\'o}, and Armero}]{oliver1993}
\bibinfo{author}{J.~Oliver}, \bibinfo{author}{J.~Sim{\'o}},
  \bibinfo{author}{F.~Armero},
\newblock \bibinfo{title}{An analysis of strong discontinuities induced by
  strain-softening in rate-independent inelastic solids},
\newblock \bibinfo{journal}{Computational Mechanics}  (\bibinfo{year}{1993}).
\bibitem[{Oliver et~al.(1999)Oliver, Cervera, and Manzoli}]{oliver1999}
\bibinfo{author}{J.~Oliver}, \bibinfo{author}{M.~Cervera},
  \bibinfo{author}{O.~Manzoli},
\newblock \bibinfo{title}{Strong discontinuities and continuum plasticity
  models: The strong discontinuity approach},
\newblock \bibinfo{journal}{International Journal of Plasticity}
  \bibinfo{volume}{15} (\bibinfo{year}{1999}) \bibinfo{pages}{319--351}.
\bibitem[{Mo{\"e}s et~al.(1999)Mo{\"e}s, Dolbow, and Belytschko}]{moes1999}
\bibinfo{author}{N.~Mo{\"e}s}, \bibinfo{author}{J.~Dolbow},
  \bibinfo{author}{T.~Belytschko},
\newblock \bibinfo{title}{A finite element method for crack growth without
  remeshing},
\newblock \bibinfo{journal}{International Journal for Numerical Methods in
  Engineering} \bibinfo{volume}{46} (\bibinfo{year}{1999})
  \bibinfo{pages}{131--150}.
\bibitem[{Chessa et~al.(2002)Chessa, Smolinski, and Belytschko}]{chessa2002}
\bibinfo{author}{J.~Chessa}, \bibinfo{author}{P.~Smolinski},
  \bibinfo{author}{T.~Belytschko},
\newblock \bibinfo{title}{The extended finite element method ({{XFEM}}) for
  solidification problems},
\newblock \bibinfo{journal}{International Journal for Numerical Methods in
  Engineering} \bibinfo{volume}{53} (\bibinfo{year}{2002})
  \bibinfo{pages}{1959--1977}.
\bibitem[{Mo{\"e}s and Belytschko(2002)}]{moes2002}
\bibinfo{author}{N.~Mo{\"e}s}, \bibinfo{author}{T.~Belytschko},
\newblock \bibinfo{title}{Extended finite element method for cohesive crack
  growth},
\newblock \bibinfo{journal}{Engineering Fracture Mechanics}
  \bibinfo{volume}{69} (\bibinfo{year}{2002}) \bibinfo{pages}{813--833}.
\bibitem[{Hillerborg et~al.(1976)Hillerborg, Mod{\'e}er, and
  Petersson}]{hillerborg1976}
\bibinfo{author}{A.~Hillerborg}, \bibinfo{author}{M.~Mod{\'e}er},
  \bibinfo{author}{P.~E. Petersson},
\newblock \bibinfo{title}{Analysis of crack formation and crack growth in
  concrete by means of fracture mechanics and finite elements},
\newblock \bibinfo{journal}{Cement and Concrete Research} \bibinfo{volume}{6}
  (\bibinfo{year}{1976}) \bibinfo{pages}{773--781}.
\bibitem[{Needleman(1987)}]{needleman1987}
\bibinfo{author}{A.~Needleman},
\newblock \bibinfo{title}{A {{Continuum Model}} for {{Void Nucleation}} by
  {{Inclusion Debonding}}},
\newblock \bibinfo{journal}{Journal of Applied Mechanics} \bibinfo{volume}{54}
  (\bibinfo{year}{1987}) \bibinfo{pages}{525--531}.
\bibitem[{Tvergaard(1990)}]{tvergaard1990}
\bibinfo{author}{V.~Tvergaard},
\newblock \bibinfo{title}{Effect of fibre debonding in a whisker-reinforced
  metal},
\newblock \bibinfo{journal}{Materials Science and Engineering: A}
  \bibinfo{volume}{125} (\bibinfo{year}{1990}) \bibinfo{pages}{203--213}.
\bibitem[{Paggi and Reinoso(2017)}]{paggi2017}
\bibinfo{author}{M.~Paggi}, \bibinfo{author}{J.~Reinoso},
\newblock \bibinfo{title}{Revisiting the problem of a crack impinging on an
  interface:{{A}} modeling framework for the interaction between the phase
  field approach for brittle fracture and the interface cohesive zone model},
\newblock \bibinfo{journal}{Computer Methods in Applied Mechanics and
  Engineering} \bibinfo{volume}{321} (\bibinfo{year}{2017})
  \bibinfo{pages}{145--172}.
\bibitem[{Pastor et~al.(1991)Pastor, Peraire, and Zienkiewicz}]{pastor1991}
\bibinfo{author}{M.~Pastor}, \bibinfo{author}{J.~Peraire},
  \bibinfo{author}{O.~C. Zienkiewicz},
\newblock \bibinfo{title}{Adaptive remeshing for shear band localization
  problems},
\newblock \bibinfo{journal}{Archive of Applied Mechanics} \bibinfo{volume}{61}
  (\bibinfo{year}{1991}) \bibinfo{pages}{30--39}.
\bibitem[{Steinmann and Willam(1992)}]{steinmann1992}
\bibinfo{author}{P.~Steinmann}, \bibinfo{author}{K.~Willam},
\newblock \bibinfo{title}{Adaptive {{Techniques}} for {{Localization
  Analysis}}},
\newblock \bibinfo{journal}{undefined}  (\bibinfo{year}{1992}).
\bibitem[{Belytschko and Tabbara(1993)}]{belytschko1993}
\bibinfo{author}{T.~Belytschko}, \bibinfo{author}{M.~Tabbara},
\newblock \bibinfo{title}{H-{{Adaptive}} finite element methods for dynamic
  problems, with emphasis on localization},
\newblock \bibinfo{journal}{International Journal for Numerical Methods in
  Engineering} \bibinfo{volume}{36} (\bibinfo{year}{1993})
  \bibinfo{pages}{4245--4265}.
\bibitem[{Peric et~al.(1994)Peric, Yu, and Owen}]{peric1994}
\bibinfo{author}{D.~Peric}, \bibinfo{author}{J.~Yu}, \bibinfo{author}{D.~Owen},
\newblock \bibinfo{title}{On error estimates and adaptivity in elastoplastic
  solids: {{Applications}} to the numerical simulation of strain localization
  in classical and {{Cosserat}} continua},
\newblock \bibinfo{journal}{International Journal for Numerical Methods in
  Engineering} \bibinfo{volume}{37} (\bibinfo{year}{1994})
  \bibinfo{pages}{1351--1379}.
\bibitem[{Zienkiewicz et~al.(1995{\natexlab{a}})Zienkiewicz, Pastor, and
  Huang}]{zienkiewicz1995}
\bibinfo{author}{O.~C. Zienkiewicz}, \bibinfo{author}{M.~Pastor},
  \bibinfo{author}{M.~Huang},
\newblock \bibinfo{title}{Softening, localisation and adaptive remeshing.
  {{Capture}} of discontinuous solutions},
\newblock \bibinfo{journal}{Computational Mechanics} \bibinfo{volume}{17}
  (\bibinfo{year}{1995}{\natexlab{a}}) \bibinfo{pages}{98--106}.
\bibitem[{Zienkiewicz et~al.(1995{\natexlab{b}})Zienkiewicz, Huang, and
  Pastor}]{zienkiewicz1995a}
\bibinfo{author}{O.~C. Zienkiewicz}, \bibinfo{author}{M.~Huang},
  \bibinfo{author}{M.~Pastor},
\newblock \bibinfo{title}{Localization problems in plasticity using finite
  elements with adaptive remeshing},
\newblock \bibinfo{journal}{International Journal for Numerical and Analytical
  Methods in Geomechanics} \bibinfo{volume}{19}
  (\bibinfo{year}{1995}{\natexlab{b}}) \bibinfo{pages}{127--148}.
\bibitem[{Ferreira et~al.(2022)Ferreira, Andrade~Pires, and
  Bessa}]{ferreira2022}
\bibinfo{author}{B.~P. Ferreira}, \bibinfo{author}{F.~M. Andrade~Pires},
  \bibinfo{author}{M.~A. Bessa},
\newblock \bibinfo{title}{Fast {{Homogenization}} through {{Clustering}}-based
  {{Reduced Order Modeling}}},
\newblock in: \bibinfo{booktitle}{Fundamentals of {{Multiscale Modeling}} of
  {{Structural Materials}}}, \bibinfo{publisher}{{Elsevier (in press)}},
  \bibinfo{year}{2022}.
\bibitem[{Pedregosa et~al.(2012)Pedregosa, Varoquaux, Gramfort, Michel,
  Thirion, Grisel, Blondel, Prettenhofer, Weiss, Dubourg, Vanderplas, Passos,
  Cournapeau, Brucher, Perrot, Duchesnay, and Louppe}]{pedregosa2012}
\bibinfo{author}{F.~Pedregosa}, \bibinfo{author}{G.~Varoquaux},
  \bibinfo{author}{A.~Gramfort}, \bibinfo{author}{V.~Michel},
  \bibinfo{author}{B.~Thirion}, \bibinfo{author}{O.~Grisel},
  \bibinfo{author}{M.~Blondel}, \bibinfo{author}{P.~Prettenhofer},
  \bibinfo{author}{R.~Weiss}, \bibinfo{author}{V.~Dubourg},
  \bibinfo{author}{J.~Vanderplas}, \bibinfo{author}{A.~Passos},
  \bibinfo{author}{D.~Cournapeau}, \bibinfo{author}{M.~Brucher},
  \bibinfo{author}{M.~Perrot}, \bibinfo{author}{E.~Duchesnay},
  \bibinfo{author}{G.~Louppe},
\newblock \bibinfo{title}{Scikit-learn: {{Machine Learning}} in {{Python}}},
\newblock \bibinfo{journal}{Journal of Machine Learning Research}
  \bibinfo{volume}{12} (\bibinfo{year}{2012}).
\bibitem[{Zhang et~al.(2019)Zhang, Tang, Yu, Zhu, and Liu}]{zhang2019a}
\bibinfo{author}{L.~Zhang}, \bibinfo{author}{S.~Tang}, \bibinfo{author}{C.~Yu},
  \bibinfo{author}{X.~Zhu}, \bibinfo{author}{W.~K. Liu},
\newblock \bibinfo{title}{Fast calculation of interaction tensors in
  clustering-based homogenization},
\newblock \bibinfo{journal}{Computational Mechanics} \bibinfo{volume}{64}
  (\bibinfo{year}{2019}) \bibinfo{pages}{351--364}.
\bibitem[{Zienkiewicz and Zhu(1987)}]{zienkiewicz1987}
\bibinfo{author}{O.~C. Zienkiewicz}, \bibinfo{author}{J.~Z. Zhu},
\newblock \bibinfo{title}{A simple error estimator and adaptive procedure for
  practical engineerng analysis},
\newblock \bibinfo{journal}{International Journal for Numerical Methods in
  Engineering} \bibinfo{volume}{24} (\bibinfo{year}{1987})
  \bibinfo{pages}{337--357}.
\bibitem[{Zienkiewicz and Zhu(1992)}]{zienkiewicz1992}
\bibinfo{author}{O.~C. Zienkiewicz}, \bibinfo{author}{J.~Z. Zhu},
\newblock \bibinfo{title}{The superconvergent patch recovery and a posteriori
  error estimates. {{Part}} 1: {{The}} recovery technique},
\newblock \bibinfo{journal}{International Journal for Numerical Methods in
  Engineering} \bibinfo{volume}{33} (\bibinfo{year}{1992})
  \bibinfo{pages}{1331--1364}.
\bibitem[{Boroomand and Zienkiewicz(1997)}]{boroomand1997}
\bibinfo{author}{B.~Boroomand}, \bibinfo{author}{O.~C. Zienkiewicz},
\newblock \bibinfo{title}{Recovery by {{Equilibrium}} in {{Patches}} (rep)},
\newblock \bibinfo{journal}{International Journal for Numerical Methods in
  Engineering} \bibinfo{volume}{40} (\bibinfo{year}{1997})
  \bibinfo{pages}{137--164}.
\bibitem[{Kelly et~al.(1983)Kelly, Gago, Zienkiewicz, and Babuska}]{kelly1983}
\bibinfo{author}{D.~W. Kelly}, \bibinfo{author}{J.~P. D. S.~R. Gago},
  \bibinfo{author}{O.~C. Zienkiewicz}, \bibinfo{author}{I.~Babuska},
\newblock \bibinfo{title}{A posteriori error analysis and adaptive processes in
  the finite element method: {{Part I}}\textemdash error analysis},
\newblock \bibinfo{journal}{International Journal for Numerical Methods in
  Engineering} \bibinfo{volume}{19} (\bibinfo{year}{1983})
  \bibinfo{pages}{1593--1619}.
\bibitem[{Verf{\"u}rth(1994)}]{verfurth1994}
\bibinfo{author}{R.~Verf{\"u}rth},
\newblock \bibinfo{title}{A posteriori error estimation and adaptive
  mesh-refinement techniques},
\newblock \bibinfo{journal}{Journal of Computational and Applied Mathematics}
  \bibinfo{volume}{50} (\bibinfo{year}{1994}) \bibinfo{pages}{67--83}.
\bibitem[{Bank and Weiser(1985)}]{bank1985}
\bibinfo{author}{R.~E. Bank}, \bibinfo{author}{A.~Weiser},
\newblock \bibinfo{title}{Some {{A Posteriori Error Estimators}} for {{Elliptic
  Partial Differential Equations}}},
\newblock \bibinfo{journal}{Mathematics of Computation} \bibinfo{volume}{44}
  (\bibinfo{year}{1985}) \bibinfo{pages}{283--301}.
\bibitem[{Ainsworth and Oden(1993)}]{ainsworth1993}
\bibinfo{author}{M.~Ainsworth}, \bibinfo{author}{J.~T. Oden},
\newblock \bibinfo{title}{A unified approach to a posteriori error estimation
  using element residual methods},
\newblock \bibinfo{journal}{Numerische Mathematik} \bibinfo{volume}{65}
  (\bibinfo{year}{1993}) \bibinfo{pages}{23--50}.
\bibitem[{Ladev{\`e}ze(1985)}]{ladeveze1985}
\bibinfo{author}{P.~Ladev{\`e}ze},
\newblock \bibinfo{title}{Sur une famille d'algorithmes en m\'ecanique des
  structures},
\newblock \bibinfo{journal}{undefined}  (\bibinfo{year}{1985}).
\bibitem[{Ladev{\`e}ze and Germain(1989)}]{ladeveze1989}
\bibinfo{author}{P.~Ladev{\`e}ze}, \bibinfo{author}{P.~Germain},
\newblock \bibinfo{title}{La m\'ethode \`a grand incr\'ement de temps pour
  l'analyse de structures \`a comportement non lin\'eaire d\'ecrit par
  variables internes},
\newblock \bibinfo{journal}{undefined}  (\bibinfo{year}{1989}).
\bibitem[{Gallimard et~al.(1996)Gallimard, Ladev{\`e}ze, and
  Pelle}]{gallimard1996}
\bibinfo{author}{L.~Gallimard}, \bibinfo{author}{P.~Ladev{\`e}ze},
  \bibinfo{author}{J.~Pelle},
\newblock \bibinfo{title}{Error estimation and adaptivity in elastoplasticity},
\newblock \bibinfo{journal}{International Journal for Numerical Methods in
  Engineering} \bibinfo{volume}{39} (\bibinfo{year}{1996})
  \bibinfo{pages}{189--217}.
\bibitem[{Eriksson and Johnson(1988)}]{eriksson1988}
\bibinfo{author}{K.~Eriksson}, \bibinfo{author}{C.~Johnson},
\newblock \bibinfo{title}{An {{Adaptive Finite Element Method}} for {{Linear
  Elliptic Problems}}},
\newblock \bibinfo{journal}{Mathematics of Computation} \bibinfo{volume}{50}
  (\bibinfo{year}{1988}) \bibinfo{pages}{361--383}.
\bibitem[{Johnson and Hansbo(1992)}]{johnson1992}
\bibinfo{author}{C.~Johnson}, \bibinfo{author}{P.~Hansbo},
\newblock \bibinfo{title}{Adaptive finite element methods in computational
  mechanics},
\newblock \bibinfo{journal}{Computer Methods in Applied Mechanics and
  Engineering} \bibinfo{volume}{101} (\bibinfo{year}{1992})
  \bibinfo{pages}{143--181}.
\bibitem[{Rannacher and Suttmeier(1998)}]{rannacher1998}
\bibinfo{author}{R.~Rannacher}, \bibinfo{author}{F.~T. Suttmeier},
\newblock \bibinfo{title}{A {{Posteriori Error Control}} and {{Mesh
  Adaptation}} for {{FE}} models in {{Elasticity}} and
  {{Elasto}}-{{Plasticity}}},
\newblock in: \bibinfo{booktitle}{Advances in Adaptive Computational Methods in
  Mechanics}, number~\bibinfo{number}{47} in \bibinfo{series}{Studies in
  Applied Mechanics}, \bibinfo{publisher}{{Elsevier}},
  \bibinfo{address}{{Amsterdam ; New York}}, \bibinfo{year}{1998}, pp.
  \bibinfo{pages}{275--292}.
\bibitem[{Jin and Wiberg(1990)}]{jin1990}
\bibinfo{author}{H.~Jin}, \bibinfo{author}{N.-E. Wiberg},
\newblock \bibinfo{title}{Two-dimensional mesh generation, adaptive remeshing
  and refinement},
\newblock \bibinfo{journal}{International Journal for Numerical Methods in
  Engineering} \bibinfo{volume}{29} (\bibinfo{year}{1990})
  \bibinfo{pages}{1501--1526}.
\bibitem[{Peri{\'c} et~al.(1998)Peri{\'c}, Dutko, and Owen}]{peric1998a}
\bibinfo{author}{D.~Peri{\'c}}, \bibinfo{author}{M.~Dutko},
  \bibinfo{author}{D.~R.~J. Owen},
\newblock \bibinfo{title}{Aspects of adaptive strategies for large deformation
  problems at finite inelastic strains},
\newblock in: \bibinfo{booktitle}{Aspects of Adaptive Strategies for Large
  Deformation Problems at Finite Inelastic Strains},
  number~\bibinfo{number}{47} in \bibinfo{series}{Studies in Applied
  Mechanics}, \bibinfo{publisher}{{Elsevier}}, \bibinfo{address}{{Amsterdam ;
  New York}}, \bibinfo{year}{1998}, pp. \bibinfo{pages}{349--363}.
\bibitem[{Lee and Bathe(1994)}]{lee1994}
\bibinfo{author}{N.-S. Lee}, \bibinfo{author}{K.-J. Bathe},
\newblock \bibinfo{title}{Error indicators and adaptive remeshing in large
  deformation finite element analysis},
\newblock \bibinfo{journal}{Finite Elements in Analysis and Design}
  \bibinfo{volume}{16} (\bibinfo{year}{1994}) \bibinfo{pages}{99--139}.
\bibitem[{Peri{\'c} et~al.(1999)Peri{\'c}, Vaz, and Owen}]{peric1999}
\bibinfo{author}{D.~Peri{\'c}}, \bibinfo{author}{M.~Vaz},
  \bibinfo{author}{D.~R.~J. Owen},
\newblock \bibinfo{title}{On adaptive strategies for large deformations of
  elasto-plastic solids at finite strains: Computational issues and industrial
  applications},
\newblock \bibinfo{journal}{Computer Methods in Applied Mechanics and
  Engineering} \bibinfo{volume}{176} (\bibinfo{year}{1999})
  \bibinfo{pages}{279--312}.
\bibitem[{Zienkiewicz and Wu(1994)}]{zienkiewicz1994}
\bibinfo{author}{O.~C. Zienkiewicz}, \bibinfo{author}{J.~Wu},
\newblock \bibinfo{title}{Automatic directional refinement in adaptive analysis
  of compressible flows},
\newblock \bibinfo{journal}{International Journal for Numerical Methods in
  Engineering} \bibinfo{volume}{37} (\bibinfo{year}{1994})
  \bibinfo{pages}{2189--2210}.
\bibitem[{Yu et~al.(1991)Yu, Peric, and Owen}]{yu1991}
\bibinfo{author}{J.~Yu}, \bibinfo{author}{D.~Peric}, \bibinfo{author}{D.~R.~J.
  Owen},
\newblock \bibinfo{title}{An {{Assessment}} of the {{Cosserat Continuum}}
  through the {{Finite Element Simulation}} of a {{Strain Localisation
  Problem}}},
\newblock in: \bibinfo{editor}{E.~O{\~n}ate}, \bibinfo{editor}{J.~Periaux},
  \bibinfo{editor}{A.~Samuelsson} (Eds.), \bibinfo{booktitle}{The Finite
  Element Method in the 1990's: {{A Book Dedicated}} to {{O}}.{{C}}.
  {{Zienkiewicz}}}, \bibinfo{publisher}{{Springer}}, \bibinfo{address}{{Berlin,
  Heidelberg}}, \bibinfo{year}{1991}, pp. \bibinfo{pages}{321--332}.
  \DOIprefix\doi{10.1007/978-3-662-10326-5_33}.
\bibitem[{Pijaudier-Cabot et~al.(1995)Pijaudier-Cabot, Bod{\'e}, and
  Huerta}]{pijaudier-cabot1995}
\bibinfo{author}{G.~Pijaudier-Cabot}, \bibinfo{author}{L.~Bod{\'e}},
  \bibinfo{author}{A.~Huerta},
\newblock \bibinfo{title}{Arbitrary {{Lagrangian}}\textendash{{Eulerian}}
  finite element analysis of strain localization in transient problems},
\newblock \bibinfo{journal}{International Journal for Numerical Methods in
  Engineering} \bibinfo{volume}{38} (\bibinfo{year}{1995})
  \bibinfo{pages}{4171--4191}.
\bibitem[{Stein et~al.(1998)Stein, Barthold, Ohnimus, and Schmidt}]{stein1998}
\bibinfo{author}{E.~Stein}, \bibinfo{author}{F.~J. Barthold},
  \bibinfo{author}{S.~Ohnimus}, \bibinfo{author}{M.~Schmidt},
\newblock \bibinfo{title}{Adaptive finite elements in elastoplasticity with
  mechanical error indicators and {{Neumann}}-type estimators},
\newblock in: \bibinfo{booktitle}{Advances in Adaptive Computational Methods in
  Mechanics}, number~\bibinfo{number}{47} in \bibinfo{series}{Studies in
  Applied Mechanics}, \bibinfo{publisher}{{Elsevier}},
  \bibinfo{address}{{Amsterdam ; New York}}, \bibinfo{year}{1998}, pp.
  \bibinfo{pages}{81--99}.
\bibitem[{Babu{\v s}ka(1988)}]{babuska1988}
\bibinfo{author}{I.~Babu{\v s}ka},
\newblock \bibinfo{title}{The p and h-p {{Versions}} of the {{Finite Element
  Method}}: {{The State}} of the {{Art}}},
\newblock in: \bibinfo{editor}{D.~L. Dwoyer}, \bibinfo{editor}{M.~Y. Hussaini},
  \bibinfo{editor}{R.~G. Voigt} (Eds.), \bibinfo{booktitle}{Finite
  {{Elements}}}, {{ICASE}}/{{NASA LaRC Series}},
  \bibinfo{publisher}{{Springer}}, \bibinfo{address}{{New York, NY}},
  \bibinfo{year}{1988}, pp. \bibinfo{pages}{199--239}.
  \DOIprefix\doi{10.1007/978-1-4612-3786-0_10}.
\bibitem[{Guo and Babu{\v s}ka(1986{\natexlab{a}})}]{guo1986}
\bibinfo{author}{B.~Guo}, \bibinfo{author}{I.~Babu{\v s}ka},
\newblock \bibinfo{title}{The h-p version of the finite element method},
\newblock \bibinfo{journal}{Computational Mechanics} \bibinfo{volume}{1}
  (\bibinfo{year}{1986}{\natexlab{a}}) \bibinfo{pages}{21--41}.
\bibitem[{Guo and Babu{\v s}ka(1986{\natexlab{b}})}]{guo1986a}
\bibinfo{author}{B.~Guo}, \bibinfo{author}{I.~Babu{\v s}ka},
\newblock \bibinfo{title}{The h-p version of the finite element method},
\newblock \bibinfo{journal}{Computational Mechanics} \bibinfo{volume}{1}
  (\bibinfo{year}{1986}{\natexlab{b}}) \bibinfo{pages}{203--220}.
\bibitem[{Demkowicz et~al.(1989)Demkowicz, Oden, Rachowicz, and
  Hardy}]{demkowicz1989}
\bibinfo{author}{L.~Demkowicz}, \bibinfo{author}{J.~T. Oden},
  \bibinfo{author}{W.~Rachowicz}, \bibinfo{author}{O.~Hardy},
\newblock \bibinfo{title}{Toward a universal h-p adaptive finite element
  strategy, part 1. {{Constrained}} approximation and data structure},
\newblock \bibinfo{journal}{Computer Methods in Applied Mechanics and
  Engineering} \bibinfo{volume}{77} (\bibinfo{year}{1989})
  \bibinfo{pages}{79--112}.
\bibitem[{Rachowicz et~al.(1989)Rachowicz, Oden, and Demkowicz}]{rachowicz1989}
\bibinfo{author}{W.~Rachowicz}, \bibinfo{author}{J.~T. Oden},
  \bibinfo{author}{L.~Demkowicz},
\newblock \bibinfo{title}{Toward a universal h-p adaptive finite element
  strategy part 3. design of h-p meshes},
\newblock \bibinfo{journal}{Computer Methods in Applied Mechanics and
  Engineering} \bibinfo{volume}{77} (\bibinfo{year}{1989})
  \bibinfo{pages}{181--212}.
\bibitem[{Zienkiewicz et~al.(1989)Zienkiewicz, Zhu, and Gong}]{zienkiewicz1989}
\bibinfo{author}{O.~C. Zienkiewicz}, \bibinfo{author}{J.~Z. Zhu},
  \bibinfo{author}{N.~G. Gong},
\newblock \bibinfo{title}{Effective and practical h\textendash p-version
  adaptive analysis procedures for the finite element method},
\newblock \bibinfo{journal}{International Journal for Numerical Methods in
  Engineering} \bibinfo{volume}{28} (\bibinfo{year}{1989})
  \bibinfo{pages}{879--891}.
\bibitem[{Nochetto et~al.(2009)Nochetto, Siebert, and Veeser}]{nochetto2009a}
\bibinfo{author}{R.~H. Nochetto}, \bibinfo{author}{K.~G. Siebert},
  \bibinfo{author}{A.~Veeser},
\newblock \bibinfo{title}{Theory of adaptive finite element methods: {{An}}
  introduction},
\newblock in: \bibinfo{editor}{R.~DeVore}, \bibinfo{editor}{A.~Kunoth} (Eds.),
  \bibinfo{booktitle}{Multiscale, {{Nonlinear}} and {{Adaptive
  Approximation}}}, \bibinfo{publisher}{{Springer}}, \bibinfo{address}{{Berlin,
  Heidelberg}}, \bibinfo{year}{2009}, pp. \bibinfo{pages}{409--542}.
  \DOIprefix\doi{10.1007/978-3-642-03413-8_12}.

\end{thebibliography}

\end{document}